\documentclass{article}
\usepackage[a4paper,top=2.6cm,bottom=2.6cm,left=2.6cm,right=2.6cm]{geometry}
\usepackage{authblk}
\usepackage{paralist}
\usepackage{graphics}
\usepackage{epsfig} 
\usepackage{graphicx}  
\usepackage{epstopdf}
\usepackage[colorlinks=true]{hyperref}
\hypersetup{urlcolor=blue, linkcolor=red, citecolor=red}

\usepackage{amsmath}
\usepackage{mathtools} 
\mathtoolsset{showonlyrefs}
\usepackage{amsfonts}
\usepackage{amssymb}
\usepackage{amsthm}
\usepackage{mathrsfs}
\usepackage{stackrel,graphicx}
\usepackage{enumitem}
\usepackage{cases}
\usepackage{cancel}
\usepackage{booktabs}
\usepackage{braket}
\usepackage[strict]{changepage}
\usepackage{subfigure}
\usepackage[small,justification=justified]{caption}
\usepackage{titlesec} 
\usepackage{multicol}
\usepackage{hyperref}
\usepackage{pgfplots}
\usepackage[cbgreek]{textgreek}
\usepackage{bbm}
\usepackage{amssymb,amsmath,amsthm}
\usepackage{multirow}
\usepackage[table]{xcolor}
\usepackage[right]{lineno}
\usepackage{verbatim}

\usepackage[sorting=none, style=chem-angew, articletitle=true, maxnames=2]{biblatex}

\makeatletter
\def\blx@maxline{77}
\makeatother
\bibliography{bibl_analysis}
\DeclareRedundantLanguages{english,german,french}{english,german,ngerman,french}

\newtheorem{theorem}{Theorem}[section]
\newtheorem{cor}{Corollary}[section]

\newtheorem{lem}[theorem]{Lemma}

\newtheorem{thm}{Theorem}[section]
\newtheorem*{rem*}{Remark}
\theoremstyle{definition}
\newtheorem{define}{Definition}[section]
\newtheorem{rem}{Remark}[section]

\newtheorem{assume}{Assumption}[section]

\providecommand{\keywords}[1]{\textbf{\textit{Keywords---}} #1}

\newcommand{\brho}{\boldsymbol{\varrho}}
\newcommand{\vrhoone}{{^1\!\varrho}}
\newcommand{\vrhotwo}{{^2\!\varrho}}

\newcommand{\cgam}{\text{\rm \textgamma}}
\newcommand{\clam}{\text{\rm \textlambda}}

\newenvironment{sistem}
{\left\lbrace\begin{array}{@{}l@{}}}
{\end{array}\right.}

\let\oldproofname=\proofname
\renewcommand{\proofname}{\rm\bf{\oldproofname}}

\newcommand{\grad}{\nabla}

\newcommand{\R}{\mathbb{R}}
\newcommand{\N}{\mathbb{N}}

\newcommand{\xbf}{{\bf x}}

\renewcommand*{\sup}[1]{\underset{#1}{\text{{\rm sup} }}}
\renewcommand*{\inf}[1]{\underset{#1}{\text{{\rm inf} }}}
\renewcommand*{\liminf}[1]{\underset{#1}{\text{{\rm liminf} }}}

\newcommand{\D}{\mathcal{D}}
\newcommand{\dom}{\mathfrak{D}}

\newcommand{\Om}{\Omega}
\newcommand{\om}{\omega}

\renewcommand{\L}{\mathcal{L}}

\newcommand{\Ra}{\Rightarrow}

\newcommand{\vsig}{\varsigma}
\newcommand{\vrho}{\varrho}

\newcommand{\G}{\mathcal{G}}

\newcommand{\CLTi}[2]{C([t_{i-1},t_i];L_2(\dom)}
\newcommand{\CHTi}[2]{C([t_{i-1},t_i];H^1(\dom)}

\newcommand{\Fbs}{\mathcal{G}^{\beta,\sigma}}

\newcommand{\X}{\mathbf{X}}
\newcommand{\loc}{{\rm{loc}}}
\newcommand{\Zb}{\mathbf{Z}}
\newcommand{\cf}{{\rm{cf}}}

\newcommand{\const}{{\mbox{\large{\rm k}}}}

\newcommand{\lbd}{{\rm m}}
\newcommand{\ubd}{{\rm M}}

\newcommand{\indc}{\mathbbm{1}}

\newcommand{\E}{\mathbb{E}}
\renewcommand{\P}{\mathbb{P}}

\newcommand{\FT}{\mathfrak{F}}
\newcommand{\IFT}{\mathfrak{F}^{-1}}

\newcommand{\SW}{\mathcal{S}}

\newcommand{\TD}{\mathcal{S}'}

\newcommand{\Fil}{\mathcal{F}}
\newcommand{\cnt}{\mathcal{C}}
\newcommand{\B}{\mathcal{B}}
\newcommand{\Y}{\mathcal{Y}}
\newcommand{\W}{\mathcal{W}}
\newcommand{\V}{\mathcal{V}}
\newcommand{\Z}{\mathcal{Z}}

\newcommand{\bU}{\mathbf{U}}

\newcommand{\K}{\mathcal{K}}
\newcommand{\wsp}{\W^s_p}
\newcommand{\wop}{\W^1_p}
\newcommand{\wtp}{\W^2_p}
\newcommand{\wrp}{\W^r_p}
\newcommand{\zsp}{\Z^s_p}
\newcommand{\zop}{\Z^1_p}
\newcommand{\ztp}{\Z^2_p}
\newcommand{\zrp}{\Z^r_p}
\newcommand{\T}{\mathcal{T}}
\newcommand{\Tst}{\T_{s,t}}
\newcommand{\Tzt}{\T_{0,t}}
\newcommand{\Tzs}{\T_{0,s}}

\newcommand{\Ttaut}{\T_{\tau,t}}
\newcommand{\Tstau}{\T_{s,\tau}}

\pgfplotsset{compat=1.6}

\newcommand{\bmat}[1]{ 
	\begin{bmatrix} #1
\end{bmatrix}}

\begin{document}

\title{Existence and uniqueness of solutions for a random PDE-ODE model describing glioma invasion and radiotherapy}
\author[1]{Martina Conte}
\author[2]{Sandesh Athni Hiremath}
\author[2]{Christina Surulescu}
\affil[1]{\centerline{\small Department of Mathematical Sciences "G.L. Lagrange", Politecnico di Torino} \newline \centerline{\small  Corso Duca degli Abruzzi, 24 - 10129 Torino, Italy}}
\affil[2]{\centerline {\small Department of Mathematics, RPTU Kaiserslautern-Landau} \newline \centerline{\small Gottlieb-Daimler-Str. 48 - 67663 Kaiserslautern, Germany}}

\setcounter{Maxaffil}{0}
\renewcommand\Affilfont{\itshape\small}
\maketitle

\begin{abstract}
We study a stochastic model for the growth and invasion of glioma cells in the anisotropic brain tissue, also incorporating effects of radiotherapeutic interventions and associated uncertainties. The PDE--ODE model couples glioma and healthy tissue dynamics in a nonlinear manner and explicitly includes treatment-related stochasticity through random depletion terms. The macroscopic system has been obtained by an upscaling procedure from a microscopic description in~\cite{conte2026randomness}. We prove the existence and uniqueness of local strong solutions to the corresponding stochastic system. The analysis addresses challenges arising from the nonlinear coupling, anisotropic diffusion, taxis mechanisms, and stochastic treatment effects, providing a first rigorous step toward the mathematical understanding of stochastic continuum models for glioma invasion and therapy response.
\end{abstract}

\keywords{Glioma treatment modeling, random reaction-diffusion-taxis equations, existence and uniqueness}

\section{Introduction}
Treatment of glioma is a challenging task, primarily due to its highly aggressive, diffusive infiltration of brain tissue, with fractal tumor margins that cannot be assessed by medical imaging. This makes exhaustive resection or radiation impossible and calls for suitable models to help evaluate the tumor dynamics and enhance therapeutic strategies~\cite{rockne2009mathematical,swanson2011quantifying,engwer2016effective,hunt2017multiscale,conte2021modeling,conte2023mathematical}.

Glioma progression is a result of complex interactions between malignant cells and their microenvironment~\cite{wrensch2002epidemiology,ostrom2014epidemiology,wick2018treatment,berens1999those}. A promising route to effectively connecting microscopic, subcellular-level dynamics with mesoscopic and macroscopic glioma behavior is provided by multiscale frameworks based on kinetic descriptions~\cite{bellomo2022towards,Bell-KTAP}. Such approaches yield macroscopic diffusion-advection-reaction partial differential equations (PDEs) by asymptotic limits and can naturally accommodate tissue anisotropy and mechanistic interactions~\cite{painter2013mathematical,Swan2017,conte2023mathematical,engwer2015glioma,hunt2017multiscale,engwer2016effective}. However, such and related deterministic settings~\cite{conte2021modeling,rockne2009mathematical,swanson2000quantitative,swanson2011quantifying,harpold2007evolution} cannot account for the large intertumor variability observed in clinical practice. Moreover, radiotherapy itself is quite complex, and its efficacy is affected by substantial uncertainties stemming from biological heterogeneity (variable radiosensitivity, oxygen and nutrient fluctuations), as well as from practical delivery issues (patient positioning, imaging and segmentation limits, missed or delayed fractions)~\cite{gondi2021radiotherapy,nieder1999hyperfractionated,beckham2024hypofractionated,trone2020survival,junttila2013influence,segedin2016uncertainties,winey2014geometric}. Stochastic PDE (SPDE) formulations involving random perturbations in motility, tumor microenvironment, and tumor–tissue interactions~\cite{hiremath2018coupled,hiremath2016stochastic,Katsaounis2023,Lima2014,Essarrout2022} are potentially better suited for capturing random effects of radiotherapy on glioma progression. Motivated by the clinical importance of treatment-related variability and by the need for models that propagate randomness from the microscale to clinically relevant macroscopic observables, in a recent work~\cite{conte2026randomness} we proposed a novel multiscale stochastic model which incorporates stochasticity at the mesoscopic level and propagates it up to the macroscopic tumor scale, resulting in a random PDE-ODE system with haptotaxis. 

In this work, we focus on the well‑posedness of the coupled tumor–tissue random PDE system in~\cite{conte2026randomness} featuring myopic diffusion and haptotaxis along with terms accounting for random effects of radiotherapy, under biologically motivated structural assumptions. The mathematical challenges addressed here stem from the strong nonlinear coupling between tumor cell density and normal tissue, the presence of anisotropic, spatially varying diffusion tensors originating from patient-specific diffusion tensor imaging (DTI) data, and stochasticity introduced by the therapy terms. These features place our problem at the interface between deterministic reaction–diffusion–taxis theory and the theory of SPDEs in infinite-dimensional spaces. To tackle these challenges, we build on a variety of analytical tools including maximal regularity and non-autonomous evolution operator methods for quasilinear parabolic problems (see e.g.~\cite{Yagi09,PAZY83,ENGEL99}), and standard references for the stochastic analysis of infinite-dimensional systems~\cite{DaPrato1992,PrevotRockner2007,oksendal2013stochastic}.

The rest of the manuscript is organized as follows. In Section~\ref{modeling} we recall the macroscopic random PDE model derived in~\cite{conte2026randomness}, restating the important structural properties of the coefficients and the stochastic radiotherapy terms, required for the subsequent analysis. Section~\ref{sec:analysis} contains the main analytical results: local existence and uniqueness of strong solutions for the stochastic coupled tumor–tissue system, along with a priori bounds and regularity properties needed for the semigroup approach. We provide clear assumptions that capture the biological parameter regimes used in the modeling study and indicate how the stochastic radiotherapy perturbations are handled mathematically. Finally, Section~\ref{sec:discussion} summarizes our findings and provides some comments on related works. 

\section{Model description and motivation}\label{modeling}
In this section, we recall the main characteristics and structural properties of the multiscale random PDE model presented in~\cite{conte2026randomness}, which describes the progression of glioma cells in the brain under different radiotherapy regimes. The model integrates three levels of description, namely microscopic, mesoscopic, and macroscopic, within the framework of the kinetic theory of active particles (see e.g.,~\cite{bellomo2022towards} and references therein). Its primary objective is to explicitly couple the dynamics of glioma cells and healthy tissue and to analyze their responses to various therapeutic protocols while accounting for randomness in treatment outcomes. Incorporating stochastic effects allows for randomness-aware predictions, which facilitate the assessment of the robustness of treatment protocols and may support the optimization of radiotherapy under realistic sources of variability.

The model describes the migration of glioma cells along the fibrous structure of brain tissue, with a tendency to avoid regions with intraspecific overcrowding. Cancer cell migration is assumed to be mediated by transmembrane receptors, which play a key role in cell–cell and cell–tissue adhesion mechanisms. Cell proliferation is modeled as dependent on the availability of surrounding healthy tissue, which supports tumor growth and provides essential resources. Finally, radiotherapy-induced damage to both tumor and normal tissue is modeled using the classical linear–quadratic (LQ) formulation, modulated by an Ornstein–Uhlenbeck (OU) process to account for uncertainty and variability in treatment effects. The model is formulated by first introducing a microscopic description of the interactions between cells and tissue, which leads to the corresponding kinetic transport equations for mesoscopic glioma cell density. By performing a parabolic scaling limit, a macroscopic random reaction–advection–diffusion equation governing tumor cell evolution is obtained and coupled with a random ordinary differential equation (ODE) for the healthy tissue (see~\cite{conte2026randomness} for a detailed derivation of the model).

In the following, we directly focus on the resulting random macroscopic system for the densities $M(t,x)$ and $Q(t,x)$ of tumor cells and normal  tissue, respectively:

\begin{align}
\begin{sistem}
\begin{aligned}
  \partial_tM &+\nabla_x\cdot \left[\dfrac{\lambda_1k^-}{(B(M,Q)+\lambda_0)B(M,Q)^2}\mathbb{D}_M(x)\,\left(k_1^+ S_Q \nabla_x Q-k_2^+ S_M \nabla_x M\right)\,M\right]\\[0.4cm]
  &=\nabla_x \nabla_x:\,\Big(\mathbb{D}_M(x) M\Big) + \Big[\mu_{M}(1 - M)Q-R_M(d_r,\xi_t)\Big]M  \,,
\end{aligned} \\[1.5cm]
\partial_t Q=-d_Q \dfrac{M Q}{1+M}-R_Q(d_r,\xi_t)Q\,,
\end{sistem}
\label{mac_sis}
\end{align}
where $t>0$ is the time variable, while $x\in D\subseteq \mathbb{R}^n$ represents the spatial position. Here, $\mathbb{D}_M(x)$ denotes the effective tumor diffusion tensor, which accounts for anisotropy in the brain tissue structure. The function $B(M,Q)$, derived from the underlying microscopic cell dynamics, is given by
\begin{equation}
    B(M,Q):=Q S_Q (d_r,\brho_Q)k^+_1+M S_M (d_r,\brho_M)k^+_2+k^-\,,
\end{equation}
while the coefficients $S_i(d_r,\brho_i)$ and $R_i(d_r,\xi_t)$, for $i=M,Q$, represent therapy-related functions. Specifically, to incorporate uncertainties that affect the effective delivered dose, such as patient setup errors, positioning variability, organ motion, biological variability in radiosensitivity, and treatment interruptions, we couple the classical LQ-based survival model with an Ornstein–Uhlenbeck process $\xi_t$. This process captures time-dependent fluctuations in treatment efficacy, reflecting both random delivery errors and biological variability. Modeling $\xi_t$ as a mean-reverting Ornstein–Uhlenbeck process is consistent with the literature on stochastic modeling of time-dependent treatment heterogeneities (see e.g.~\cite{oksendal2013stochastic,hlatky1994influence,badri2016optimal}) and further supports the representation of these perturbations as a stochastic process acting directly on the radiotherapy term. Accordingly, the radiotherapy operator is defined as
\begin{equation}
    \label{eq:radiotherapy_term}
R_i(d_r,\xi_t) = \xi_t \sum_{j=1}^{n_{\mathrm f}} \big(1-S_i(d_r,\brho_i) \big) f_\delta(t-t_j),\qquad \text{for } i \in \{M,Q\}\,. 
\end{equation}
Here, similarly to~\cite{conte2023mathematical,hunt2017multiscale}, the temporal administration of radiation is modeled by smooth, compactly supported pulses $f_\delta(t - t_j)$, representing the $j$-th fraction delivered at time $t_j$. The terms $R_i$ describe the direct cytotoxic effect of ionizing radiation on tumor and normal tissue populations and are grounded in well-established radiobiological principles and clinical considerations. The survival functions $S_i(d_r,\brho_i)$, expressed as 

$$
S_i(d_r,\brho_i) = e^{-(\vrhoone_i + \vrhotwo_i) d_r}, \qquad \text{for } i \in \{M,Q\}\,,
$$ 
represent the fraction of cells of the $i$-population surviving a single radiation dose. The linear ($\vrhoone_i$) and quadratic ($\vrhotwo_i$) parameters correspond to single- and double-track lethal events, respectively (see~\cite{fowler1989linear} and the discussion in~\cite{joiner2009}). Finally, ${n_{\mathrm f}}$ denotes the number of fractions into which the total dose $d_r$ is divided. In clinical practice, fractionated schedules, where the total therapeutic dose is delivered in a sequence of smaller fractions, are adopted to widen the therapeutic window and reduce potential treatment-related side effects (see, e.g.,~\cite{o1997response, wheldon1998linear}). This formulation also allows for a spatially dependent prescribed dose, i.e., $d_r=d_r(x)$, thereby enabling the modeling of different dose painting strategies based on heterogeneous dose prescriptions. For completeness, the meaning of the remaining parameters employed in model~\eqref{mac_sis} is reported in Table~\ref{tab_parameter}.

\begin{table} [!h]
\begin{center}
   \begin{tabular}{|l|l|} 
   \toprule  
   \rule{0pt}{1ex}Parameter & Description \\
  \midrule
   \rule{0pt}{2.5ex}$\lambda_0$ & basal turning rate \\
\rule{0pt}{2.5ex}$\lambda_1$ & sensitivity to environmental cues \\
\rule{0pt}{2.5ex}$k^+_1$ & tumor-ECM attachment rate \\
\rule{0pt}{2.5ex}$k^+_2$ & cell-cell attachment rate \\
\rule{0pt}{2.5ex}$k^-$ & detachment rate \\
\rule{0pt}{2.5ex}$\mu_{M}$ & tumor proliferation rate  \\
\rule{0pt}{2.5ex}$d_Q$ & tissue degradation rate (by acidity) \\
\rule{0pt}{2.5ex}$\mu_R$ & expectation of the OU process \\
\rule{0pt}{2.5ex}$\sigma_R$ & diffusivity of the OU process \\
\bottomrule
    \end{tabular}
\end{center}
\caption{{ Model parameters of system~\eqref{mac_sis}}}
 \label{tab_parameter}
\end{table}

\begin{figure}[h!]
       \centering
       \includegraphics[width=\linewidth]{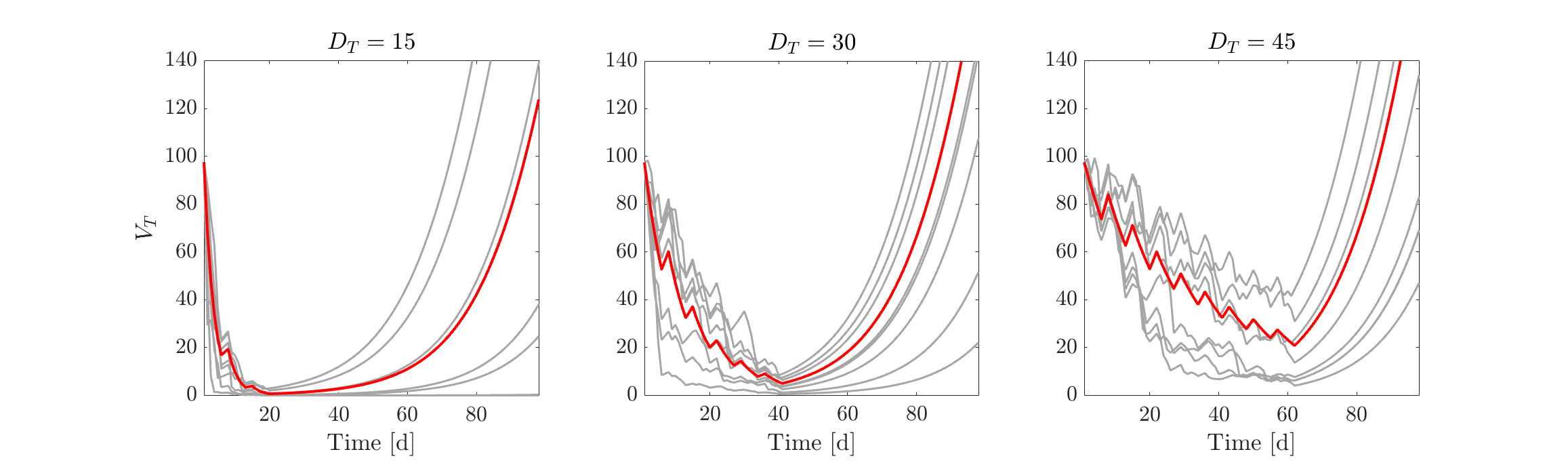}
       \caption{{\bf Deterministic vs stochastic trajectories.} Comparison of simulated tumor volume ($V_T$) evolution over time of the stochastic model (10 different trajectories in gray) and the deterministic model (red line) in the case of $D_T=15,30,45$, as indicated in the panels' title. $D_T$ indicates the number of days during which the treatment is administered.}
       \label{fig:Test1_Vol_SM1}
\end{figure}

\begin{figure}[h!]
       \centering
       \includegraphics[width=\linewidth]{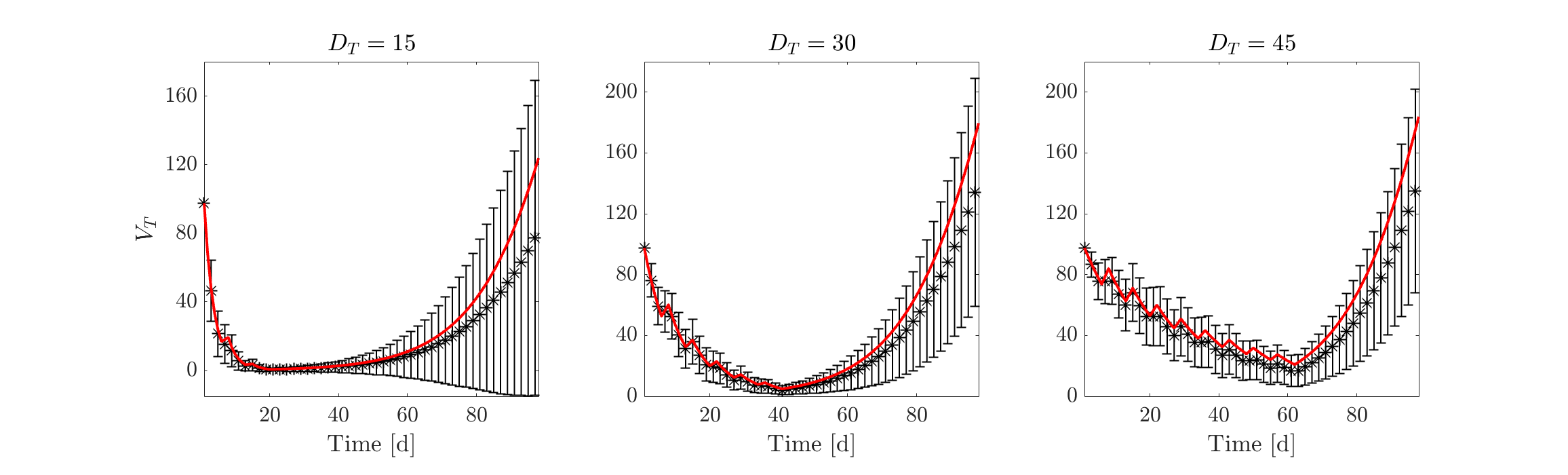}
       \caption{{\bf Mean and standard deviation for the stochastic model.} Comparison of simulated tumor volume ($V_T$) evolution over time of the deterministic model (red line) and the mean (black star markers) and standard deviation (black error bars) obtained from 200 simulation of the stochastic model in the case of the $D_T=15,30,45$, as indicated in the panels' title. $D_T$ indicates the number of days during which the treatment is administered.}
       \label{fig:Test2_STD_VT}
\end{figure}

The proposed stochastic framework enables a range of analyses concerning both the temporal evolution of tumor mass and the effects of different treatment strategies. In particular, we emphasize that incorporating stochasticity allows for the representation of a much richer spectrum of possible treatment outcomes. As it is well known, patient-specific tumor responses to therapy can vary significantly, a variability that deterministic models are often unable to capture. As an example, in Figure~\ref{fig:Test1_Vol_SM1}, we can observe that system~\eqref{mac_sis} is capable of reproducing diverse scenarios, including tumor growth, eradication, and controlled progression (gray trajectories). In contrast, the deterministic model yields a single trajectory (red trajectory), which may be qualitatively similar, but can exhibit substantial quantitative differences with respect to tumor volume evolution obtained in the stochastic case. On average, as shown in Figure~\ref{fig:Test2_STD_VT}, the deterministic model and its stochastic counterpart exhibit similar behavior. However, the stochastic model exhibits substantially larger standard deviations, particularly for shorter, higher-daily-dose therapy schedules ($D_T=15$). Overall, these results highlight the importance and relevance of considering stochastic approaches that incorporate randomness in treatment response, as opposed to relying solely on classical deterministic models. Moreover, this framework enables a deeper analysis of radiotherapy outcomes through key metrics such as Tumor Control Probability (TCP), Normal Tissue Complication Probability (NTCP), and Uncomplicated Tumor Control Probability (UTCP), as well as R-score analysis~\cite{conte2026randomness}. These metrics quantify, respectively, the probability of effective tumor eradication (TCP), the risk of normal tissue complications (NTCP), and the probability of achieving tumor control without severe side effects (UTCP). In addition, the R-score provides a practical criterion for ranking and comparing different treatment protocols.

To facilitate the subsequent existence analysis, we introduce the following equivalent formulation of the macroscopic stochastic PDE system.

\begin{define}[{\it Multiscale model}]\label{def:mult_glioma_sys}
Let $u:=M(t,x)$ and $v:=Q(t,x)$ denote the macroscopic variables. Set $\kappa_1(x):=k_2^+ S_M(d_r(x),\brho_M)$, $\kappa_2(x):=k_1^+ S_Q(d_r(x),\brho_Q)$, and let $\kappa_3:=d_Q$ be a positive constant. Further, let $\dom \subset \mathbb{R}^n$, $n \in \{1,2,3\}$, be an open bounded spatial domain with sufficiently smooth boundary {and $T<\infty $}. Then, the macroscopic system of stochastic PDEs describing tumor and normal tissue dynamics, together with the corresponding initial and boundary conditions, is defined as
{\allowdisplaybreaks
\begin{subequations}
\label{eq:mult_glioma_sys}
\begin{align}
\label{eq:pre_cancer}
&\left.\begin{aligned}
\textbf{Cancer cell}&\textbf{ density model:} \\[0.2cm]
&\partial_t u - \grad \cdot  \left(D_1 u+ D_0 \grad u\right)+\grad \cdot\left(g \mathbb{D}_M \grad v\right) = \left[\mu_M(1-u) v - R_M(\xi_t) \right]u \quad \text{in }(0,T]\times \dom,\\[0.2cm]
&\left(D_1 u + D_0 \grad u\right) \cdot {\nu}\, { = }\, \left(g \,\mathbb{D}_M \grad v\right) \cdot  {\nu }= 0\quad \text{on }(0,T]\times \partial \dom\\[0.2cm] 
& u(0) = u_0 \quad \text{in }\dom\\[0.2cm]
&\text{\rm {with $\nu $ denoting the outer unit normal on $\partial \dom$, and}} \qquad \quad&\\[0.2cm]
&D_0 := \mathbb{D}_M( 1 + \kappa_1 \zeta(u,v) ), \quad {D_1 := \nabla \cdot \mathbb{D}_M}, \quad  g(u,v) := {\kappa_2 \zeta(u,v)},\\[0.2cm]
&\zeta(u,v) := \frac{\lambda_1 k^- u}{B(u,v)^2(\lambda_0 + B(u,v))}, \quad B(u,v) := \kappa_2 v +\kappa_1 u + k^-, 
\end{aligned}\qquad \right\}\\[0.5cm]
\label{eq:pre_normal}
&\left.\begin{aligned}
\textbf{Normal tissue}& \textbf{ density model:} \\[0.2cm]
&\partial_t v =  -\kappa_3 \psi(u) v-R_Q(\xi_t) v, \quad \text{in }(0,T)\times \dom\hspace*{2cm} \\[0.2cm]
&v(0) = v_0\quad \text{in } \dom\\[0.2cm]
\text{\rm with} \qquad \quad& \\[0.2cm]
&\psi(u) := \frac{u}{1 + u}, \quad 
\end{aligned}\qquad \right\}\\[0.5cm]
\label{eq:radiation_term}
&\left.\begin{aligned}
\textbf{Radiation}& \textbf{ model:} \\[0.2cm]
&R_M(\xi_t) := \xi_t \sum_{i=1}^{n_{\mathrm f}} \Big(1-S_M(d_r,\brho_M) \Big) f_\delta(t-t_i) \\
&R_Q(\xi_t) := \xi_t \sum_{i=1}^{n_{\mathrm f}} \Big(1-S_Q(d_r,\brho_Q) \Big) f_\delta(t-t_i)
\end{aligned}\qquad \right\} \\[2ex]
\label{eq:noise}
&\left.\begin{aligned} 
    &\qquad \qquad \quad d \xi_t = \left(\mu_R-\xi_t\right) d t+\sigma_R \xi_t \: d W_t \\ 
    &\qquad \qquad \quad \xi(0)  = \xi_0 
\end{aligned} \quad \right\}
\end{align}
\end{subequations}
}
\end{define}

\begin{rem}[Derivative of radiation terms]\label{rem:radiation_term_dt}
With respect to equations \eqref{eq:radiation_term} and \eqref{eq:noise}, It\^o's formula yields the stochastic differentials of $R_M$ and $R_Q$ as follows:
\begin{align}
\label{eq:radiation_term_dt}
\left.\begin{aligned}
dR_M(\xi_t) &:= \sum_{i=1}^{n_{\mathrm f}} \Big(1-S_M(d_r,\brho_M) \Big) f_\delta(t-t_i) ~ d\xi_t  + \xi_t \sum_{i=1}^{n_{\mathrm f}} \Big(1-S_M(d_r,\brho_M) \Big) f'_\delta(t-t_i) ~ dt\\
dR_Q(\xi_t) &:= \sum_{i=1}^{n_{\mathrm f}} \Big(1-S_Q(d_r,\brho_Q) \Big) f_\delta(t-t_i) ~ d\xi_t + \xi_t \sum_{i=1}^{n_{\mathrm f}} \Big(1-S_Q(d_r,\brho_Q) \Big) f'_\delta(t-t_i) ~ dt
\end{aligned}\qquad \right\}
\end{align}
\end{rem}
\noindent In the following sections, we adopt this reformulated version of the macroscopic system.

\section{Analysis of the macroscopic model}\label{sec:analysis}
In this section, we address the question of local existence and uniqueness of strong solutions to the macroscopic system \eqref{eq:mult_glioma_sys} subjected to appropriate boundary conditions. To this end, we first introduce some notations and functional spaces that will be used throughout the analysis. Then, we state the problem setting along with the main assumptions on the model parameters. Finally, we present the main results regarding existence and uniqueness of solutions to the stochastic macroscopic system \eqref{eq:mult_glioma_sys}.

\subsection{Notations}\label{sec:notations}
We first introduce some basic notation and conventions. Let $I=(0, T] \subset \mathbb{R}_{+}$ be a finite time interval. Let $\dom_I:= I \times \dom$ denote the time-space cylinder on which the random PDE-ODE system is set. Throughout the analysis, $n\in\{1,2,3\}$ denotes the spatial dimension of $\dom$ and Sobolev embeddings are understood with respect to this dimension. Let $(\Om,\Fil,\P)$ be a complete probability space endowed with a normal filtration $(\Fil_t)_{t\geq0}$,  generated by the standard Brownian motion $(W_t)_{t\ge0}$, i.e.,
\[ 
\Fil_0\supseteq\{A\in\Fil:\P(A)=0\}, \qquad \Fil_t=\Fil_{t+}=\bigcap_{s>t}\Fil_s, \quad t\geq0. 
\] 
The exponential Ornstein-Uhlenbeck process $(\xi_t)_{t\ge0}$ in \eqref{eq:noise} is driven by $(W_t)_{t\ge0}$.
Partial derivatives, in both classical and distributional sense, with respect to time $t$ and space $x_i \in \dom$ variables, are denoted by $\partial_t$ and $\partial_{x_i}$, respectively, where $x_i$ is the $i$-th spatial coordinate. Moreover, $\nabla, \nabla \cdot$ and $\Delta$ stand for the spatial gradient, divergence and Laplace operators, respectively, while $\partial_\nu$ denotes the derivative in the direction of the outward unit normal $\nu$ on $\partial \dom$. We denote the Lebesgue measure of a set $A$ by $|A|$ and by int $A$ its interior. We assume familiarity with the standard $L_p$ spaces, Sobolev spaces $W^k_p$, H\"older spaces, and their standard properties, as well as with more general $L_p$ spaces of functions with values in general Banach spaces and with anisotropic Sobolev spaces. For the more specialized spaces used in the analysis, including the Bessel potential space $\wsp$, bounded Bochner spaces $\B$, time H\"older spaces $\cnt^{\vrho}$ with $\vrho<1$, and weighted H\"older spaces $\G^{\beta,\vsig}$, we refer to the Appendix~\ref{sec:appendix} (in particular see Definitions~\ref{def:BesPotSpace},~\ref{def:bdBochSp},~\ref{def:timeHoldSp}, and~\ref{def:weighHol}). Based on the introduced notation, we consider for our analysis the following Lebesgue and Sobolev spaces:
\begin{align}\label{eq:functional_spaces}
\begin{aligned}
{\V_p} &:= \left(L_p(\mathfrak{D}),\|\cdot\|_{L_p}\right), \qquad\quad{\V_p}(T):= \left(L_p(I \times \mathfrak{D}),\|\cdot\|_{L_{p}(T)}\right) \equiv L_p\left(I ; L_p(\mathfrak{D})\right)\,, \\[0.2cm] 
\Z_p(T) &:= \B([0,T]; \V_p), \qquad\qquad\,\,\, \Z^s_p(T) := \B([0,T]; \wsp), \\[0.2cm] 
\Y^{\beta,\vsig} &:= \G^{\beta,\vsig}((0,T];\V_p), \qquad\qquad\,\,\quad \Y := \Y^{1,\vsig}\,,\\[0.2cm] 
\cnt^{\vrho}_{s,p}(T) &:= \cnt^{\vrho}([0,T];\wsp) \text{ for } \vrho < 1, \quad \W^s_p :=\wsp(\mathfrak{D})\,,\\[0.2cm]
\mathbf{K}^{s,p}_{\rm b} & := \Big \{ w \in \wsp: {\exists\ \const_w>0\ \text{s.t. }}\|w\|_{\wsp} \le \const_w\Big \},\\[0.2cm]  
\bU(T) &:= \Big \{u \in \ztp(T): u \in \mathcal{C}([0,T];\W^s_p) \cap \mathcal{C}^{\vrho}([0,T];\W^{r}_p) \cap \mathcal{C}^{1}([0,T];\V_p) \Big \}, \\[0.2cm]
\X(T) &:= L_p\left(\Om; \zsp(T) \right) ~ \bigcap ~ L_p\left(\Om; \cnt^{\vrho}_{r,p}(T) \right), \quad \\[0.2cm] \K^p_{\loc}(\Zb) &:= \Big \{ w \in \Zb(T) :\E \left[\|w\|^p_{\zsp} \right] \le \const^p_{\Zb} \Big \}, \\[0.2cm] 
\Zb(T) &:= \left \{w \in \X(T) \cap L_p\left(\Om; \ztp \right):w(0) = u_0, ~ \E \left[\left(\sup{0\le s < t \le T} \frac{\|w_t - w_s\|_{\wrp}}{(t-s)^{\vsig}}\right)^p\right] \le 1 \right \},
\end{aligned}
\end{align}
{{where $\const_w, \const_\Zb, \const_{\X}$ are sufficiently large positive constants and the H\"older exponent $\vsig$ satisfies $\vsig < \min\{1/p,\,1/2\}$. The spaces $\X(T)$ and $\mathcal{Y}(T) := \X(T) \cap L_p(\Om;\ztp(T))$ are endowed with the norms

$$
\|w\|_{\X(T)} := \|w\|_{L_p(\Om;\zsp(T))} + \|w\|_{L_p(\Om;\cnt^{\vrho}_{r,p}(T))}
$$
and

$$\|w\|_{\mathcal{Y}(T)} := \|w\|_{\X(T)} + \|w\|_{L_p(\Om;\ztp(T))}$$
respectively. 
Both $(\X(T),\|\cdot\|_{\X(T)})$ and $(\mathcal{Y}(T),\|\cdot\|_{\mathcal{Y}(T)})$ are Banach spaces, since $\zsp(T)$, $\cnt^{\vrho}_{r,p}(T)$, and $\ztp(T)$ are Banach (Definitions~\ref{def:bdBochSp} and~\ref{def:timeHoldSp} in Appendix~\ref{sec:appendix}). For $w \in \Zb(T)$ we denote the pathwise $\vsig$-H\"older seminorm by

$$
[w]_{\vsig} := \sup{0\le s < t \le T} \frac{\|w_t - w_s\|_{\wrp}}{(t-s)^{\vsig}}.
$$}}

\noindent For convenience we summarize the main spaces and exponents used throughout Section~\ref{sec:analysis}.

\begin{table}[h!]
\centering
\small
\begin{tabular}{|l|l|l|}
\toprule
Symbol & Meaning & Typical range / role \\
\midrule
$n$ & spatial dimension of $\dom$ & $n\in\{1,2,3\}$ \\[0.05cm]
$p$ & Lebesgue/Sobolev integrability & $p>n$ (ensures $1+n/p<2$) \\[0.05cm]
$r,s,s_1,s_2$ & Sobolev regularity exponents & $1+\tfrac{n}{2p}<r<s<s_1<s_2=2$ \\[0.05cm]
$\vsig$ & time H\"older exponent (paths / data) & $\vsig<\min\{1/p,1/2\}$; further restricted by~\eqref{eq:U_Soln_vsig_Yagi} \\[0.05cm]
$\vrho$ & solution H\"older exponent in $\W^r_p$ & $\vrho=1-\alpha$ (cf.\ Theorem~\ref{thm:U_Soln}) \\[0.05cm]
$\V_p$ & $L_p(\dom)$ & base Banach space \\[0.05cm]
$\W^s_p$ & Bessel potential / Sobolev space & spatial regularity \\[0.05cm]
$\Z^s_p(T)$ & $\B([0,T];\W^s_p)$ & bounded Bochner paths \\[0.05cm]
$\X(T)$ & $L_p(\Om;\Z^s_p)\cap L_p(\Om;\cnt^{\vrho}_{r,p})$ & stochastic path space \\[0.05cm]
$\mathcal{Y}(T)$ & $\X(T)\cap L_p(\Om;\Z^2_p(T))$ & fixed-point ambient space \\[0.05cm]
$\Zb(T)$ & closed subset of $\mathcal{Y}(T)$ & $w(0)=u_0$ and stochastic $\vsig$-H\"older constraint \\[0.05cm]
$A(w,v)$ & sectorial elliptic operator & frozen variables $(w,v)$ \\[0.05cm]
$F(u,v)$ & semilinear reaction--taxis map & right-hand side of~\eqref{eq:ACP_Can} \\
\bottomrule
\end{tabular}
\caption{Main symbols used in the analysis in Section~\ref{sec:analysis}.}
\label{tab:notation_analysis}
\end{table}

\subsection{Problem setting and assumptions}\label{sec:problem_setting}
Before stating the formal hypotheses, we briefly summarize the main modelling and analytical requirements underlying the subsequent analysis. As mentioned above, we consider a fixed bounded spatial domain with a sufficiently smooth boundary and a finite time interval. The diffusion tensor describing tissue anisotropy is assumed to possess suitable spatial regularity and to be uniformly elliptic, ensuring coercivity of the second‑order part of the differential operator and boundedness of related multiplication operators. The taxis and decay coefficients are sufficiently smooth and uniformly bounded; moreover, appropriate positivity assumptions are imposed whenever required by the biological interpretation of the corresponding mechanisms, preventing degeneracies in the diffusion and taxis contributions.

Radiotherapy effects are represented through smooth temporal fractionation profiles with compact support and spatially regular coefficients. The stochastic modulation $\xi_t$ is the Ornstein--Uhlenbeck process introduced above. It has uniformly bounded moments and almost surely H\"older continuous sample paths on compact time intervals, allowing the radiotherapy operators to be interpreted pathwise. The initial data are assumed to be nonnegative and to belong to the Sobolev spaces specified below. We impose homogeneous no-flux boundary conditions for the tumor equation, excluding boundary transport and ensuring consistency with the mass conservation at the boundary.

From an analytical viewpoint, we work in Lebesgue and Sobolev spaces $L_p$ and $\W^{s,p}$ with $p>n$ and $s>1+n/p$, which guarantee the embeddings and maximal regularity properties needed for the semigroup and evolution operator techniques used in the existence and uniqueness proofs. All assumptions are biologically reasonable and consistent with the parameter regimes used in the companion modelling manuscript~\cite{conte2026randomness}. They are stated formally as follows.

\begin{assume}[Initial data]\label{asm:init_cond}
\begin{enumerate}
    \item[]
    \item $u_0 \ge 0$ and $u_0 \in \wtp(\dom)$,
    \item $v_0 \in [0,1]$ and $v_0 \in \wtp(\dom)$.
\end{enumerate}    
\end{assume}

\begin{assume}[Model parameters]\label{asm:model_params}
The coefficients appearing in the above model \eqref{eq:mult_glioma_sys}  satisfy the following conditions. 
\begin{enumerate}
\item \label{Assumpt-3.4.1.} Noise coefficient: $\sigma_R < 1$, with mean-reversion level $\mu_R>0$ and $\Fil_0$-measurable initial value $\xi_0\ge 0$ for the exponential Ornstein-Uhlenbeck process \eqref{eq:noise}, driven by $(W_t)_{t\ge0}$. {The SDE \eqref{eq:noise} admits a unique strong solution $(\xi_t)_{t\ge0}$, which is $(\Fil_t)$-adapted with continuous sample paths. In particular, $\xi_t\ge 0$ almost surely for all $t\ge 0$.} 
{\allowdisplaybreaks
\item \label{Assumpt-3.4.2.} Diffusion coefficients: 
\begin{align*}
    \mathbb{D}_M(x) &= \bmat{c^{11}(x) & \cdots~ & c^{1n}(x) \\
                \vdots & \ddots & \vdots \\
                c^{n1}(x) & \cdots~ & c^{nn}(x) \\
                } \in \R^{n\times n}, \quad 
    \mathbb{D}_1(x) = \bmat{\grad c^{11}(x) & \cdots~ & \grad c^{1n}(x) \\
                \vdots & \ddots & \vdots \\
                \grad c^{n1}(x) & \cdots~ & \grad c^{nn}(x) \\
                } \in \R^{n^2\times n}, \quad \\[2ex]
   \mathbb{D}_2(x) &= \bmat{\grad^2 c^{11}(x) & \cdots~ & \grad^2 c^{1n}(x) \\
                \vdots & \ddots & \vdots \\
                \grad^2 c^{n1}(x) & \cdots~ & \grad^2 c^{nn}(x) \\
                } \in \R^{n^3\times n}, 
\end{align*}}
where 
{\allowdisplaybreaks
\begin{align}
\label{eq:diff_ten_const}
    \hspace*{-0.5cm}\begin{array}{l l l l r}
        &\|c^{kl}(\cdot)\|_{\cnt^2} \le \tilde \ubd_{c}, &\overset{n}{\underset{k,l}{\sum}}~\|c^{kl}\|_{\cnt^2} \le \ubd_c, & c^{kl}(\cdot) \in \cnt^2(\dom), \: k,l \in \{1,\dots,n\},\\[2ex]
        &\lbd_D |\cgam|^2 < \cgam^{\top} \mathbb{D}_M(x) \cgam < \widetilde \ubd_D |\cgam|^2, & |\cgam_1^{\top} \mathbb{D}_M(x) \cgam_0| < \widetilde \ubd_D |\cgam_1 \cdot \cgam_0|, &\cgam, \cgam_0, \cgam_1 \in \R^n,\: x \in \dom, &\\[2ex]
        &|\cgam^{\top}_2 \mathbb{D}_1(x) \cgam_2| < \widetilde \ubd_D |\cgam_2|^2,& |\cgam^{\top}_3 \mathbb{D}_2(x) \cgam_3| < \widetilde \ubd_D |\cgam_3|^2, &\cgam_2 \in \R^{{n^2}}, \cgam_3 \in \R^{{n^3}}, &\\[2ex]
        &|\cgam^{\top}_{2_1} \mathbb{D}_1(x) \cgam_{2_2}| < \widetilde \ubd_D |\cgam_{2_1} \cdot \cgam_{2_2}|,  &|\cgam^{\top}_{3_1} \mathbb{D}_2(x) \cgam_{3_2}| < \widetilde \ubd_D |\cgam_{3_1} \cdot \cgam_{3_2}|, & \cgam_{2_1}, \cgam_{2_2} \in \R^{{n^2}}, \cgam_{3_1}, \cgam_{3_2} \in \R^{{n^3}}, & \\[2ex]
        &\ubd_D := \max(\ubd_c, \widetilde \ubd_D) < \infty. & & \\[-4ex]
    \end{array}
\end{align}
}
\item \label{Assumpt-3.4.5.} The spatial dose $d_r \in \wsp(\dom)$ and the temporal fractionation profile $f_\delta \in \cnt^{\vsig}([0,\infty);\R)$, for {{$\vsig < \min\{1/p,\,1/2\}$}}. Consequently, the survival factors $S_M(d_r,\brho_M)$, $S_Q(d_r,\brho_Q)$ and the coefficients
$$
\kappa_1(x) := k_2^+ S_M(d_r(x),\brho_M),\qquad \kappa_2(x) := k_1^+ S_Q(d_r(x),\brho_Q)
$$
from Definition~\ref{def:mult_glioma_sys} belong to $\wsp(\dom)$. In particular, since $s>1+n/p$, one has $\kappa_1,\kappa_2\in C^1(\bar\dom)$ with $0\le \kappa_i \le \ubd_{\kappa}:= \max\{k_1^+,k_2^+\}$ (using $0<S_i\le 1$).
\item \label{Assumpt-3.4.3.} With $\kappa_1,\kappa_2\in\wsp(\dom)$ as in Definition~\ref{def:mult_glioma_sys} and Assumption~\ref{asm:model_params}.\ref{Assumpt-3.4.5.}, the maps $(u,v)\mapsto\zeta(u,v;x)$ and $(u,v)\mapsto g(u,v;x):=\kappa_2(x)\zeta(u,v;x)$ from~\eqref{eq:pre_cancer} are $C^{\infty}$ in $(u,v)$ for a.e.\ $x\in\dom$. Moreover, there exists $\ubd_g<\infty$ (depending on $\|\kappa_1\|_{\wsp}$, $\|\kappa_2\|_{\wsp}$, $\lambda_1$, $\lambda_0$, and $k^-$) such that

$$
0\le g(\cdot,\cdot;x) \le \ubd_g,\qquad
\|\partial_u g(\cdot,\cdot;x)\|_{L_{\infty}(\R^2)},\, \|\partial_v g(\cdot,\cdot;x)\|_{L_{\infty}(\R^2)} \le \ubd_g
\quad\text{for a.e.\ }x\in\dom,
$$
and the Nemytskii operators $(u,v)\mapsto\zeta(u,v;\cdot)$ and $(u,v)\mapsto g(u,v;\cdot)$ are locally Lipschitz from $\W^{r}_p(\dom)\times\W^{r}_p(\dom)$ into $\W^{r}_p(\dom)$ for every $r\in\bigl(\tfrac{n}{p},s\bigr]$, with Lipschitz constants controlled by $\ubd_g$. Writing $\nabla g=\partial_u g\,\nabla u+\partial_v g\,\nabla v+\nabla^{\kappa}g$ for the full spatial gradient of the composed map (where $\nabla^{\kappa}g$ collects derivatives through the spatial coefficients $\kappa_1,\kappa_2$), one has $\|\nabla^{\kappa}g(u,v;\cdot)\|_{L_{\infty}(\dom)}\le \ubd_g$ whenever $u,v$ remain bounded in $\W^{r}_p$.
\item \label{Assumpt-3.4.6.} The constants $\mu_M, \kappa_3, \lambda_1, \lambda_0, k^-, k_1^+, k_2^+ \in \R^+$ are strictly positive. 
\item \label{Assumpt-3.4.7.} In order to ensure the boundary condition for \eqref{eq:pre_cancer} we make the following assumptions for $D_1$, $v_0$ and $R_Q$:

$$ 
D_1 = 0, \quad \nu \cdot (g \,\mathbb{D}_M \nabla v_0) = 0 \quad \text{and} \quad \nu \cdot (g\, \mathbb{D}_M \nabla R_Q) = 0 \quad \text{on } \partial \dom.
$$
\end{enumerate}
\end{assume}
\noindent
An immediate consequence of Assumption \ref{asm:model_params}.\ref{Assumpt-3.4.1.} is that the Ornstein–Uhlenbeck process $\xi_t$ has uniformly bounded moments and H\"older continuous sample paths on compact time intervals, as summarized in Corollary \ref{cor:lp_hold_cont_proc} provided in the Appendix \ref{sec:appendix}. As a result, we can deduce the following boundedness property for the radiation terms $R_M$ and $R_Q$ defined in \eqref{eq:radiation_term} and serving as decay rates for glioma and tissue.

\begin{lem}[]\label{rad_asure_bded}
Let $R_M$ and $R_Q$ be as defined in \eqref{eq:radiation_term}. Then $R_M$ and $R_Q$ are pathwise bounded and are thus elements of $L_p\big(\Om, \cnt^{\vsig}([0,T];\wsp)\big)$.
In particular, we have that
\begin{align}
    \|R_M\|_{\cnt^{\vsig}} \le \const_{R_M} \quad ~~\text{\rm and }\quad \|R_Q\|_{\cnt^{\vsig}} \le \const_{R_Q} \quad ~~\text{\rm for a.a. $\om \in \Om$.}
\end{align}
\end{lem}
\begin{proof}
    By Corollary~\ref{cor:lp_hold_cont_proc} we have that $\xi \in L_p(\Om; \cnt^{\vsig}([0,T];\R))$ and, from Assumption \ref{asm:model_params}.\ref{Assumpt-3.4.5.}, we have that the spatial function $d_r \in \wsp(\dom)$ and the temporal function $f_\delta \in \cnt^{\vsig}([0,\infty),\R)$. Since the radiotherapy terms are obtained as products of H\"older continuous functions of time and spatial coefficients belonging to $\wsp(\dom)$, then $R_M,R_Q\in L_p\big(\Omega;\cnt^{\vsig}([0,T];\wsp(\dom))\big)$, which proves the claim.
\end{proof}

\newcommand{\sbm}{{\rm{\W}_{t\ge0}}}
\newcommand{\prt}{(\Om,\Fil,\P)}
\newcommand{\nfil}{(\Fil_t)_{t\ge0}}
\newcommand{\Ft}{\Fil_t}
\newcommand{\avery}{\text{\rm \rm{a.e.}~}}
\newcommand{\asure}{\text{\rm \rm{a.s.}}}
\newcommand{\aall}{\text{\rm \rm{a.a.}}~}
\newcommand{\rv}{\text{\rm \rm{r.v.}}~}
\newcommand{\and}{\text{\rm \rm{and}}~}
\newcommand{\psip}{\psi^{\prime}}
\newcommand{\psipp}{\psi^{\prime \prime}}
\newcommand{\cntvtp}{\cnt^{\vsig}_{2,p}}
\newcommand{\cntrrp}{\cnt^{\rho}_{r,p}}

\subsection{Apriori estimates}
Starting from the multiscale glioma model formulation~\eqref{eq:mult_glioma_sys}, we derive in this section a collection of apriori estimates that play a twofold role. From the modeling perspective, they ensure the preservation of biologically meaningful properties of the solutions, such as nonnegativity of the cell densities and uniform bounds preventing unphysical blow‑up, and control the influence of the stochastic radiotherapy term. From the analytical viewpoint, these estimates establish the quantitative regularity and Lipschitz bounds required to analyze the cancer cell equation within the framework of non-autonomous semigroup/evolution operators. In particular, uniform bounds on the normal tissue density $v$ and on the radiotherapy terms guarantee the uniform ellipticity and temporal regularity of the coefficients of the parabolic operator, while the Lipschitz estimates for the nonlinear terms allow to close the fixed‑point argument and prove local existence and uniqueness. Thus, the following estimates provide the link between the biological structure of the model and the functional-analytic assumptions required for the subsequent well-posedness 
results.

\subsubsection{Positivity and boundedness}
We first establish the positivity and boundedness properties of the normal tissue density $v$.

\begin{lem}[Solution to normal tissue equation]\label{lem:sol_normal}
Let $s\in (1,2)$, $p>n$, and let $u \in \cnt([0,T];\wsp)$ be any arbitrary function  (more generally, $u \in L_p(\Om;\cnt([0,T];\wsp))$). Let $\xi \in L_p\left(\Om; \cnt([0,T];\R)\right)$ be the stochastic process defined by~\eqref{eq:noise}, and let $v_0$ satisfies    Assumption~\ref{asm:init_cond}. Then, for a.a.\ $\om\in\Om$, the solution $v$ to \eqref{eq:pre_normal} takes the following form:
    \begin{align}
    \label{eq:soln_normal}
        v(t) = v_0 \exp \left(-\int_0^t \left(\kappa_3 \psi(u_s) + R_Q(\xi_s) \right)  \: ds \right) \quad ~~\text{\rm a.a. } \om \in \Om.
    \end{align}
Furthermore, if $v_0 \in [0,1]$ $\asure$, then $v(t) \in [0,1]$ $\asure$ for $\avery t \in [0,T]$.
\end{lem}

\begin{proof}
Since $u \in \cnt([0,T];\wsp)$ and $\wsp \hookrightarrow \cnt(\bar \dom)$, and since $\xi \in L_p(\Om; \cnt([0,T];\R))$, the equation for the normal tissue component can be written in a pathwise and pointwise manner as
\begin{equation}\label{eq:normal_rode}
\begin{aligned}
        \frac{d}{dt} v(t;\om,x) &= -\kappa_3 ~ \psi(u(t;\om,x)) v(t;\om,x) - R_Q(\xi_t(\om)) v(t;\om,x),~\avery x \in \dom,~\aall \om \in \Om,~\and t > 0 \\
        v(0) &= v_0.
\end{aligned}
\end{equation}
Thus, Eq.~\eqref{eq:normal_rode} can be seen as a random non-autonomous ODE, parameterized by $x \in \dom$. The solution to it can be obtained in a pathwise manner, i.e., for a.a.\ $\om\in\Om$ with $\xi(\om)\in\cnt([0,T];\R)$. Thus, $v(t;\om,x)$ takes the following explicit form:
\begin{align}\label{eq:rode_sol}
   v(t) &= v_0 e^{- \chi_{_t} {(u)}}, \quad \chi_{_t} {(u)} := \int_0^t \left(\kappa_3 \psi(u_s) + R_Q(\xi_s) \right) \: ds.\:
\end{align}
By Assumption~\ref{asm:model_params}.\ref{Assumpt-3.4.1.}, $\xi_t \ge 0$ almost surely, hence $R_Q(\xi_t) \ge 0$ almost surely. Since $\kappa_3 > 0$ and $\psi\ge 0$, we get $\chi_{_t} \ge 0$ almost surely. Thus, $v_0 \in [0,1]$ implies $v(t;\om,x) \in [0,1]$ for $\avery x \in \dom$, $\aall \om \in \Om$, and $t > 0$.
\end{proof}

\noindent To derive the apriori estimates for the solution to the normal tissue equation, we first compute the first and second spatial derivatives of $v$ and $\chi_{_t}$. By direct differentiation of~\eqref{eq:rode_sol} and for $i,j\in\{1,\dots,n\}$, we obtain
    \begin{align}
        \partial_{x_i} v&= \partial_{x_i} v_0 e^{-\chi_{_t}} - \partial_{x_i} \chi_{_t} v, \:\\[0.2cm] 
        \partial_{x_i} \chi_{_t}(u) &= \int_0^t \big({\kappa_3}\psi^{\prime}(u_s) \partial_{x_i} u_s + \partial_{x_i} R_Q\big) ~ ds,  \label{eq:grad_rode_sol} \\[0.2cm]
       \partial_{x_j} (\partial_{x_i} v) &= \partial^2_{x_j x_i} v_0 e^{-\chi_{_t}} - \partial_{x_j x_i}^2 \chi_{_t} v - \partial_{x_j} \chi_{_t} \partial_{x_i} v_0e^{-\chi_{_t}} - \partial_{x_i} \chi_{_t} \partial_{x_j} v_0e^{-\chi_{_t}} + \partial_{x_i} \chi_{_t} \partial_{x_j} \chi_{_t} v, \label{eq:lap_rode_sol} \\[0.2cm]
       \partial^2_{x_j x_i} \chi_{_t} &=  \int_0^t \Big({\kappa_3}\left (\psi^{\prime \prime}(u_s) \partial_{x_j} u_s \partial_{x_i} u_s + \psi^{\prime}(u_s) \partial^2_{x_j x_i} u_s\right ) + \partial^2_{x_j x_i} R_Q\Big) ~ ds. \label{eq:lap_chi_rode_sol}
\end{align}
Next, we derive some $L_p$ estimates for the normal tissue density $v$ that will be useful in the subsequent analysis.

\begin{lem}[$L_p$ estimate of $v$]\label{lem:lp_estimate_v} 
Let $v_0 \in \wtp,$ $0 \le v_0(x) \le 1$, $\chi_{_t} \ge 0\ \asure$. Let $s \in (1,2]$, and let 

$$u \in L_p(\Om;\ztp(T)) \quad \text{ and } \quad R_Q \in L_p(\Om;\cnt^{\vsig}([0,T];\wsp)) ~~~ (cf.\ Lemma~\ref{rad_asure_bded}).$$ 
Then, for a.a.\ $\om \in \Om$,

$$ \ztp(T) \ni u(\om) \mapsto v(u(\om)) \in \cnt^{\vsig}\left([0,T];\wtp\right), \quad \vsig < {1 \over p}, $$
such that 
{\allowdisplaybreaks
\begin{align}
\label{eq:v_nu_lp_est}
\|v\|_{\cntvtp} &\le \Bigl(\const^p(v_0) + \const^p_v(T,u,R_Q) T\Bigr)^{1 \over p} + \const_{v,tL}(T,u,R_Q) T^{1-\vsig}\\[0.2cm]
{\const^p(v_0)} &:= {\|v_0\|^p_{\wtp} \left( 1 + 2n\, C_p\, C_{\mathrm{emb}}^{2p}\, \|v_0\|^{p}_{\wtp}\right)}   \\[0.2cm]   
{\const^p_v(T,u,R_Q)} &:= \Bigg( C_{\chi,1}(T)\, \ubd^p_{v}\, \kappa^p_3\, \ubd^p_{\psi}\, \|u\|^{p}_{\ztp} + C_{\chi,2}(T)\, \ubd^p_{v}\, \kappa^p_3\, \ubd^p_{\psi}\, \biggl( \|u\|^{p}_{\ztp} + 2n\, C_{\mathrm{emb}}^{2p}\, \|u\|^{2p}_{\ztp}\biggr) \\[0.2cm]
&\hspace*{0.5cm} + C_{\chi,2p}(T)\, 2n\,(1 + \ubd^p_v)\, \bigl( \kappa^{2p}_3\, \ubd^{2p}_{\psi}\, \|u\|^{2p}_{\ztp} + \|R_Q\|^{2p}_{\cnt_{2,p}} \bigr) \\[0.2cm]
&\hspace*{0.5cm} + \bigl(C_{\chi,1}(T) + C_{\chi,2}(T)\bigr)\, \ubd^p_{v}\, \|R_Q\|^p_{\cnt_{2,p}} \Bigg),\\[0.2cm]
{\const_{v,tL}(T,u,R_Q)} &:= C_{\chi,*}(T)\, \Bigl( \kappa_3\, \ubd_{\psi}\, \|u\|_{\ztp} + \|R_Q\|_{\cnt_{2,p}} \Bigr),
\end{align}
}
where $\ubd_{\psi}$ and $\ubd_{v}$ are the $L_{\infty}$ bounds on $\psi(u)$ (and its derivatives) and on $v_0$, respectively. The factors $C_{\mathrm{Y},p}$, $C_{\mathrm{H}}(T)$, $C_p$, $C_{\mathrm{emb}}$, $C_{\mathrm{Y},p}^{(3)}$, $C_{\mathrm{Y},2p}$, $C_{\mathrm{H},2p}(T)$ and the composites $C_{\chi,1}(T)$, $C_{\chi,2}(T)$, $C_{\chi,2p}(T)$ are defined in the proof below (cf.~\eqref{eq:v_abs_const}--\eqref{eq:v_comp_const}), and $C_{\chi,*}(T)$ is introduced with the H\"older bound~\eqref{eq:v_time_lip}.
\end{lem}

\begin{proof}\label{lem:proof_lp_estimate_v}
As per \eqref{eq:pre_normal}, $u \mapsto \psi(u):=u/(1+u)$ based on which we have that $\psi,\psip,\psipp\in L_{\infty}(\R)$. Let

$$\ubd_{\psi}:= \max\left(\|\psi\|_{L_{\infty}(\R)}, \|\psip\|_{L_{\infty}(\R)}, \|\psipp\|_{L_{\infty}(\R)} \right) \text{ and } \ubd_{v}:=\|v_0\|_{L_{\infty}(\dom)}.$$ 
Let the constants appearing from Young's inequality, H\"older's inequality, and the multinomial inequality be denoted as 
\begin{align}
\label{eq:v_abs_const}
\begin{aligned}
C_{\mathrm{Y},p}&:=2^{p-1},&
C_{\mathrm{H}}(t)&:=t^{p-1},&
C_p&:=5^{p-1},\\[0.1cm]
C_{\mathrm{Y},p}^{(3)}&:=3^{p-1},&
C_{\mathrm{Y},2p}&:=2^{2p-1},&
C_{\mathrm{H},2p}(t)&:=t^{2p-1}.
\end{aligned}
\end{align}
Further, let the composites appearing in~\eqref{eq:v_lp_est_1}-\eqref{eq:v_lp_est_2} be denoted as
\begin{align}
\label{eq:v_comp_const}
\begin{aligned}
C_{\chi,1}(T)&:=C_{\mathrm{H}}(T)\,C_{\mathrm{Y},p},&
C_{\chi,2}(T)&:=C_p\,C_{\mathrm{H}}(T)\,C_{\mathrm{Y},p}^{(3)},&
C_{\chi,2p}(T)&:=C_p\,C_{\mathrm{H},2p}(T)\,C_{\mathrm{Y},2p}\,C_{\mathrm{emb}}^{2p}
\end{aligned}
\end{align}
where $C_{\mathrm{emb}}=C_{\mathrm{emb}}(\dom,n,p)>0$ is the constant appearing in the embedding $\wtp\hookrightarrow\W^1_{2p}$ ($p>n/2$).
Based on this, applying the Minkowski (triangle) inequality and H\"older's inequality in $\V_p$, we estimate the norm of $v$ and its first derivatives as:
\begin{align}
    \label{eq:v_lp_est_1}
    \begin{aligned}
    \|v\|^p_{\V_p} &\leq \|v_0\|^p_{\V_p} \|e^{-\chi_t}\|^p_{L_\infty(\dom)} \le \|v_0\|^p_{\V_p} \\[0.2cm]
    \|\partial_{x_i} v\|^p_{\V_p} &\leq \|\partial_{x_i} v_0\|^p_{\V_p} \|e^{-\chi_t}\|^p_{L_\infty(\dom)} + \|\partial_{x_i} \chi_t\|^p_{\V_p} \|v_0\|^p_{L_\infty(\dom)} \|e^{-\chi_t}\|^p_{L_\infty(\dom)} \\[0.2cm]
    &\leq \|\partial_{x_i} v_0\|^p_{\V_p} + \ubd^p_{v} \|\partial_{x_i} \chi_t\|^p_{\V_p} \\[0.2cm]
\implies \sum_{i=1}^n \|\partial_{x_i} v\|^p_{\V_p} &\leq \sum_{i=1}^n \|\partial_{x_i} v_0\|^p_{\V_p} + \ubd^p_{v} \sum_{i=1}^n \|\partial_{x_i} \chi_t\|^p_{\V_p}.
    \end{aligned}
\end{align}
By H\"older (w.r.t. time) and Young's inequality on the integrand of~\eqref{eq:grad_rode_sol},
\begin{align*}
\|\partial_{x_i} \chi_t\|^p_{\V_p}
&\le C_{\mathrm{H}}(t)\int_0^t \bigl\|{\kappa_3}\psi'(u_s)\partial_{x_i}u_s+\partial_{x_i}R_Q\bigr\|^p_{\V_p}\,ds \\[0.2cm]
&\le C_{\mathrm{H}}(t)\,C_{\mathrm{Y},p}\int_0^t \Bigl({\kappa_3}^p\|\psi'\|_{L_\infty(\dom)}^p\|\partial_{x_i}u_s\|^p_{\V_p}+\|\partial_{x_i}R_Q\|^p_{\V_p}\Bigr)\,ds,
\end{align*}
and therefore, using $C_{\mathrm{H}}(t)\le C_{\mathrm{H}}(T)$ for $t\in[0,T]$,
\begin{align*}
\sum_{i=1}^n \|\partial_{x_i} v\|^p_{\V_p}
&\le \sum_{i=1}^n \|\partial_{x_i} v_0\|^p_{\V_p} + \ubd^p_{v}\, C_{\mathrm{H}}(T)\, C_{\mathrm{Y},p} \int_{0}^t \left({\kappa^p_3} \ubd^p_{\psi} \sum_{i=1}^n \|\partial_{x_i} u_s\|^p_{\V_p} + \sum_{i=1}^n \|\partial_{x_i} R_Q\|^p_{\V_p}\right) ~ ds.
\end{align*}
Since $W^{k}_p(\dom) \equiv \W^k_p(\dom)$ for $k \in \N$ and since $\|v\|^p_{\W^1_p} \equiv \|v\|^p_{W^1_p} := \|v\|^p_{\V_p} + \sum \limits_{i=1}^n \|\partial_{x_i} v\|^p_{\V_p}$, we can combine the above estimates to get that
\begin{align*}
\|v(t)\|^p_{\W^1_p} &\leq \|v_0\|^p_{\wop} + \ubd^p_{v}\, C_{\mathrm{H}}(T)\, C_{\mathrm{Y},p} \int_{0}^t \big({\kappa^p_3} \ubd^p_{\psi} \|u_s\|^p_{\wop} + \|R_Q\|^p_{\wop}\big) ~ ds, \\[0.2cm]
    &\leq \|v_0\|^p_{\V_p} + \underbrace{\ubd^p_{v}\, C_{\chi,1}(T)  \Big( {\kappa^p_3} \ubd^p_{\psi} ~ \|u\|^p_{\zop(T)} + \|R_Q\|^p_{\cnt_{1,p}(T)} \Big)t}_{\text{from }\|\chi_t\|^p_{W^{1,p}}}.
\end{align*}
Similarly, for the second derivatives, we apply the multinomial inequality with $C_p=5^{p-1}$ to the five terms in~\eqref{eq:lap_rode_sol} and H\"older's inequality on products to obtain:
{\allowdisplaybreaks 
\begin{align}\label{eq:v_lp_est_2}
        \|\partial^2_{x_j x_i} v\|^p_{\V_p} &\leq C_p\Big( \|\partial^2_{x_j x_i} v_0\|^p_{\V_p} + \|\partial^2_{x_j x_i} \chi_t\|^p_{\V_p} \|v_0\|^p_{L_\infty(\dom)} \\[0.2cm]
    &\quad + \|\partial_{x_j} \chi_t\|^{2p}_{\V_{2p}} + \|\partial_{x_i} v_0\|^{2p}_{\V_{2p}} + \|\partial_{x_i} \chi_t\|^{2p}_{\V_{2p}} + \|\partial_{x_j} v_0\|^{2p}_{\V_{2p}} \\[0.2cm]
    &\quad + (\|\partial_{x_i} \chi_t\|^{2p}_{\V_{2p}} + \|\partial_{x_j} \chi_t\|^{2p}_{\V_{2p}}) \|v\|^p_{L_\infty(\dom)} \Big) \\[0.2cm]
\implies \sum_{i,j=1}^n \|\partial^2_{x_j x_i} v\|^p_{\V_p} &\leq C_p\Big( \sum_{i,j=1}^n \|\partial^2_{x_j x_i} v_0\|^p_{\V_p} +  \ubd_v^p \sum_{i,j=1}^n \|\partial^2_{x_j x_i} \chi_t\|^p_{\V_p}  \\[0.2cm]
    &\quad + \sum_{i,j=1}^n \|\partial_{x_j} \chi_t\|^{2p}_{\V_{2p}} + \sum_{i,j=1}^n \|\partial_{x_i} v_0\|^{2p}_{\V_{2p}} + \sum_{i,j=1}^n \|\partial_{x_i} \chi_t\|^{2p}_{\V_{2p}} + \sum_{i,j=1}^n \|\partial_{x_j} v_0\|^{2p}_{\V_{2p}} \\[0.2cm]
    &\quad + \ubd^p_{v} \sum_{i,j=1}^n (\|\partial_{x_i} \chi_t\|^{2p}_{\V_{2p}} + \|\partial_{x_j} \chi_t\|^{2p}_{\V_{2p}})  \Big),
\end{align}
}
where, 
\begin{align}\label{eq:lap_chi_est}
    \begin{aligned}
&\|\partial^2_{x_j x_i} \chi_t\|^p_{\V_p} \leq C_{\mathrm{H}}(t)\,C_{\mathrm{Y},p}^{(3)}\int_0^t \Bigg({\kappa^p_3}\ubd^p_{\psi} \Big((\|\partial_{x_j} u_s\|^{2p}_{\V_{2p}} + \|\partial_{x_i} u_s\|^{2p}_{\V_{2p}}) + \|\partial^2_{x_j x_i} u_s\|^p_{\V_p} \Big) + \|\partial^2_{x_j x_i} R_Q\|^p_{\V_p}\Bigg) ~ ds \\[0.2cm]
&\implies \sum_{i,j} \|\partial^2_{x_j x_i} \chi_t\|^p_{\V_p} \leq C_{\mathrm{H}}(T)\,C_{\mathrm{Y},p}^{(3)}\int_0^t \Bigg({\kappa^p_3} \ubd^p_{\psi} \Big( 2n \sum_{j=1}^n \|\partial_{x_j} u_s\|^{2p}_{\V_{2p}} + \sum_{i,j=1}^n \|\partial^2_{x_j x_i} u_s\|^p_{\V_p} \Big) + \sum_{i,j=1}^n \|\partial^2_{x_j x_i} R_Q\|^p_{\V_p}\Bigg) ~ ds\,.
    \end{aligned}
\end{align}
Likewise, the $L_{2p}$ bound on first derivatives of $\chi_t$ reads
\begin{align*}
\sum_{i=1}^n \|\partial_{x_i} \chi_t\|^{2p}_{\V_{2p}}
&\le C_{\mathrm{H},2p}(T)\,C_{\mathrm{Y},2p}\int_0^t \Bigl({\kappa_3}^{2p}\ubd_{\psi}^{2p}\sum_{i=1}^n\|\partial_{x_i}u_s\|^{2p}_{\V_{2p}}+\sum_{i=1}^n\|\partial_{x_i}R_Q\|^{2p}_{\V_{2p}}\Bigr)\,ds.
\end{align*}
So, combining the estimates~\eqref{eq:v_lp_est_1}, \eqref{eq:v_lp_est_2}, and \eqref{eq:lap_chi_est} we get that 
{\allowdisplaybreaks
\begin{align}
    \|v(t)\|^p_{\wtp} &= \|v(t)\|^p_{\V_p} + \sum_{i=1}^n \|\partial_{x_i} v(t)\|^p_{\V_p} + \sum_{i,j=1}^n \|\partial^2_{x_j x_i} v(t)\|^p_{\V_p} \\[0.2cm]
    &\le \|v_0\|^p_{\wtp} + C_{\chi,1}(T)\, \ubd^p_{v} \int_0^t \Big( \kappa^p_3\, \ubd^p_{\psi}\, \|u_s\|^p_{\wtp} + \|R_Q\|^p_{\cnt_{2,p}} \Big)\,ds \\[0.2cm]
    &+ C_{\chi,2}(T)\, \ubd^p_v \int_0^t \Big( \kappa^p_3\, \ubd^p_{\psi}\, \|u_s\|^p_{\wtp} + 2 n \kappa^p_3\, \ubd^p_{\psi}\, \|u_s\|^{2p}_{\W^1_{2p}} + \|R_Q\|^p_{\wtp}\Big)\,ds \\[0.2cm]
    &+ C_p\, 2n (1 + \ubd^p_v)\, C_{\mathrm{H},2p}(T)\,C_{\mathrm{Y},2p} \int_0^t \big({\kappa^{2p}_3} \,\ubd^{2p}_{\psi} \,\|u_s\|^{2p}_{\W^1_{2p}} + \|R_Q\|^{2p}_{\W^1_{2p}}\big)\,ds \\[0.2cm]
    &+ C_p\, 2n \sum_{i=1}^n \|\partial_{x_i} v_0\|^{2p}_{\V_{2p}} \\[0.2cm]
\implies \|v(t)\|^p_{\wtp} &\le \|v_0\|^p_{\wtp} + C_p\, 2n \|v_0\|^{2p}_{\W^1_{2p}}   \\[0.2cm]
    & + C_{\chi,1}(T)\, \ubd^p_{v} \int_0^t \Big( \kappa^p_3\, \ubd^p_{\psi}\, \|u_s\|^p_{\wtp} + \|R_Q\|^p_{\wtp} \Big)\,ds \\[0.2cm]
    & + C_{\chi,2}(T)\, \ubd^p_{v} \int_0^t \Big( \kappa^p_3\, \ubd^p_{\psi}\, \|u_s\|^p_{\wtp} + 2n\kappa^p_3\, \ubd^p_{\psi}\, \|u_s\|^{2p}_{\W^1_{2p}} + \|R_Q\|^p_{\wtp}\Big)\,ds  \\[0.2cm]
    &+ C_p\, C_{\mathrm{H},2p}(T)\,C_{\mathrm{Y},2p}\, 2n (1 + \ubd^p_v) \int_0^t \big({\kappa^{2p}_3} \,\ubd^{2p}_{\psi} \,\|u_s\|^{2p}_{\W^1_{2p}} + \|R_Q\|^{2p}_{\W^1_{2p}}\big)\,ds.
\end{align}}
Since $\wtp \hookrightarrow \W^1_{2p}$, for $p > n/2$ we have $\|w\|_{\W^1_{2p}}\le C_{\mathrm{emb}}\|w\|_{\wtp}$, and therefore
\begin{align}
    \|v(t)\|^p_{\wtp} &\le \|v_0\|^p_{\wtp} + C_p\, 2n\, C_{\mathrm{emb}}^{2p}\, \|v_0\|^{2p}_{\W^2_{p}}   \\[0.2cm]
    & + C_{\chi,1}(T)\, \ubd^p_{v} \int_0^t \Big( \kappa^p_3\, \ubd^p_{\psi} \,\|u_s\|^p_{\wtp} + \|R_Q\|^p_{\wtp} \Big) ds \\[0.2cm]
    & + C_{\chi,2}(T)\, \ubd^p_{v} \int_0^t \Big( \kappa^p_3\, \ubd^p_{\psi} \,\|u_s\|^p_{\wtp} + 2n\kappa^p_3\, \ubd^p_{\psi}\, C_{\mathrm{emb}}^{2p}\, \|u_s\|^{2p}_{\wtp} + \|R_Q\|^p_{\wtp}\Big) ds  \\[0.2cm]
    &+ C_{\chi,2p}(T)\, 2n (1 + \ubd^p_v) \int_0^t \big({\kappa^{2p}_3} \,\ubd^{2p}_{\psi}\, \|u_s\|^{2p}_{\wtp} + \|R_Q\|^{2p}_{\wtp}\big)ds.
\end{align}
Grouping the terms, we get that
\begin{align}
    \|v(t)\|^p_{\wtp} \le& \,\|v_0\|^p_{\wtp} \left( 1 + 2n\, C_p\, C_{\mathrm{emb}}^{2p}\, \|v_0\|^{p}_{\wtp}\right) \\[0.2cm]
    &+ C_{\chi,1}(T)\, \ubd^p_{v}\, \kappa^p_3 \,\ubd^p_{\psi} \int_0^t \|u_s\|^{p}_{\wtp}\, ds + C_{\chi,2}(T)\, \ubd^p_{v}\, \kappa^p_3 \,\ubd^p_{\psi} \int_0^t \Big( \|u_s\|^{p}_{\wtp} + 2n\, C_{\mathrm{emb}}^{2p}\, \|u_s\|^{2p}_{\wtp}\Big)\, ds \\[0.2cm]
    &+ C_{\chi,2p}(T)\, 2n (1 + \ubd^p_v) \int_0^t \big({\kappa^{2p}_3} \,\ubd^{2p}_{\psi}\, \|u_s\|^{2p}_{\wtp} + \|R_Q\|^{2p}_{\wtp}\big)\,ds \\[0.2cm]
    &+ \bigl(C_{\chi,1}(T)+C_{\chi,2}(T)\bigr)\, \ubd^p_{v} \int_0^t \|R_Q\|^p_{\wtp}\, ds.
\end{align}
Thus, using the assumptions $u \in \ztp$ and $R_Q \in \cnt_{2,p}$ we get that
\begin{align}
    \|v(t)\|^p_{\wtp} &\le \overbrace{\|v_0\|^p_{\wtp} \left( 1 + 2n\, C_p\, C_{\mathrm{emb}}^{2p}\, \|v_0\|^{p}_{\wtp}\right)}^{\const^p(v_0)} \\[0.2cm]
    &\quad + \underbrace{\begin{pmatrix}
    C_{\chi,1}(T)\, \ubd^p_{v}\, \kappa^p_3\, \ubd^p_{\psi}\, \|u\|^{p}_{\ztp} + C_{\chi,2}(T)\, \ubd^p_{v}\, \kappa^p_3\, \ubd^p_{\psi}\, \biggl( \|u\|^{p}_{\ztp} + 2n\, C_{\mathrm{emb}}^{2p}\, \|u\|^{2p}_{\ztp}\biggr) + \\[0.1cm]
    C_{\chi,2p}(T)\, 2n\,(1 + \ubd^p_v)\, \bigl( \kappa^{2p}_3\, \ubd^{2p}_{\psi}\, \|u\|^{2p}_{\ztp} + \|R_Q\|^{2p}_{\cnt_{2,p}} \bigr) + \bigl(C_{\chi,1}(T) + C_{\chi,2}(T)\bigr)\, \ubd^p_{v}\, \|R_Q\|^p_{\cnt_{2,p}}
    \end{pmatrix}}_{\const^p_v(T,u,R_Q)}~t \\[-0.2cm]
\implies \|v(t)\|^p_{\wtp} &\le \const^p(v_0) + \const^p_v(T,u,R_Q) t, \text{ for a.a. $\om \in \Om$ fixed. }
\end{align}

\noindent Based on this, the same differentiation formulae~\eqref{eq:grad_rode_sol}-\eqref{eq:lap_chi_rode_sol} on a subinterval $[s,t]\subset[0,T]$ yield a Lipschitz bound in time for the path $v$. Indeed, writing $\chi_t-\chi_s=\int_s^t\bigl(\kappa_3\psi(u_\tau)+R_Q(\xi_\tau)\bigr)\,d\tau$ and repeating the $W^2_p$ estimates of the preceding display with integrals restricted to $[s,t]$, one obtains
\begin{align}
\|v(t) - v(s)\|_{\wtp} &\le C_{\chi,*}(T)\, \Bigl( \kappa_3\, \ubd_{\psi}\, \|u\|_{\ztp} + \|R_Q\|_{\cnt_{2,p}} \Bigr)\, |t-s| = \const_{v,tL}(T,u,R_Q)\, |t-s|. \label{eq:v_time_lip}
\end{align}
where $C_{\chi,*}(T):=C_{\chi,*}(\ubd_v,C_p,C_{\mathrm{emb}},T)>0$ collects the H\"older's, Young's, multinomial and embedding prefactors similar to those in \eqref{eq:v_comp_const}.

\noindent Consequently, for every $\vsig\in\bigl(0,1/p\bigr)$ and $|t-s|\le T$,
\begin{align}
\frac{\|v(t) - v(s)\|_{\wtp}}{|t-s|^{\vsig}} &\le \const_{v,tL}(T,u,R_Q)\, |t-s|^{1-\vsig} \le \const_{v,tL}(T,u,R_Q)\, T^{1-\vsig}, \notag\\[0.2cm]
\implies \| v \|_{\cntvtp} &\le \left(\const^p(v_0) + \const^p_v(T,u,R_Q) T\right)^{1 \over p} + \const_{v,tL}(T,u,R_Q)\, T^{1-\vsig} \label{eq:v_cvsig_h2p}.
\end{align}
Since \eqref{eq:v_cvsig_h2p}  holds for \aall $\om \in \Om$ and $\vsig < {1 \over p}$ we have that 

$$ \B\left([0,T];\wtp\right) \ni u \mapsto v(u) \in \cnt^{\vsig}\left([0,T];\wtp\right).$$
\end{proof}

\noindent Next, we derive some Lipschitz estimates for the normal tissue density $v$ that will be useful later in the fixed point argument underlying the existence and uniqueness analysis.

\begin{lem}[Lipschitz inequality]\label{lem:Lip_ineq_normal}
Let $v_0 \in \wtp$, $0 \le v_0(x) \le 1$, $\chi_{_t} \ge 0 ~ \asure$, $s \in (1,2]$ and 
$$
u(\om) \in \B([0,T];\wtp), \quad R_Q(\om) \in \cnt^{\vsig}([0,T];\wtp) \quad \text{ \rm{for a.a.} $\om \in \Om$ }.
$$
Write $\ubd_{\psi}$ as in Lemma~\ref{lem:lp_estimate_v}, $\ubd_{v}:=\|v_0\|_{L_{\infty}(\dom)}\le 1$, and $ \const_{v,L}(t)\,:=\,\ubd_{v}\,t.$
Then there exists a constant $C_{v,L}\!\ge\!\kappa_3$, depending on $\kappa_3$, $C_p$, $C_{\mathrm{emb}}$, $\|v_0\|_{\wtp}$, $\|u\|_{\B([0,T];\wtp)}$ such that
    \begin{align}
    \label{eq:v_l2_lipschitz}
    \begin{aligned}
     \|v^1(t) - v^2(t)\|_{\wtp} &\le C_{v,L}\,\ubd_{\psi}\, \const_{v,L}(t) ~ \|u^1 - u^2\|_{\B([0,T];\wtp)} \\[0.2cm]
\|v^1 - v^2\|_{\cnt^{\vsig}([0,T];\wtp)} &\le C_{v,L}\,\ubd_{\psi}\, \ubd_{v}\, T^{1-\vsig} ~ \|u^1 - u^2\|_{\B([0,T];\wtp)}.
    \end{aligned}
    \end{align}
    
\end{lem}
\begin{proof}\label{lem:proof_Lip_ineq_normal}
Starting from \eqref{eq:rode_sol} and noting that $R_Q$ is independent of $u$, we have the identities
\begin{align*}
v^1(t)-v^2(t)
&= v_0\bigl(e^{-\chi^1_t}-e^{-\chi^2_t}\bigr),\\[0.1cm]
\chi^1_t-\chi^2_t
&= \kappa_3\int_0^t\bigl(\psi(u^1_s)-\psi(u^2_s)\bigr)\,ds.
\end{align*}
By the mean value theorem applied pathwise and pointwise, there exists an intermediate field $\tilde\chi_t$ between $\chi^1_t$ and $\chi^2_t$ such that

$$
e^{-\chi^1_t}-e^{-\chi^2_t}=-e^{-\tilde\chi_t}\,(\chi^1_t-\chi^2_t),
$$
and therefore

$$
v^1(t)-v^2(t)=-v_0\,e^{-\tilde\chi_t}\,(\chi^1_t-\chi^2_t).
$$
Since $\chi^i_t\ge 0$ and $0\le v_0\le\ubd_{v}\le 1$, one has $0\le e^{-\tilde\chi_t}\le 1$. Consequently,
\begin{align*}
\|v^1(t)-v^2(t)\|_{\V_p}
&\le \ubd_{v}\,\|\chi^1_t-\chi^2_t\|_{\V_p}
= \ubd_{v}\,\kappa_3\,\Bigl\|\int_0^t\bigl(\psi(u^1_s)-\psi(u^2_s)\bigr)\,ds\Bigr\|_{\V_p}\\[0.1cm]
&\le \ubd_{v}\,\kappa_3\int_0^t\|\psi(u^1_s)-\psi(u^2_s)\|_{\V_p}\,ds,\\[-5ex]
\end{align*}
where the last step is Minkowski's integral inequality.
Using the pointwise bound $|\psi(a)-\psi(b)|\le\|\psip\|_{L_{\infty}}|a-b|$ and the definition ${\ubd_{\psi}\ge\|\psip\|_{L_{\infty}}}$ from Lemma~\ref{lem:lp_estimate_v}, we obtain

$$
\|v^1(t)-v^2(t)\|_{\V_p}\le \kappa_3\,\ubd_{v}\,\ubd_{\psi}\, t\,\|u^1-u^2\|_{\B([0,T];\wtp)}
= \kappa_3\,\ubd_{\psi}\,\const_{v,L}(t)\,\|u^1-u^2\|_{\B([0,T];\wtp)}.
$$
For the first spatial derivatives we start from~\eqref{eq:grad_rode_sol}. Since $R_Q$ is independent of $u$,
\begin{align*}
\partial_{x_i}\chi^1_t-\partial_{x_i}\chi^2_t
&= \kappa_3\int_0^t\Bigl(\psi'(u^1_s)\partial_{x_i}u^1_s-\psi'(u^2_s)\partial_{x_i}u^2_s\Bigr)\,ds \\[0.1cm]
&= \kappa_3\int_0^t\Bigl(\psi'(u^1_s)\,\partial_{x_i}(u^1_s-u^2_s)+\bigl(\psi'(u^1_s)-\psi'(u^2_s)\bigr)\partial_{x_i}u^2_s\Bigr)\,ds.
\end{align*}
Minkowski's integral inequality in $\V_p$, the pointwise bound $|\psi'(a)-\psi'(b)|\le\|\psipp\|_{L_{\infty}}|a-b|$, and the embedding $\wtp\hookrightarrow\W^1_{2p}$ (with constant $C_{\mathrm{emb}}$ as in Lemma~\ref{lem:lp_estimate_v}) yield

$$
\|\partial_{x_i}\chi^1_t-\partial_{x_i}\chi^2_t\|_{\V_p}
\le \kappa_3\,\ubd_{\psi}\,t\,\Bigl(1+C_{\mathrm{emb}}\|u^2\|_{\B([0,T];\wtp)}\Bigr)\,\|u^1-u^2\|_{\B([0,T];\wtp)}.
$$
Using the mean-value representation of $e^{-\chi^1_t}-e^{-\chi^2_t}$, the bound $0\le v^2(t)\le\ubd_{v}$, the $L_p$ estimate already proved, and writing $C_{v,L}$ for a constant that absorbs $\kappa_3$, $C_{\mathrm{emb}}$ and the uniform $\B$-bounds on $u$, we obtain
\begin{align}
\label{eq:v_lip_grad}
\|\partial_{x_i}(v^1(t)-v^2(t))\|_{\V_p}
&\le C_{v,L}\,\ubd_{\psi}\,\const_{v,L}(t)\,\|u^1-u^2\|_{\B([0,T];\wtp)}.
\end{align}

\noindent Likewise for the second derivative we obtain
\begin{align}
\label{eq:v_lip_hess}
\|\partial_{x_ix_j}^2(v^1(t)-v^2(t))\|_{\V_p}
&\le C_{v,L}\,\ubd_{\psi}\,\const_{v,L}(t)\,\|u^1-u^2\|_{\B([0,T];\wtp)},
\end{align}
where $C_{v,L}$ may additionally depend on $C_p$ and $\|v_0\|_{\wtp}$. Combining the $L_p$ bound with~\eqref{eq:v_lip_grad}--\eqref{eq:v_lip_hess} we conclude
\begin{align}
\label{eq:v_lipschitz}
\begin{aligned}
\|v^1(t) - v^2(t)\|_{\wtp} &\le C_{v,L}\,\ubd_{\psi}\, \const_{v,L}(t) \|u^1 - u^2\|_{\B([0,T];\wtp)}.
\end{aligned}
\end{align}
Moreover, the same argument applied to the signed increments of $\chi$ and its derivatives over an arbitrary subinterval $[s,t]\subset[0,T]$ yields

$$
\|(v^1-v^2)(t)-(v^1-v^2)(s)\|_{\wtp}
\le C_{v,L}\,\ubd_{\psi}\,\ubd_{v}\,|t-s|\,\|u^1-u^2\|_{\B([0,T];\wtp)}.
$$
Consequently, for every $\vsig\in(0,1)$,
\begin{align*}
\frac{\|(v^1-v^2)(t)-(v^1-v^2)(s)\|_{\wtp}}{|t-s|^{\vsig}}
\le C_{v,L}\,\ubd_{\psi}\,\ubd_{v}\,|t-s|^{1-\vsig}\,\|u^1-u^2\|_{\B([0,T];\wtp)}\\
\le C_{v,L}\,\ubd_{\psi}\,\ubd_{v}\, T^{1-\vsig}\,\|u^1-u^2\|_{\B([0,T];\wtp)}.
\end{align*}
Together with~\eqref{eq:v_lipschitz} this gives
\begin{align}
\label{eq:v_lipschitz_hold}
\|v^1 - v^2\|_{\cnt^{\vsig}([0,T];\wtp)} &\le C_{v,L}\,\ubd_{\psi}\, \ubd_{v}\, T^{1-\vsig} \|u^1 - u^2\|_{\B([0,T];\wtp)},
\end{align}
which is the second claim in~\eqref{eq:v_l2_lipschitz}.
\end{proof}

\newcommand{\tfor}{~~~\text{\rm for}~}

\subsection{Abstract reformulation}\label{sec:abstract_reform}
The estimates derived above provide the analytical foundation for the functional framework adopted throughout the remainder of this section, while also reflecting the underlying biological assumptions incorporated into the model. In particular, the positivity and uniform boundedness of $v$, together with the pathwise H\"older continuity and moment estimates for the radiotherapy process $\xi_t$ (hence for $R_M, R_Q$), and the Lipschitz continuity and regularity properties of the solution operator $u\mapsto v(u)$ ensure that the coefficients appearing in the tumor equation are uniformly bounded, spatially regular, and temporally H\"older continuous (for almost every realization of the noise). Together with the uniform ellipticity and smoothness of the diffusion tensor, these properties guarantee that the principal part of the associated parabolic operator is uniformly coercive and that the resulting non‑autonomous family of operators arising satisfies the resolvent and H\"older continuity assumptions required by the theory of evolution operators.  

These results allow us to reformulate the tumor growth equation as a non-autonomous Cauchy problem in $L_p$, thereby providing the appropriate setting for applying the evolution operator machinery. This framework naturally motivates the use of Sobolev and weighted H\"older spaces in the search for local solutions. We now proceed to make this reformulation precise and subsequently establish the local existence and uniqueness of solutions.

To this end, let $p \in (n,\infty)$ and $s \in (1 + {n \over p}, 2)$. For given functions $w$, $v$ and $R_M$, the tumor evolution equation~\eqref{eq:pre_cancer} can be reformulated as the following semilinear non-autonomous Cauchy problem:
\begin{subequations}
\begin{align}
\label{eq:ACP_Can}
\begin{aligned}
    {d \over dt} u(t) + A(w(t), v(t))\, u(t) &= F(u(t),v(t)) ~~\text{\rm in } \V_p \\
    u(0) &= u_0,
\end{aligned}
\end{align}
with
\begin{align}
\label{eq:ACP_Ops}
\begin{aligned}
    A(w,v) u &:= -\grad \cdot (D_0(w,v) \grad u) + {u}, \quad \forall w \in \wsp(\dom) \hookrightarrow \cnt^1(\bar \dom), \\[0.2cm]
    \D(A(w,v)) &= \Big\{ u \in \W^2_p(\dom):{\nu} \cdot (D_0(w,v) \grad u) = 0 \ {\text{ on }\partial \dom }\Big \}, \\[0.2cm]
    F(u,v) &:= {u + \grad \cdot (D_1  u) - \grad \cdot (g(u,v) \mathbb{D}_M \grad v ) } + {{\mu_M u(1-u)v}}-R_M(\om) u, \quad \forall u,v \in \W^s_p(\dom),
\end{aligned}
\end{align}
where $D_0(w,v):=\mathbb{D}_M(1+\kappa_1\zeta(w,v))$ denotes the diffusion coefficient from~\eqref{eq:pre_cancer} with $u$ replaced by the frozen cancer density $w$.
\end{subequations}
\newcommand{\pit}{{\pi \over 2}}
To proceed, we need to ensure that the mapping $u\mapsto D_0(w,v)\,u$ defines a bounded linear operator from $\W^s_p(\dom)$ into $\wrp(\dom)$ for given functions $w$ and $v$. This is established in the following lemma. 

\begin{lem}[Operator norm of diffusion tensor]\label{lem:diff_ten}
Let $\zeta$ be as in~\eqref{eq:pre_cancer}, let $\kappa_1,\kappa_2\in\wsp(\dom)$ be as in Assumption~\ref{asm:model_params}.\ref{Assumpt-3.4.5.}, and let 

$$ w(\om),v(\om) \in \cnt^{\vsig}([0,T];\W^{r}_p(\dom)), ~~\text{\rm for a.a. $\om \in \Om$} $$
with $p \in \big(n,\: \infty\big)$, $r \in \Big({n \over p}, \infty \Big)$ and $\mathbb{D}_M=(D_M^{i,j})_{i,j=1,\dots,n} \in \R^{n\times n}$ as in Assumption \ref{asm:model_params}. Then for a.a. $\om \in \Om$
\begin{align*} &D^{i,j}_0(w(\om),v(\om)) := D_{ M}^{i,j}\Big(1 + \kappa_1 \zeta\big(w(\om),v(\om)\big)\Big) \in \L(\wsp,\wrp),\quad {\frac{n}{p} < }  ~r \le s, ~{s \ge 1}\,, ~~\forall i,j \in \{1, \dots, n\}.
\end{align*}
\end{lem}

\begin{proof}\label{lem:proof_diff_ten}
Since $r>\frac np$, the Sobolev embedding theorem yields $\W^{r}_p(\dom) \hookrightarrow C(\bar \dom)$. Because $\zeta$ is $C^{\infty}$ in $(u,v)$ and $\kappa_1,\kappa_2\in\wsp\subset\W^{r}_p$ (Assumption~\ref{asm:model_params}.\ref{Assumpt-3.4.5.}), the Nemytskii composition
$x \mapsto \zeta\big(w({t,}x),v({t,}x);x\big)$ belongs to $\W^{r}_p(\dom)$, and therefore $D^{i,j}_0 \in \W^{r}_p(\dom)$. Thus, the multiplication operator $u \mapsto D^{i,j}_0 u$ is a bounded linear operator, i.e., $D^{i,j}_0 \in \L(\wsp,\wrp)$ for any $s \ge 1$, ${n \over p} < r \le s$. Indeed, we have that the inequality
    \begin{align*}
    \|D^{i,j}_0 u\|_{\wrp} &\le C \|D^{i,j}_0\|_{\W^{r}_p}\|u\|_{\wsp} \le C \|D^{i,j}_0\|_{\W^{r}_p}\|u\|_{\wsp}
    \end{align*}
holds in the almost sure sense. On the other hand, for $w_1, w_2, v_1, v_2 \in \W^{r}_p$ and using that $\W^{r}_p$ is a Banach algebra for $r>n/p$, 
{\allowdisplaybreaks
    \begin{align*}
        \Big\|D^{i,j}_0(w_1, v_1) - D^{i,j}_0(w_2,v_2) \Big\|_{\W^{r}_p}& = \Big\|\kappa_1 D^{i,j}_M [\zeta(w_1, v_1) - \zeta(w_2,v_2)] \Big\|_{\W^{r}_p}\\[0.2cm]
        &\hspace*{-2cm}\le C \|\kappa_1\|_{\W^{r}_p} \big\|D^{i,j}_M \big\|_{\cnt^{1}} ~ \Big\|\zeta(w_1, v_1) - \zeta(w_2,v_2) \Big\|_{\W^{r}_p} \\[0.2cm]
        &\hspace*{-2cm}\le C \|\kappa_1\|_{\W^{r}_p} \big\|D^{i,j}_M \big\|_{\cnt^{1}} ~ \left(\Big\|\partial_w \zeta(\tilde w,\tilde v) (w_1 - w_2)\Big\|_{\W^{r}_p}  + \Big\|\partial_v \zeta(\tilde w, \tilde v) (v_1 - v_2) \Big\|_{\W^{r}_p}\right ) \\[0.2cm]
        &\hspace*{-2cm}\le C \|\kappa_1\|_{\W^{r}_p} \big\|D^{i,j}_M \big\|_{\cnt^{1}} ~ \bigg(\Big\|\partial_w \zeta(\tilde w,\tilde v)\Big\|_{L_{\infty}} 
        \big\|w_1 - w_2\big\|_{\W^{r}_p}  + \Big\|\partial_v \zeta(\tilde w, \tilde v)\Big\|_{L_{\infty}} \big\| v_1 - v_2 \big\|_{\W^{r}_p}\bigg) \\[0.2cm]
        &\hspace*{-2cm}\le {C^{i,j}_{\mathbb{D}_M}} \Big (C_{\zeta_w} \big\|w_1 - w_2 \big\|_{\W^{r}_p}  +  C_{\zeta_v} \big\|v_1 - v_2 \big\|_{\W^{r}_p}\Big ),
    \end{align*}
}
where $C_{\zeta_w},C_{\zeta_v}$ absorb the uniform bounds on $\partial_w\zeta,\partial_v\zeta$ (which depend only on $\|\kappa_1\|_{L_{\infty}}$, $\|\kappa_2\|_{L_{\infty}}$, $\lambda_1$, $\lambda_0$, and $k^-$ by Assumption~\ref{asm:model_params}.\ref{Assumpt-3.4.3.}--\ref{Assumpt-3.4.5.}), and $C^{i,j}_{\mathbb{D}_M}$ absorbs $\|\kappa_1\|_{\W^{r}_p}$ together with $\|D_M^{i,j}\|_{\cnt^1}$.
Additionally, if $ w \mapsto v(w) $ is a bounded nonlinear operator with 

$$
\|v(w)\|_{\W^{r}_p} \le C_v \|w\|_{\W^{r}_p}
$$ 
and 

$$
\big\|v_1 - v_2\big\|_{\W^{r}_p} = \big\|v(w_1) - v(w_2)\big\|_{\W^{r}_p} \le C_v \big\|w_1 - w_2\big\|_{\W^{r}_p},
$$ 
then we get that 
\begin{align*}
\big\|D^{i,j}_0(w_1, v_1(w_1)) - D^{i,j}_0(w_2,v_2(w_2))\big\|_{\W^{r}_p} &\le  {C^{i,j}_{\mathbb{D}_M}} (C_{\zeta_w} +  C_{\zeta_v} C_v ) \big\|w_1 - w_2\big\|_{\W^{r}_p} \\[0.2cm]
\implies \sup{i,j} \big\|D^{i,j}_0(w_1, v_1(w_1)) - D^{i,j}_0(w_2,v_2(w_2))\big\|_{\W^{r}_p} &\le  \sup{i,j} C^{i,j}_{\mathbb{D}_M} (C_{\zeta_w} +  C_{\zeta_v} C_v) \big\|w_1 - w_2\big\|_{\W^{r}_p} \\[0.2cm]
\implies \sup{i,j} \big\|D^{i,j}_0(w_1, v_1(w_1)) - D^{i,j}_0(w_2,v_2(w_2))\big\|_{\L(\wsp,\wrp)} &\le  \sup{i,j} C^{i,j}_{\mathbb{D}_M} (C_{\zeta_w} +  C_{\zeta_v} C_v) \big\|w_1 - w_2\big\|_{\W^{r}_p} \\[0.2cm]
\implies \big\|D_0(w_1, v_1(w_1)) - D_0(w_2,v_2(w_2))\big\|_{\L(\wsp,\wrp)} &\le  C_{\mathbb{D}_M} (C_{\zeta_w} +  C_{\zeta_v} C_v) \big\|w_1 - w_2\big\|_{\W^{r}_p}.
\end{align*}
\end{proof}

\noindent As an immediate consequence, we obtain the following remark, which establishes that the diffusion coefficient $D_0(w,v(w))$ is Lipschitz with respect to $w$.
\begin{rem}{}{}
    Specializing Lemma~\ref{lem:diff_ten} to any $r$ with $\frac{n}{p}<r\le 1$ (in particular $r=1$ when $p>n$), one obtains, for $v\mapsto v(w)$ as above,
    
    $$ 
    \big\|D_0(w_1, v_1(w_1)) - D_0(w_2, v_2(w_2))\big\|_{\L(\wop,\wop)} \le  C_{\mathbb{D}_M} \, (C_{\zeta_w} + C_{\zeta_v} C_v ) \, \|w_1 - w_2\|_{\wop}.
    $$
    In the main analysis we take $s>1+{n\over p}$ (so $\kappa_1,\kappa_2\in C^1$) and $r\in\bigl(1+{n\over 2p},s\bigr)$. The above result is only to clarify the $W^{1}_p$ Lipschitz bound used in intermediate estimates.
\end{rem}
\newcommand{\constloc}{\const_{\rm{loc}}}

\noindent In order to establish the existence of an evolution operator, we need to verify that the family of elliptic operators $A(w,v)$ satisfies the resolvent estimates and H\"older continuity properties. These properties are summarized in the following lemma.

\begin{lem}[Properties of $A(w,v)$]\label{lem:Au_property}
Let $(\Om,\Fil,\P)$ be a probability space generated by $(W_t)_{t\ge0}$. Let ${\om \mapsto v(\om,t) \in \W^2_p(\dom)}$ and  $\om \mapsto w(\om,t) \in \wsp$ be given functions with $v \in [0,1]$ and $w \ge 0$ almost surely, $s \ge 1 + {n \over p}$. Then, for $r \in \big({n \over p}, s \big)$, the operator 

$$A\big(w(\om), v(\om)\big) \in \L\bigl(\D(A(w,v));L_p \bigr)$$ 
satisfies the following properties for almost all $\om \in \Om$:
\begin{align}\label{eq:Au_struc_property}
\begin{aligned}
\sigma(A(w,v)) \subset \Sigma_{\pit} &:= \Big\{\clam \in \mathbb{C}: |\text{\rm arg} \clam| < {\pi \over 2} \Big\}, \\[0.2cm]
\Big\|(\clam - A(w,v))^{-1} \Big\|_{\L(\V_p)} &\le {\const_A \over |\clam| }, \quad \clam \notin \Sigma_{\pit}, \\[0.2cm]
\D(A(w_1,v_1)) = \D(A(w_2,v_2)) &~~\text{\rm for}\,\, w_i\in \wsp,\; v_i \in \W^2_p(\dom),~ i \in \{1,2\}, \\[0.2cm]
\D(A(w_1,v_1)) \subset \D(A^{\vartheta}(w_2,v_2)) &~~ \text{\rm for } \vartheta \in (0,1], \\[0.2cm]
\Big\|A(w_1,v_1)^{\vartheta} \big[ A(w_1,v_1)^{-1} - A(w_2,v_2)^{-1}\big] \Big\|_{\L(\V_p)} &\le \const_{A} \Big( \|w_1 - w_2\|_{\W^{r}_p} + \|v_1 - v_2\|_{\W^{r}_p} \Big), ~~s > r > {n \over p}, 
\end{aligned}
\end{align}
Additionally, for $v_i = v(w_i)$ and  $w \mapsto v(w)$ such that Lemma~\ref{lem:diff_ten} (cf.\ also Lemma~\ref{lem:Lip_ineq_normal}) holds, we get
\begin{align}
\begin{aligned}
\Big\|A(w_1,v(w_1))^{\vartheta} \big[ A(w_1,v(w_1))^{-1} - A(w_2,v(w_2))^{-1}\big] \Big\|_{\L(\V_p)} &\le \const_{A,w} \|w_1 - w_2\|_{\W^{r}_p}.
\end{aligned}
\end{align}
\end{lem}

\begin{proof}\label{lem:proof_Au_property}
{\allowdisplaybreaks 
Since $\W^s_p \hookrightarrow \cnt^1(\bar \dom)$ for $s \ge 1 + {n \over p}$, for given $w \in \W^s_p$ and $v \in \W^2_p(\dom)$, and since $\kappa_1\in\wsp(\dom)$ by Assumption~\ref{asm:model_params}.\ref{Assumpt-3.4.5.}, we have that $D_0(w,v) \in \cnt^1(\dom;\R^{n \times n})$. Additionally, since the diffusion matrix $D_0(w,v) := \mathbb{D}_M(1 + \kappa_1\zeta(w,v))$ is strongly elliptic and $1 \cdot D_0(w,v) \cdot 1 \in L_{\infty}(\dom)$, Theorem~2.12 of \cite{Yagi09} (Chapter~2, \S2.2, p.~58) is applicable, based on which we get that}
\begin{align*}
    \sigma(A(w,v)) \subset \Sigma_{\pit} &:= \Big\{ \clam \in \mathbb{C}: |\text{\rm arg} \clam| < {\pi \over 2} \Big\}, \\
\| (\clam - A(w,v))^{-1}\|_{\L(\V_p)} &\le {\const_A \over |\clam| }, \: \clam \notin \Sigma_{\pit}, \: \const_A \ge 1.
\end{align*}
Due to strong ellipticity of $\mathbb{D}_M$ and Lemma \ref{lem:diff_ten}, we also have that 

$$ 
\D(A(w_1,v_1)) = \D(A(w_2,v_2))
$$
for arbitrary $w_1, w_2 \in {\wsp}$ and $v_1, v_2 \in \W^2_p(\dom)$ (in particular also for $w_1, w_2 \in {\mathbf{K}^{s,p}_{\rm b}}$,cf.~\eqref{eq:functional_spaces}).
Because ${D^{i,j}_0(w,v) \in \mathcal{L}(\wsp,\wrp)}$, with $s,r$ as in Lemma~\ref{lem:diff_ten}, there exists $\vartheta \in (0,1]$ such that 

$$ 
\D(A(w_1,v_1)) \subset \D(A^{\vartheta}(w_2,v_2))
$$
for arbitrary $w_1, w_2 \in {\wsp}$ and $v_1, v_2 \in \W^2_p(\dom)$ (in particular also for $w_1, w_2 \in {\mathbf{K}^{s,p}_{\rm b}}$, cf.~\eqref{eq:functional_spaces}). In fact, $\D(A(w_1,v_1)) = \D(A(w_2,v_2)) = \W^2_{p,N} \subset \D(A^{\vartheta}(w,v)) \subset \W^{2\vartheta'}_{p,N}, ~ {(w,v)\in \{(w_1,v_1),(w_2,v_2)\},\ {\frac{p+1}{2p}<}\vartheta' < \vartheta}$ (cf.~\eqref{eq:dom_embds}).

\noindent Now let $p,q \in \N$ such that ${1 \over p } + {1 \over q} = 1$, $p>\max(2,n)$, and $\vartheta \in \Big (\frac{1}{2}+\frac{n}{2p},1\Big ]$. Then  for $f_1 \in L_p$ and $f_2 \in L_q$ we have that
{\allowdisplaybreaks
\begin{align*}
    &\Big| \Big\langle A(w_1,v_1)^{\vartheta} \big[A(w_1,v_1)^{-1} - A(w_2,v_2)^{-1}\big] f_1, f_2 \Big\rangle_{(L_p,L_q)} \Big | \\[0.2cm]
    &= {\Big| \Big\langle A(w_1,v_1)^{\vartheta} \big[A(w_1,v_1)^{-1}A(w_2,v_2)A(w_2,v_2)^{-1} - A(w_1,v_1)^{-1}A(w_1,v_1)A(w_2,v_2)^{-1}\big] f_1, f_2 \Big\rangle_{(L_p,L_q)} \Big |} \\[0.2cm]
    &= \Big\langle A(w_1,v_1)^{\vartheta-1} \big[A(w_2,v_2) - A(w_1,v_1)\big] A(w_2,v_2)^{-1} f_1, f_2 \Big\rangle_{(L_p,L_q)} \\[0.2cm]
    &= \Big\langle \big[A(w_2,v_2) - A(w_1,v_1)\big] A(w_2,v_2)^{-1} f_1, A^*(w_1,v_1)^{(\vartheta-1)}f_2 \Big\rangle_{(\W^{2\vartheta-2}_p, \W^{2-2\vartheta}_q)} \\[0.2cm]
    &= \Big\langle \big[D_{0}(w_2,v_2) - D_{0}(w_1,v_1)\big] \grad {\Big(}A(w_2,v_2)^{-1} f_1{\Big )}, \grad {\Big(}A^*(w_1,v_1)^{(\vartheta-1)}f_2{\Big )} \Big\rangle_{(\W^{2\vartheta-1}_p, \W^{1-2\vartheta}_q)} \\[0.2cm]
    &\le \Big\|[D_{0}(w_2,v_2) - D_{0}(w_1,v_1)]  \grad {\Big(}A(w_2,v_2)^{-1} f_1{\Big)}\Big\|_{\W^{2\vartheta-1}_p}\Big\|\grad {\Big(}A^*(w_1,v_1)^{(\vartheta-1)}f_2 {\Big)}\Big\|_{\W^{1-2\vartheta}_q} \\[0.2cm]
    &\le \Big\|D_{0}(w_2,v_2) - D_{0}(w_1,v_1)\Big\|_{\L(\W^1_p,\W^{2\vartheta-1}_p)} \Big\|\grad {\Big(}A(w_2,v_2)^{-1} f_1{\Big)}\Big\|_{\W^1_p} \Big\|\grad {\Big(}A^*(w_1,v_1)^{(\vartheta-1)}f_2 {\Big)}\Big\|_{\W^{1-2\vartheta}_q}.
\end{align*}}
As $\grad \in {\L(\W^{\vartheta}_q,\W^{\vartheta-1}_q)}$ and $ \D((A^*)^{1-\vartheta}) \subset \W^{2(1-\vartheta)}_q $ for $\vartheta \in (0,1)$ \text{(cf. (3.76) in \cite{Yagi09}, Chapter~3, \S3.3, p.~88)}, now we have that 
\begin{align*}
    \Big\|\grad \Big( A^*(w_1,v_1)^{(\vartheta-1)}f_2 \Big) \Big\|_{\W^{1-2\vartheta}_{q}}  &= \Big\|A^*(w_1,v_1)^{(\vartheta-1)}f_2 \Big\|_{\W^{2(1-\vartheta)}_q} \\[0.2cm]
    &\le C_{\vartheta} \Big\|A^*(w_1,v_1)^{(1 - \vartheta)} A^*(w_1,v_1)^{-(1 - \vartheta)}f_2 \Big\|_{L_q} = C_{\vartheta} \| f_2 \|_{L_q}.
\end{align*}
Similarly, since $\|u\|_{\W^2_p} \le C \|A(w,v) u\|_{L_p}$ we have that 
\begin{align*} 
\|\grad A(w,v)^{-1} f \|_{\W^1_p} = \| A(w,v)^{-1} f \|_{\W^2_p} \le C \|f\|_{L_p}.
\end{align*}
Altogether, this leads to 
\begin{align*}
    \Big| \Big\langle A(w_1,v_1)^{\vartheta} [A(w_1,v_1)^{-1} - A(w_2,v_2)^{-1}] f_1, f_2 \Big\rangle_{L_p \times L_q} \Big | &\le  C \Big\|D_{0}(w_2,v_2) - D_{0}(w_1,v_1)\Big\|_{\L(\W^1_p,\W^{2\vartheta-1}_{p})} \|f_1\|_{L_p} \|f_2\|_{L_q} \\[0.2cm]
    \Ra \|A(w_1,v_1)^{\vartheta} [A(w_1,v_1)^{-1} - A(w_2,v_2)^{-1}]\|_{\L(L_p)} &\le C  \Big\|D_{0}(w_2,v_2) - D_{0}(w_1,v_1)\Big\|_{\L(\W^1_p,\W^{2\vartheta-1}_{p})} \\[0.2cm]
    &\le C \ubd_D \Big( C_{\zeta_w} \|w_1 - w_2\|_{\W^{r}_p} + C_{\zeta_v} \|v_1 - v_2\|_{\W^{r}_p} \Big)\,, 
\end{align*}
where we have used Lemma~\ref{lem:diff_ten} and the definition of $D_0(w,v)$.

Additionally, if $w \mapsto v(w)$ is a bounded nonlinear operator with 

$$
\|v(w)\|_{\W^{r}_p} \le C_v \|w\|_{\W^{r}_p}
\quad\text{\rm and}\quad
\big\|v(w_1) - v(w_2)\big\|_{\W^{r}_p} \le C_v \big\|w_1 - w_2\big\|_{\W^{r}_p},
$$ 
then, writing $v_i := v(w_i)$, $i=1,2$, we have $\|v_1 - v_2\|_{\W^{r}_p} \le C_v \|w_1 - w_2\|_{\W^{r}_p}$ (cf.\ Lemma~\ref{lem:diff_ten} and Lemma~\ref{lem:Lip_ineq_normal}), and therefore
\begin{align*}
\Big\|A(w_1,v(w_1))^{\vartheta} \big[ A(w_1,v(w_1))^{-1} - A(w_2,v(w_2))^{-1}\big] \Big\|_{\L(\V_p)}
&\le C \,\ubd_D \,(C_{\zeta_w} + C_{\zeta_v} C_v) \|w_1 - w_2\|_{\W^{r}_p} \\
&= \const_{A,w} \|w_1 - w_2\|_{\W^{r}_p}.
\end{align*}
\end{proof}

\noindent As a consequence of Lemma~\ref{lem:Au_property}, we obtain the following corollary, which establishes the H\"older continuity of the resolvent mapping $t \mapsto A(w_t,v_t)^{-1}$ in the operator norm.

\begin{cor}[]\label{cor:Au_Holder}
Let $s\ge 1+{n\over p}$, $r \in \big({n \over p}, s \big)$, and let $w \in \B([0,T];\wsp)\cap\cnt^{\vsig}([0,T];\wrp)$ and $v \in \cnt^{\vsig}([0,T];\W^{2}_p)$. Then, there exists a constant $\const_{A,\vsig}$ such that 

$$\Big\|A(w_t,v_t)^{\vartheta} \big[ A(w_t,v_t)^{-1} - A(w_s,v_s)^{-1}\big] \Big\|_{\L(\V_p)} \le \const_{A,\vsig} ~ |t - s|^{\vsig}. $$
\end{cor}
\begin{proof}
    Denoting $A(t) := A(w_t,v_t)$ and applying Lemma \ref{lem:Au_property} we obtain 
    \begin{align*}
        \Big\|A(w_t,v_t)^{\vartheta} \big[ A(w_t,v_t)^{-1} - A(w_s,v_s)^{-1}\big] \Big\|_{\L(\V_p)} &\le \const_{A} \Big( \|w_t - w_s\|_{\W^{r}_p} + \|v_t - v_s\|_{\W^{r}_p} \Big) \\[0.2cm]
    \implies \Big\|A(t)^{\vartheta} \big[ A(t)^{-1} - A(s)^{-1}\big] \Big\|_{\L(\V_p)} &\le \const_{A,\vsig} ~ |t-s|^{\vsig},
    \end{align*}
    where the final inequality follows from the H\"older regularity of $w$ and $v$.
\end{proof}

Next, we establish the regularity properties of the nonlinear operator $F(u,v)$ appearing on the right-hand side of~\eqref{eq:ACP_Can}.

\begin{lem}[Properties of $F(u,v)$]\label{lem:Fuv_property}
For a fixed $\om \in \Om$, let $u(\om) \in \W^r_p(\dom)$, with $r \in \left(1, \infty \right)$, and $v \in \W^2_p(\dom)$ with $v(x) \in [0,1]$. Then the mapping 

$$
\Big(\W^r_p(\dom) \times \W^2_p(\dom) \Big) \ni (u,v) \mapsto F(u,v) \in \V_p(\dom) 
$$ 
defined by
\begin{align*}
F(u,v;\om) &:= {{u + \grad \cdot (D_1 u)}}   - \grad \cdot (g(u{,v}) \mathbb{D}_M \grad v) + {{\mu_M u(1-u)v}}-R_M(\om) u,   
\end{align*}
is a bounded nonlinear operator. More precisely, for \aall $\omega\in\Omega$,  
\begin{align*}
    {{\|F(u,v)\|_{\V_p}}} &{{\le  \Big(1 + \const_{R_M} +\ubd_D + (2\mu_M+\ubd_g\ubd_D)\big\|v\big\|_{\W^{2}_p}\Big) \big\|u\big\|_{\wrp}  +  \Big({2}\ubd_g { {\ubd_D}}+ \ubd_g \ubd_D \big\|v\big\|_{\wtp}\Big) \big\|v\big\|_{\wtp}.}} 
\end{align*}
where $\const_{R_M} := \|R_M\|_{\wsp}$  and $\ubd_{D}, \ubd_g$ are introduced in Assumptions~\ref{asm:model_params}.\ref{Assumpt-3.4.2.} and~\ref{asm:model_params}.\ref{Assumpt-3.4.3.}, respectively. Furthermore if $s \ge 1 + {n \over p} > r \ge 1 + {n \over 2p}$, then for \aall $\om \in \Om$ we obtain
\begin{align*}
   & {{\big\|F(u_1,v_1) - F(u_2,v_2)\big\|_{\V_p}}}\\[0.2cm] 
   &{{\le 
            \Big(\const_{R_M} + \ubd_{D} + 3 \ubd_{D} \ubd_g + \|v_1\|_{\wtp}\big(1 + \mu_M(1 + \|u_1\|_{\wrp}+\|u_2\|_{\wrp})\big)\Big) \big\|u_1 - u_2\big\|_{\wrp}}} \\[0.2cm]
        &{{\hspace*{.5cm}+  \Big(\big\|u_1\big\|_{\wrp} + \ubd_g \ubd_D (1 + 3\big\|v_1\big\|_{\wtp} + \big\|u_2\big\|_{\wrp} + \big\|v_2\big\|_{\wtp}) + 2\mu_M\|u_2\|_{\wrp}\Big) \big\|v_1 - v_2\big\|_{\wtp}\,.}}
\end{align*}
Additionally, for $u, R_M \in \cnt^{\vsig}([0,T];\wrp)$ and $v \in \cnt^{\vsig}([0,T];\wtp)$, the mapping $t \mapsto F(t) \in \V_p(\dom)$ belongs to $\cnt^{\vsig}([0,T];\V_p)$. 
\end{lem}

\begin{proof}\label{lem:proof_Fuv_property}
We first estimate the norm of $F(u,v)$ for given $u \in \W^r_p(\dom)$ and $v \in \W^2_p(\dom)$. Specifically,
{\allowdisplaybreaks
\begin{align*}
    \|F(u,v)\|_{\V_p} &= \|{u} + \grad \cdot {(D_1 u)} - \grad \cdot (g(u{,v}) \mathbb{D}_M \grad v) + \mu _{ M}u\,v -  {\mu_M}u^2v - R_M u\|_{\V_p} \\[0.3cm]
    &\le {{\|u\|_{\V_p}\ }} + \|\grad \cdot  D_1  u + D_1 \grad u\|_{\V_p} + \|\grad \cdot (g(u{,v}) \mathbb{D}_M \grad v)\|_{\V_p} + {\mu_M}\|u v\|_{\V_p} + {\mu_M}\|u^2v\|_{\V_p} + \|R_M u\|_{\V_p} \\[0.2cm]
    &\le {{\|u\|_{\V_p}\ }} + \|\grad \cdot D_1\|_{L_{\infty}} \|u\|_{p} + \|D_1\|_{L_{\infty}} \big\|\sum_{i=1}^n |\partial_{i} u| \big\|_{\V_p} + \|\grad \cdot (g(u{,v}) \mathbb{D}_M \grad v)\|_{\V_p} \\
    &\hspace*{1cm} + {\mu_M}\|uv\|_{\V_p} + {\mu_M}\|u^2v\|_{\V_p} + \|R_M\|_{L_{\infty}} \| u\|_{\V_p} \\[0.3cm]
    &\le {{\|u\|_{\V_p}\ }} + \ubd_{D} \|u\|_{\wop} + \|\grad \cdot (g(u{,v}) \mathbb{D}_M \grad v)\|_{\V_p} + {\mu_M}\|uv\|_{\V_p} + {\mu_M}\|u^2v\|_{\V_p} + \|R_M\|_{L_{\infty}} \| u\|_{\V_p} \\[0.3cm]
    &\le {{(1+\const_{R_M})\|u\|_{\V_p}}} + \ubd_D \|u\|_{\wop} + {\mu_M}\|uv\|_{\V_p} + {\mu_M}\|u^2v\|_{\V_p} \\[0.3cm]
    &\hspace*{1cm} + \|\grad g(u{,v})\cdot  \mathbb{D}_M \grad v \|_{\V_p} + \| g(u{,v}) (\grad \cdot \mathbb{D}_M) \cdot \grad v \|_{\V_p} + \|g(u{,v})  \mathbb{D}_M : \nabla ^2 v \|_{\V_p} \\[0.2cm]
    &\le {{(1+\const_{R_M})\|u\|_{\V_p}}} +  \ubd_D \|u\|_{\wop} + {\mu_M}\|uv\|_{\V_p} + {\mu_M}\|u^2v\|_{\V_p} {+\|g_u\|_{L_{\infty}}\ubd_D\|\grad u \cdot \grad v \|_{\V_p}} \\[0.2cm]
 &\hspace*{1cm} {+\|g_v\|_{L_{\infty}}\ubd_D\| |\grad v|^2 \|_{\V_p}+\ubd_D\|\nabla^{\kappa}g\cdot\grad v\|_{\V_p}+ \ubd_g  {\ubd_D} \big\| \sum_i |\partial_i v| \big \|_{\V_p}}  {+ \ubd_g \ubd_D \big \|\sum \limits _{i,j}|\partial^2_{ij}v| \big\|_{\V_p}}  \\[0.15cm]
 &\le {{(1+\const_{R_M})\|u\|_{\V_p}}} + \ubd_D \|u\|_{\wop} + {\mu_M}\|uv\|_{\V_p} + {\mu_M}\|u^2v\|_{\V_p} {+ \ubd_g \ubd_D\|\grad u \cdot \grad v \|_{\V_p}} \\[0.2cm]
 &\hspace*{1cm} + \ubd_g \ubd_{D}\| \grad v \|^2_{\V_{2p}}+ \ubd_g\ubd_D\|\grad v\|_{\V_{2p}}+ \ubd_g {\ubd_D} \big\| \sum_i |\partial_i v| \big \|_{\V_p}  + {\ubd_g \ubd_D \big \|\sum \limits _{i,j}|\partial^2_{ij}v| \big\|_{\V_p}},
\end{align*}
}
where we used $\grad g=\partial_u g\,\grad u+\partial_v g\,\grad v+\nabla^{\kappa}g$ and $\|\nabla^{\kappa}g\|_{L_{\infty}}\le\ubd_g$ from Assumption~\ref{asm:model_params}.\ref{Assumpt-3.4.3.}, and $\ubd_D$ as in Assumption~\ref{asm:model_params}.\ref{Assumpt-3.4.2.}. Since $\W^{1+r_1}_p \hookrightarrow \W^1_{2p}(\dom)$ for $r_1 > \frac{n}{2p}$, taking $r \in \left(1+\frac{n}{2p}, 2\right]$ we get that
\begin{align*}
    \big\|\grad u \cdot \grad v\big\|_{\V_p} &\le \big\|\grad u\big\|_{\V_{2p}} \big\|\grad v\big\|_{\V_{2p}} \le \big\|u\big\|_{\W^{1}_{2p}} \big\|v\big\|_{\W^{1}_{2p}} \le \big\|u\big\|_{\W^{r}_p} \big\|v\big\|_{\W^{r}_p} \\[0.2cm]
    \big\|u^2 v\big\|_{\V_p} &\le \big\|u\big\|_{\V_{2p}} \big\|v\big\|_{L_{\infty}} \le \big\|u\big\|_{\W^{1}_{2p}} \big\|v\big\|_{\W^{1}_{2p}} \le \big\|u\big\|_{\W^{r}_p} \big\|v\big\|_{\W^{r}_p}.
\end{align*}
Therefore, we obtain
\begin{align*}
\big\|F(u,v)\big\|_{\V_p} &\le {{(1 + \const_{R_M})\big\|u\big\|_{\V_p}}} + \ubd_{D} \big\|u\big\|_{\wop} + {\mu_M}\|u\|_{\W^r_p}\|v\|_{\W^2_p} + {\mu_M}\|u\|_{\W^r_p}\|v\|_{\W^r_p} \\[0.2cm]
 &\hspace*{1cm} {+ \ubd_g \ubd_D \big\|u\big\|_{\W^{r}_p} \big\|v \big\|_{\W^{r}_p}} {+ \ubd_g \ubd_D \big\|v \big\|^2_{\W^{r}_p}+ \ubd_g { {\ubd_D}} \Big\| \sum_i |\partial_i v| \Big\|_{\V_p}}  + {\ubd_g \ubd_D \Big\|\sum \limits _{i,j}|\partial^2_{ij}v| \Big\|_{\V_p}}\\[0.3cm]
 &\le {{(1 + \const_{R_M} )\big\|u\big\|_{\V_p}}} + (\mu_M + {2}\ubd_g {\ubd_D}) \big\|v \big\|_{\wtp}  + \ubd_D \big\|u\big\|_{\wop} {+ (\mu_M + \ubd_g \ubd_D) \big\|u\big\|_{\W^{r}_p} \big\|v \big\|_{\W^{2}_p}} \\[0.2cm]
 &\hspace*{1cm} {+ \ubd_g \ubd_D \big\|v \big\|^2_{\W^{2}_p}} {+ \mu_M\big\|u\big\|_{\W^{r}_p} \big\|v\big\|_{\W^{2}_p}} \\[0.3cm]
&\le  \Big(1 + \const_{R_M} +\ubd_D + (2\mu_M+\ubd_g\ubd_D)\big\|v\big\|_{\W^{2}_p}\Big) \big\|u\big\|_{\wrp}  +  \Big({2}\ubd_g { {\ubd_D}}+ \ubd_g \ubd_D \big\|v\big\|_{\wtp}\Big) \big\|v\big\|_{\wtp}. 
\end{align*}
Now we estimate the Lipschitz continuity of $F(u,v)$ for $u_1, u_2 \in \W^s_p(\dom)$, with $s \ge 1 + {n \over p}$ and ${v_1, v_2 \in \W^2_p(\dom)}$. Accordingly, we have that
\begin{align}
    \label{eq:F_diff_pre_estimate}
    &\big\|F(u_1,v_1) - F(u_2,v_2)\big\|_{\V_p} \\[0.2cm]
    &\le \big\|u_1 - u_2\big\|_{\V_p} + \big\|\grad \cdot (D_1 (u_1 - u_2))\big\|_{\V_p} + \big\|\grad \cdot ((g(u_1,v_1) - g(u_2,v_2)) \mathbb{D}_M \grad v_1)\big\|_{\V_p}\notag\\[0.2cm]
    &\hspace*{1cm} + \big\|\grad \cdot (g(u_2,v_2) \mathbb{D}_M \grad (v_1 - v_2))\big\|_{\V_p} + {\mu_M} \big\|u_1 v_1 - u_2 v_2\big\|_{\V_p} + {\mu_M} \big\|u_1^2 v_1 - u_2^2 v_2\big\|_{\V_p} + R_M \big\|u_1 - u_2\big\|_{\V_p} \notag\\[0.3cm]
    &\le  (1 + \const_{R_M} + \ubd_{D}){{\big\|u_1 - u_2\big\|_{\wop}}} + \big\|\grad \cdot ((g(u_1,v_1) - g(u_2,v_2)) \mathbb{D}_M \grad v_1)\big\|_{\V_p}  \notag\\[0.2cm]
    &\hspace*{1cm} + \big\|\grad \cdot (g(u_2,v_2) \mathbb{D}_M \grad (v_1 - v_2))\big\|_{\V_p} + {\mu_M} \big\|u_1 v_1 - u_2 v_2\big\|_{\V_p} + {\mu_M} \big\|u_1^2 v_1 - u_2^2 v_2\big\|_{\V_p} .\notag
\end{align}
Moreover, we have
{\allowdisplaybreaks
\begin{align}\label{eq:grad_g_term_diff_estimate_1}
    &\big\|\grad \cdot ((g(u_1,v_1) - g(u_2,v_2)) \mathbb{D}_M \grad v_1)\big\|_{\V_p} \notag\\[0.3cm]
    &= \Big\|\grad (g(u_1,v_1) - g(u_2,v_2)) \cdot \mathbb{D}_M \grad v_1 + (g(u_1,v_1) - g(u_2,v_2)) (\grad \cdot \mathbb{D}_M) \cdot \grad v_1 \notag\\[0.2cm]
    &\hspace*{2.5cm}+ (g(u_1,v_1) - g(u_2,v_2))  \mathbb{D}_M : \nabla ^2 v_1 \Big\|_{\V_p} \notag\\[0.3cm]
    &\le \big\|\grad (g(u_1,v_1) - g(u_2,v_2)) \cdot \mathbb{D}_M \grad v_1\big\|_{\V_p} + \big\|(g(u_1,v_1) - g(u_2,v_2)) (\grad \cdot \mathbb{D}_M) \cdot \grad v_1\big\|_{\V_p} \notag\\[0.2cm]
    &\hspace*{2.5cm} + \big\|(g(u_1,v_1) - g(u_2,v_2))  \mathbb{D}_M : \nabla ^2 v_1\big\|_{\V_p} \notag\\[0.3cm]
    &\le \big\|\mathbb{D}_M \big\|_{\infty} \big\|\grad (g(u_1,v_1) - g(u_2,v_2))\big\|_{\V_{2p}} \big\|\grad v_1\big\|_{\V_{2p}}+ \big\|g(u_1,v_1) - g(u_2,v_2)\big\|_{L_{\infty}} \big\|\grad \cdot \mathbb{D}_M\big\|_{L_{\infty}} \big\|\grad v_1\big\|_{\V_p} \notag\\[0.2cm]
    &\hspace*{2.5cm} + \big\|g(u_1,v_1) - g(u_2,v_2)\big\|_{L_{\infty}} \big\|\mathbb{D}_M\big\|_{L_{\infty}} \big\|\nabla^2 v_1\big\|_{\V_p} \notag\\[0.3cm]
    &\le \ubd_{D} \Big(\big\|\grad (g(u_1,v_1) - g(u_2,v_2))\big\|_{\V_{2p}} \big\|\grad v_1\big\|_{\V_{2p}} + \big\|g(u_1,v_1) - g(u_2,v_2)\big\|_{L_{\infty}} \big\|\grad v_1\big\|_{\V_p} \\[0.2cm]
    &\hspace*{2.5cm} + \big\|g(u_1,v_1) - g(u_2,v_2)\big\|_{L_{\infty}}  \big\|\nabla^2 v_1\big\|_{\V_p} \Big)\,. \notag
\end{align}}
\noindent For the first component in the first term of the above estimate, the local Lipschitz continuity of the Nemytskii map $(u,v)\mapsto g(u,v;\cdot)$ in $\W^{r}_p$ (Assumption~\ref{asm:model_params}.\ref{Assumpt-3.4.3.}, which already accounts for the spatial coefficients $\kappa_1,\kappa_2$) yields
\begin{align}\label{eq:grad_g_diff_estimate}
\big\|\grad (g(u_1,v_1) - g(u_2,v_2))\big\|_{\V_{2p}} 
&\le C\big\|g(u_1,v_1) - g(u_2,v_2)\big\|_{\W^{1}_{2p}} \notag\\[0.2cm]
&\le \ubd_g \Big( \big\|u_1 - u_2\big\|_{\W^r_p} + \big\|v_1 - v_2\big\|_{\W^2_p} \Big)\,.
\end{align}
The second and third terms in~\eqref{eq:grad_g_term_diff_estimate_1} can be estimated as follows 
\begin{align}\label{eq:g_diff_estimate}
\big\|g(u_1,v_1) - g(u_2,v_2)\big\|_{L_{\infty}} &\le \big\|g_u\big\|_{L_{\infty}} \big\|u_1 - u_2\big\|_{L_{\infty}} + \big\|g_v\big\|_{L_{\infty}} \big\|v_1 - v_2\big\|_{L_{\infty}} \\[0.2cm]
&\le \ubd_g \Big( \big\|u_1 - u_2\big\|_{\W^r_p} + \big\|v_1 - v_2\big\|_{\W^2_p} \Big)\,.
\end{align}
Plugging \eqref{eq:grad_g_diff_estimate} and \eqref{eq:g_diff_estimate} into \eqref{eq:grad_g_term_diff_estimate_1} we get that
\begin{align}\label{eq:grad_g_term_diff_estimate_2}
   & \big\|\grad \cdot ((g(u_1,v_1) - g(u_2,v_2)) \mathbb{D}_M \grad v_1)\big\|_{\V_p} \notag \\[0.2cm]
    &\le \ubd_{D} \ubd_g  \Big( \big\|\grad v_1\big\|_{\V_{2p}} + \big\|\grad v_1\big\|_{\V_p} + \big\|\nabla^2 v_1\big\|_{\V_p} \Big) \Big( \big\|u_1 - u_2\big\|_{\W^r_p} + \big\|v_1 - v_2\big\|_{\W^2_p} \Big) \notag\\[0.2cm]
&\le 3 \ubd_{D} \ubd_g  \big\|v_1\big\|_{\wtp} \Big( \big\|u_1 - u_2\big\|_{\W^r_p} + \big\|v_1 - v_2\big\|_{\W^2_p} \Big).
\end{align}
\noindent Similarly,
{\allowdisplaybreaks
\begin{align}\label{eq:grad_g_term_diff_estimate_3}
&\big\|\grad \cdot (g(u_2,v_2) \mathbb{D}_M \grad (v_1 - v_2))\big\|_{\V_p}\\[0.2cm] 
&\le \big\|\grad g(u_2,v_2) \cdot \mathbb{D}_M \grad (v_1 - v_2)\big\|_{\V_p} + \big\|g(u_2,v_2) (\grad \cdot \mathbb{D}_M) \cdot \grad (v_1 - v_2)\big\|_{\V_p}  + \big\|g(u_2,v_2)  \mathbb{D}_M : \nabla ^2 (v_1 - v_2)\big\|_{\V_p} \notag\\[0.2cm]
&\le \big\|\grad g(u_2,v_2)\big\|_{\V_{2p}} \big\|\mathbb{D}_M\big\|_{L_{\infty}} \big\|\grad (v_1 - v_2)\big\|_{\V_{2p}} + \big\|g(u_2,v_2)\big\|_{L_{\infty}} \big\|\grad \cdot \mathbb{D}_M\big\|_{L_{\infty}} \big\|\grad (v_1 - v_2)\big\|_{\V_p} \notag\\[0.15cm]
&\hspace*{2.5cm} + \big\|g(u_2,v_2)\big\|_{L_{\infty}} \big\|\mathbb{D}_M\big\|_{L_{\infty}} \big\|\nabla^2 (v_1 - v_2)\big\|_{\V_p} \notag\\[0.2cm]
&\le \big(\big\|g_u\big\|_{L_{\infty}} \big\|\grad u_2\big\|_{\V_{2p}} + \big\|g_v\big\|_{L_{\infty}} \big\|\grad v_2\big\|_{\V_{2p}} + \|\nabla^{\kappa}g(u_2,v_2)\|_{L_{2p}}\big) \big\|\mathbb{D}_M\big\|_{L_{\infty}} \big\|\grad (v_1 - v_2)\big\|_{\V_{2p}}\notag \\[0.15cm]  
&\hspace*{1cm}+ \big\|g(u_2,v_2)\big\|_{L_{\infty}} \big\|\grad \cdot \mathbb{D}_M\big\|_{L_{\infty}} \big\|\grad (v_1 - v_2)\big\|_{\V_p} + \big\|g(u_2,v_2)\big\|_{L_{\infty}} \big\|\mathbb{D}_M\big\|_{L_{\infty}} \big\|\nabla^2 (v_1 - v_2)\big\|_{\V_p} \notag\\[0.2cm]
&\le \ubd_g \ubd_D \big(1 + \big\|u_2\big\|_{\W^1_{2p}} + \big\|v_2\big\|_{\W^1_{2p}}\big) \big\|\grad (v_1 - v_2)\big\|_{\W^1_{2p}} + \ubd_g \ubd_D \big\|\grad (v_1 - v_2)\big\|_{\V_p} + \ubd_g \ubd_D \big\|\nabla^2 (v_1 - v_2)\big\|_{\V_p} \notag\\[0.2cm]
&\le \ubd_g \ubd_D \Big(1 + ( \big\|u_2\big\|_{\wrp} + \big\|v_2\big\|_{\wtp})\Big) \big\|v_1 - v_2\big\|_{\W^2_p}\,,
\end{align}}
where $\|\nabla^{\kappa}g\|_{L_{2p}}\le\ubd_g$ follows from Assumption~\ref{asm:model_params}.\ref{Assumpt-3.4.3.} ($\|\nabla^{\kappa}g\|_{L_{\infty}}\le\ubd_g$).
Thus, plugging \eqref{eq:grad_g_term_diff_estimate_2} and \eqref{eq:grad_g_term_diff_estimate_3} into \eqref{eq:F_diff_pre_estimate} and using 

\begin{align*}
\|u_1 v_1 - u_2 v_2\|_{\V_p} &\le \|v_1\|_{\wtp}\|u_1-u_2\|_{\wrp} + \|u_2\|_{\wrp}\|v_1-v_2\|_{\wtp},\\
\|u_1^2 v_1 - u_2^2 v_2\|_{\V_p} &\le \|v_1\|_{\wtp}(\|u_1\|_{\wrp}+\|u_2\|_{\wrp})\|u_1-u_2\|_{\wrp} + \|u_2\|_{\wrp}\|v_1-v_2\|_{\wtp},
\end{align*}
we get that
\begin{align*}
    &\big\|F(u_1,v_1) - F(u_2,v_2)\big\|_{\V_p}\\[0.2cm] 
    &\le  \Big(\const_{R_M} + \ubd_{D} + 3 \ubd_{D} \ubd_g + \|v_1\|_{\wtp}\big(1 + \mu_M(1 + \|u_1\|_{\wrp}+\|u_2\|_{\wrp})\big)\Big) \big\|u_1 - u_2\big\|_{\wrp} \\[0.2cm]
    &\,\,\,\,\,\,\,\,\,+\Big( \big\|u_1\big\|_{\wrp} + \ubd_g \ubd_D (1 + 3\big\|v_1\big\|_{\wtp} + \big\|u_2\big\|_{\wrp} + \big\|v_2\big\|_{\wtp}) + 2\mu_M\|u_2\|_{\wrp}\Big) \big\|v_1 - v_2\big\|_{\wtp} .
\end{align*}
This completes the proof of the Lipschitz continuity of $F(u,v)$.
\end{proof}

\noindent If $v$ is regarded as a given (known) function, then $F$ can be viewed as an operator depending only on $u$. In this case, the previous lemma directly yields the following corollary.

\begin{cor}[Properties of {$F(\cdot ,v)$}]\label{cor:Fu_property}
Let $v \in \cnt^{\vsig}([0,T];\W^2_p(\dom))$ be a given function satisfying $\|v\|_{\cnt^{\vsig}} \le \const_v$. Then, for $s \in \left(1 + {n \over p}, \infty \right)$, $s_1 > s > r > 1 + {n \over 2p}$, and for fixed $\om \in \Om$, the mapping 

$$
\W^{s_1}_p(\dom) \ni u \mapsto F(u,{v; \om }) \in \V_p(\dom), 
$$ 
where 
\begin{align*}
F(u,{v;\om}):= {{u + \grad \cdot (D_1 u)}}   - \grad \cdot (g(u{,v}) \mathbb{D}_M \grad v) + {{\mu_M u(1-u)v}}-R_M(\om) u,  
\end{align*}
is a bounded nonlinear operator. Moreover,
\begin{align*}
    {{\|F(u,v)\|_{\V_p}}} &{{\le  \Big(1 + \const_{R_M} +\ubd_D + (2\mu_M+\ubd_g\ubd_D) \const_v \Big) \big\|u\big\|_{\W^{s_1}_p}  +  \Big({2}\ubd_g { {\ubd_D}}+ \ubd_g \ubd_D \const_v\Big) \const_v}} 
\end{align*}
for \aall $\om \in \Om$. Furthermore, it holds that
\begin{align}
\label{eq:Lip_Fu_op}
\begin{aligned}
{{\big\|F(u_1, v) - F(u_2, v)\big\|_{\V_p}}} &{{\le  \Big(\const_{R_M}\!\! + \const_v\!\! +  \ubd_{D}(1 + 3\ubd_g)\!\! + \mu_M\const_v\big(\|u_1\|_{\W^{r}_p}+\|u_2\|_{\W^{r}_p}\big)\Big) \big\|u_1 - u_2\big\|_{\W^{r}_p}}} \\
    &{{\le \Big(\const_{R_M}\!\! + \const_v\!\! + \ubd_{D}(1 + 3\ubd_g)\!\! + \mu_M\const_v\big(\|u_1\|_{\W^{r}_p}+\|u_2\|_{\W^{r}_p}\big)\Big) \big\|u_1 - u_2\big\|_{\W^{s_1}_p}.}}
\end{aligned}
\end{align}
\end{cor}

\noindent The results established above provide the required properties to state and prove one of the main theorems regarding the existence and regularity of the solution to the cancer cell density equation \eqref{eq:ACP_Can}.

\begin{thm}[Solution to cancer equation]\label{thm:U_Soln}
Assume that Assumption~\ref{asm:model_params} is satisfied. Let $r, s, s_1, s_2 \in \R^+$ be such that 
\begin{align}
    \label{eq:u_soln_dom_expos}
0 < r < s < s_1 < s_2  = 2 ~~\text{\rm with }  r > 1 + {n \over 2p}, \: s > 1 + {n \over p}, \:\: \vartheta \in \left({s_1 \over 2}, 1 \right].
\end{align}
Let intermediate exponents $r',s',s_1'$ satisfy $0<r<r'<s<s'<s_1<s_1'<s_2=2$ and $s_1'\in(s_1,\min\{s_2,2\vartheta\})$, and set $\alpha:=r'/2$, $\beta:=s'/2$, $\nu:=\vartheta$, and $\eta:=s_1'/2$. 
Assume moreover that the time H\"older exponent obeys
\begin{equation}\label{eq:U_Soln_vsig_Yagi}
\vsig < \min\Bigl\{\tfrac1p,\tfrac12,\,\beta-\alpha+\nu-1,\,1-\eta\Bigr\}.
\end{equation}
Let $u_0 \in \W^{2}_p(\dom)$. Then, for any given mappings
\begin{align}
    \label{eq:u_soln_v_w_premise}
    &w \in \B([0,T];\wsp)\cap\cnt^{\vsig}([0,T];\wrp) ~ \text{\rm with } \|w\|_{\B([0,T];\wsp)}+\|w\|_{\cnt^{\vsig}([0,T];\wrp)} \le \const_w,\\
    &v \in \cnt^{\vsig}([0,T];\W^{2}_p) ~ \text{\rm with } \|v\|_{\cnt^{\vsig}([0,T];\W^2_p)} \le \const_v,
\end{align}
the problem~\eqref{eq:ACP_Can} has a unique local solution $u$ on the time interval $[0,T]$ for \aall $\om \in \Om$. Furthermore, for \aall $\om \in \Om$, $u(\om)$ satisfies 
\begin{align}\label{eq:loc_sln_spc}
    &u \in \mathcal{C}([0,T];\W^s_p) \cap \mathcal{C}^{\vrho}([0,T];\W^{r}_p) \cap \mathcal{C}^{1}([0,T];L_p), \\[0.15cm]
    &A(w,v) u \in \cnt([0,T];\V_p); \\[0.15cm]
    &A(w,v) u \in \Y^{1,\vsig}((0,T];\V_p)\,.
\end{align}
\end{thm}
\begin{proof}


For prescribed functions $w,v$ satisfying~\eqref{eq:u_soln_v_w_premise}, the cancer equation~\eqref{eq:ACP_Can}-\eqref{eq:ACP_Ops} takes the form

$$
{d \over dt} u(t) + A(w(t),v(t))\, u(t) = F(u(t),v(t)) \quad \text{in } \V_p, \qquad u(0) = u_0,
$$
with $A(w,v)$ and $F(u,v)$ as in~\eqref{eq:ACP_Ops}. Since the coefficient paths $w$ and $v$ are given (frozen), the linear part $A(t):=A(w(t),v(t))$ is a prescribed time-dependent family of sectorial operators on $\V_p$, while $u\mapsto F(u,v(t))$ is a semilinear perturbation. 

Now the idea is to apply Theorem 5.6 of~\cite{Yagi09}. Although the theorem is stated for a quasilinear problem, it applies here as a degenerate special case: freezing $(w,v)$ makes the map $\tilde u\mapsto A(w(t),v(t))$ independent of the unknown $\tilde u$, while the time dependence of the sectorial family is carried entirely by the prescribed H\"older paths. Equivalently, one may view~\eqref{eq:ACP_Can} as a non-autonomous semilinear problem whose structural hypotheses coincide with those of Theorem~5.6 along the frozen coefficients. 

Based on this we proceed to verify the assumptions of Theorem~5.6. Following the notation in~\cite{Yagi09}, set $X:=\V_p$, $Y:=\W^{r}_p$, $Z:=\W^{s}_p$, $W:=\W^{s_1}_p$, and

$$
\nu := \vartheta \le 1,\qquad
\gamma := {s_2 \over 2} = 1,\qquad
\beta := {s^{\prime} \over 2},\qquad
\alpha := {r^{\prime} \over 2},\qquad
\eta := {s_1^{\prime} \over 2},
$$
where the intermediate Sobolev exponents $r^{\prime},s^{\prime},s_1^{\prime}$ are chosen so that

$$
0 < r < r^{\prime} < s < s^{\prime} < s_1 < s_1^{\prime} < s_2 = 2
\quad\text{and}\quad
s_1^{\prime}\in\bigl(s_1,\min\{s_2,2\vartheta\}\bigr).
$$
Consequently, the following embeddings hold:
\begin{align} \label{eq:U_Soln_embds}
\W^{s_2}_{p,N} \subset \D(A^{\eta}_p) \subset \W^{s_1}_{p,N} \subset \D(A^{\beta}_p) \subset \W^{s}_{p,N} \subset \D(A^{\alpha}_p) \subset \W^{r}_{p,N} \subset \V_p \quad \text{(cf. (1.70), (16.36) in~\cite{Yagi09})} \\[0.2cm]
\W^{s_2}_{p,N} \hookrightarrow \W^{s_1}_{p,N} \hookrightarrow \W^{s}_{p,N} \hookrightarrow \W^{r}_{p,N}  \quad \text{(cf.  (1.70), (16.36) in~\cite{Yagi09}).}
\end{align}

\noindent For $s > 1 + {n \over p} > r > 1 + {n \over 2p}$ and $\vartheta \le 1$, Lemma~\ref{lem:Au_property} yields the sectoriality, uniform resolvent bound, and domain comparison for $A(w,v)$, which are precisely conditions~(5.2)--(5.5) of~\cite{Yagi09} for the coefficient-parameterized family (with~\eqref{eq:Au_struc_property} giving the Lipschitz estimate in the $Y$-norm of the coefficients). Since $A$ is independent of the unknown $u$, the $u$-Lipschitz constant in~(5.5) may be taken as zero; the required time regularity of $t\mapsto A(t)^{-1}$ then follows from Corollary~\ref{cor:Au_Holder} and the coefficient regularity~\eqref{eq:u_soln_v_w_premise}. Furthermore, $r < r^{\prime}=2\alpha$ and $s < s^{\prime}=2\beta$ give the embeddings~(5.6). Condition~(5.7) is the fractional-power relation

$$
0\le \alpha < \beta < \nu=\vartheta\le 1
\quad\text{and}\quad
1+\alpha < \beta+\nu
$$
which is compatible with~\eqref{eq:u_soln_dom_expos}; e.g., take $\vartheta$ close to~$1$, $s^{\prime}=1+n/p+\epsilon$, and $r^{\prime}=1+n/(2p)+\epsilon$ for small $\epsilon>0$ and $p>n$. Thus, the structural conditions for non-autonomous analytic evolution operators are satisfied.

On the other hand, since $\W^{s_1}_p \hookrightarrow \W^{s}_p$, for $s_1 > s > 1 + {n \over p}$ Corollary~\ref{cor:Fu_property} yields the Lipschitz estimate~\eqref{eq:Lip_Fu_op} for the nonlinear mapping $u\mapsto F(u,v)$ with frozen path~$v$, which is condition~(5.44) of~\cite{Yagi09} (taking $G\equiv 0$, or absorbing the explicit H\"older time dependence of $v$ and $R_M$ into the moduli of continuity of the data). Moreover, $s_1 < s_1^{\prime}=2\eta$ gives $\D(A^{\eta}_p)\subset\W^{s_1}_p$ (condition~(5.46)), and the choice $s_1^{\prime}\in(s_1,\min\{s_2,2\vartheta\})$ ensures $\beta<\eta<\nu$ (condition~(5.47)).

Finally, $u_0\in\W^{2}_p(\dom)\subset\D(A^{\gamma})$ with $\gamma=1$ (cf.~\eqref{eq:dom_embds}) and $R_M\in\cnt^{\vsig}([0,T];\wsp)$ (Lemma~\ref{rad_asure_bded}) give the regular data~(5.29)--(5.30), under the H\"older restriction~\eqref{eq:U_Soln_vsig_Yagi} already imposed in the theorem statement (which implies in particular $\vsig<\min\{\beta-\alpha+\nu-1,\,1-\eta\}$).
Theorem~5.6 of~\cite{Yagi09} therefore applies and yields the asserted existence, uniqueness, and regularity, with H\"older exponent $\vrho=\gamma-\alpha=1-\alpha$ in $\W^{r}_p$. 
\end{proof}

Before proceeding with the construction of local solutions, we first outline the strategy. The apriori estimates and regularity results established in the previous subsection, namely the positivity and uniform bounds for $v$, the pathwise H\"older continuity and moment estimates for the radiotherapy process, the Lipschitz estimates of the nonlinearities, and the uniform ellipticity and smoothness of the diffusion tensor, ensure that, for almost every realization of the noise, the coefficients of the tumor equation are spatially regular, uniformly bounded, and H\"older continuous in time. This justifies the reformulation of the tumor PDE as a non‑autonomous abstract Cauchy problem in $L_p$ driven by a family of sectorial
operators $A(w_t,v_t)$ for given mappings $w$ and $v$ as mentioned above. 

The next series of lemmas establishes the operator norm, resolvent, and H\"older continuity estimates needed to construct the associated evolution operator and to close a fixed‑point (contraction) argument. This framework yields local (pathwise) solutions that are compatible with the biological assumptions and stochastic structure developed above. As a preliminary step, we first collect several technical results that enable the definition of an approximate solution mapping for the cancer cell density equation~\eqref{eq:ACP_Can}.
 
\begin{lem}[Almost sure boundedness]\label{lem:asure_bded}
Let $w \in \K^p_{\loc}(\Zb)$. Then $\big\|w\big\|_{\B([0,T];\wrp)}<\infty$ for a.a.\ $\om\in\Om$. In particular, there exists a $(\Fil)$-measurable finite random variable $\ubd_w:\Om\to[0,\infty)$ such that $|w(t,x)|\le \ubd_w$ for a.e.\ $x\in\dom$, all $t\in[0,T]$, and a.a.\ $\om\in\Om$.
\end{lem}

\begin{proof}
Since $w \in \K^p_{\loc}(\Zb)$, we have $\E \bigl[\big\|w\big\|^p_{\zrp} \bigr] \le \const^p_{\Zb}<\infty$, hence $\big\|w\big\|_{\B([0,T];\wrp)}<\infty$ for a.a.\ $\om\in\Om$.
Because $r>n/p$, the embedding $\wrp\hookrightarrow L_{\infty}(\dom)$ yields, for a.a.\ $\om$,
\begin{align*}
\sup{t \in [0,T]} \underset{x \in \dom}{\rm ess\,sup}\, |w(t,x)|
= \sup{t \in [0,T]} \big\|w(t)\big\|_{L_{\infty}(\dom)}
\le C_{\mathrm{emb}}\,\big\|w\big\|_{\B([0,T];\wrp)} < \infty.
\end{align*}
Taking $\ubd_w(\om):=C_{\mathrm{emb}}\,\big\|w(\om)\big\|_{\B([0,T];\wrp)}$ (and $\ubd_w:=0$ on the null set where the norm is infinite) gives the claim.
\end{proof}

\newcommand{\kploc}{\K^p_{\loc}(\Zb)}

Based on the above lemma, we can define an approximate solution mapping for \eqref{eq:ACP_Can} as follows.
\begin{lem}[Approximate solution mapping]\label{lem:sem_lin_aprx_soln}
Let $w \in \K^p_{\loc}(\Zb)$ and let $v(w)$ denote the solution of~\eqref{eq:pre_normal}. Then, the problem 
{\allowdisplaybreaks
\begin{align}\label{eq:explicit_pre_cancer}
\partial_t u-\grad \cdot  \left(D_1 u+ D_0(w,v(w)) \grad u \right)+\grad \cdot\left(g(u,v(w)) \mathbb{D}_M \grad v(w) \right)&=\big[\mu_M(1-u) v(w) - R_M(\xi_t)\big] u \\[0.2cm]
\left(D_1 u + D_0(w,v(w)) \grad u\right) \cdot {\nu} = \left(g(u,v(w)) \mathbb{D}_M \grad v(w)\right) \cdot  {\nu }&= 0 \\[0.2cm]
u(0) &= u_0 \\
v(w) = v_0 \exp \Big(-\int_0^t \big(\kappa_3 \psi(w_s) &+ R_Q(\xi_s) \big)  \: ds \Big)
\end{align}}
admits a unique solution $u(w)$. Consequently, the solution $u(w)$ is represented as a mapping $w \mapsto \Phi(w)$ with $\Phi :\K^p_{\loc}(\Zb) \to \X(T)$.
\end{lem}

\begin{proof} 
The proof is divided into three steps.
\begin{enumerate}
    \item Since $w \in \kploc$, we can apply Lemma~\ref{lem:asure_bded} to obtain that $w \in \ztp(T)\cap\zrp(T)$ almost surely, i.e., $w(\om) \in \zrp(T)$ for $\aall \om \in \Om$.
    \item With $w(\om) \in \zrp(T)$ we apply Lemma \ref{lem:sol_normal} to obtain $\ztp(T) \ni w \mapsto v(w) \in \cntvtp$ as a solution to \eqref{eq:pre_normal}. In particular, we obtain the pathwise mapping $\om \mapsto v(w(\om))$. Furthermore, since $\ubd_w(\om)<\infty$ for a.a.\ $\om\in\Om$ by Lemma~\ref{lem:asure_bded}, the estimate~\eqref{eq:v_cvsig_h2p} yields a finite (path-dependent) bound
    
    $$ 
    \const_v(\om) := \big\|v(\om)\big\|_{\cntvtp} < \infty \quad \text{for a.a.\ $\om\in\Om$.}
    $$
    \item Therefore, all assumptions of Theorem~\ref{thm:U_Soln} are satisfied pathwise with the given coefficient paths $w$ and $v=v(w)$ (taking $\const_w(\om):=\|w(\om)\|_{\B([0,T];\wsp)}+\|w(\om)\|_{\cnt^{\vsig}_{r,p}}$, which is finite a.s.\ since $w\in\K^p_{\loc}(\Zb)\subset\X(T)$, and $\const_v(\om)$ as above). Hence, for a.a.\ $\om\in\Om$, problem~\eqref{eq:explicit_pre_cancer} (equivalently,~\eqref{eq:ACP_Can} with $A(w,v(w))$ and $F(u,v(w))$) admits a unique solution  $u(\om) \in \bU(T)$, where
    
    $$ 
    \bU(T) := \left \{u \in \ztp :  u \in \mathcal{C}([0,T];\W^s_p) \cap \mathcal{C}^{\vrho}([0,T];\W^{r}_p) \cap \mathcal{C}^{1}([0,T];\V_p) \right\}. 
    $$
    Since $u(\om) \in \ztp \cap \cntrrp$ we also get that $\E\left[\|u\|^p_{\ztp}\right] + \E\left[\|u\|^p_{\cntrrp}\right] < \infty$, thus, $u \in \X(T).$
\end{enumerate}
Consequently, we may define the mapping $w \mapsto \Phi(w) := u(w)$ with $\Phi:\K^p_{\loc}(\Zb) \to \X(T).$
\end{proof}

We next prove that $\Phi$ is a contraction operator, which will imply the existence of a unique local solution. The following technical estimate is the key ingredient. 

\begin{rem}[Dependence of $A_t$ on the input $w$]\label{rem:A_of_w}
For $w\in\K^p_{\loc}(\Zb)$, write $v(w)$ for the normal-tissue component in Lemma~\ref{lem:sem_lin_aprx_soln} and $u(w)$ for the corresponding cancer solution. We use the shorthand
\[
A[w]_t \;:=\; A(w_t,\,v(w)_t),
\qquad
F[w]_t \;:=\; F\big(u(w)_t,\,v(w)_t\big),
\]
so that the evolution operators $\Tst^i$ in Lemma~\ref{lem:Op_Sg_Lip_ineq} are generated by the families $A[w^i]_t$, $i\in\{1,2\}$. In particular, $A[w]$ depends on the auxiliary input $w$ through both the cancer coefficient $w_t$ and the normal-tissue path $v(w)_t$ in $D_0(w,v)$. At the fixed point $w=\Phi(w)$ one has $u(w)=w$ and $v(w)=v(u)$.

\end{rem}

\begin{lem}[Lipschitz inequality]\label{lem:Op_Sg_Lip_ineq} 
Let $w^1, w^2 \in \K^p_{\loc}(\Zb)$, let $u^i:=u(w^i)$ and $v^i:=v(w^i)$ be as in Remark~\ref{rem:A_of_w}, and set $A^i_t := A[w^i]_t=A(w^i_t,v^i_t)$, $F^i:=F[w^i]=F(u^i,v^i)$. Let $U_0 \in \wtp$.
Furthermore, let $\Tst^1$ and $\Tst^2$ denote the evolution operators generated by the families $A^1_t$ and $A^2_t$, respectively (cf.\ Theorem~\ref{thm:non_auto_evol_op} in Section~\ref{sec:SG_Prop}). Then, for $0 \leq \alpha < \vartheta = 1$ and $0 \leq t \leq T$ the following estimates hold:
{\allowdisplaybreaks
\begin{align}\label{eq:Aalp_T12_U0_lip}
\left\|(A^1_t)^{\alpha}\left[\Tzt^1-\Tzt^2\right] U_0\right\|_{\V_p}\leq C_\alpha t^{1-\alpha+\vrho}\left\|U_0\right\|_{\wtp}\left\|w^1-w^2\right\|_{\cnt^\vrho_{r,p}(T)},  
\end{align}
}
\begin{align}\label{eq:Aalp_T12_F_lip}
\left\|(A^1_t)^\alpha \int_0^t\left[\Tst^1-\Tst^2\right] F(s) ~ ds\right\|_{\V_p}\leq C_\alpha t^{1-\alpha+\vrho}\|F\|_{\Y^{1,\vsig}}\left\|w^1-w^2\right\|_{\cnt^\vrho_{r,p}(T)},   
\end{align}
and
\begin{align}\label{eq:Aalp_T2_F12_lip}
\left\|(A^2_t)^\alpha \int_0^t \Tst^2\left[F^1-F^2\right] d s\right\|_{\V_p}  \leq C_\alpha t^{1-\alpha+\vrho} \Big(\left\|u^1-u^2\right\|_{\cnt^\vrho_{r,p}(T)}+\left\|w^1-w^2\right\|_{\cnt^\vrho_{r,p}(T)}\Big).
\end{align}
\end{lem}

\begin{proof}
To prove the estimates, we first consider the evolution operators $\Tst^{i,n}$, $i=1,2$, generated by the family of bounded linear operators, i.e., the families of Yosida approximations $A^i_{n,t}$ corresponding to $A^i_{t}$, $i\in\{1,2\}$.
We begin by observing that
\begin{align*}
(A^1_{n, t})^{\alpha}\left[\Tzt^{1,n}-\Tzt^{2,n}\right] (A^2_{n,0})^{-1} &= - (A^1_{n, t})^{\alpha} \Bigg[ \int_0^t \frac{\partial}{\partial s} \left[ \T^{1,n}_{s,t} \T^{2,n}_{0,s} \right] ds \Bigg](A^2_{n,0})^{-1}\\[0.1cm]
& = (A^1_{n, t})^{\alpha} \Bigg[\int_0^t \Tst^{1,n} A^1_{n,s}\left[(A^1_{n,s})^{-1}-(A^2_{n,s})^{-1}\right] A^2_{n,s} \Tzs^{2,n} ds \Bigg](A^2_{n,0})^{-1}. \\[0.1cm]
& =\int_0^t (A^1_{n, t})^{^\alpha} \Tst^{1,n} A^1_{n,s}\left[(A^1_{n,s})^{-1}-(A^2_{n,s})^{-1}\right] A^2_{n,s} \Tzs^{2,n} (A^2_{n,0})^{-1} d s .
\end{align*}

\noindent Letting $n\rightarrow\infty$ gives
\begin{align*}
\begin{aligned}
(A^1_t)^{\alpha}\left[\Tzt^1-\Tzt^2\right] (A^2_0)^{-1} &= \int_0^t (A^1_t)^\alpha \Tst^1 (A^1_s)^{1-\vartheta} (A^1_s)^{\vartheta}\left[(A^1_s)^{-1}- (A^2_s)^{-1}\right] A^2_s \Tzs^2 (A^2_0)^{-1} d s .
\end{aligned}
\end{align*}
Hence,
{\allowdisplaybreaks
\begin{align*}
&\left\|[(A^1_t)^\alpha\left[\Tzt^1-\Tzt^2\right]] U_0\right\|_{\V_p}= \left\|[(A^1_t)^\alpha\left[\Tzt^1-\Tzt^2\right] (A^2_0)^{-1} ] A^2_0 U_0\right\|_{\V_p} \\[0.2cm]
& \quad \leq \int_0^t \big\|(A^1_t)^\alpha \Tst^1 (A^1_s)^{1-\vartheta}\big\|_{\L(\V_p)} \big\|(A^1_s)^{\vartheta}\left[(A^1_s)^{-1}- (A^2_s)^{-1}\right]\|_{\L(\V_p)} \big\|A^2_s \Tzs^2 (A^2_0)^{-1}\big\|_{\L(\V_p)} \big\| A^2_0U_0 \big\|_{\V_p} d s \\[0.2cm]
& \quad \leq C \big\|A^2_0 U_0\big\|_{\V_p} \int_0^t(t-s)^{\vartheta-\alpha-1} 
 \Big( \big\|w^1_s - w^2_s \big\|_{\wrp} + \big\|v^1_s - v^2_s \big\|_{\wtp} \Big)  d s, \quad (\rm{cf.} ~ \eqref{eq:Atheta_T_Aphi_theta_bdd} ~~\text{\rm and } ~ Lemma~\ref{lem:Au_property}) \\[0.2cm]
& \quad \leq C \|U_0\|_{\wtp} \left\|w^1-w^2\right\|_{\cnt^\vrho_{r,p}(T)}  \int_0^t(t-s)^{\vartheta-\alpha-1} s^{\vrho} d s \quad  \text{(by Lemma~\ref{lem:Lip_ineq_normal} and the $\cnt^\vrho$-norm of $w^1-w^2$)} \\[0.2cm]
& \quad \leq C t^{1-\alpha+\vrho+\vartheta-1} \|U_0\|_{\wtp} \left\|w^1-w^2\right\|_{\cnt^\vrho_{r,p}(T)} \\[0.2cm]
& \quad \leq C t^{1-\alpha+\vrho} \big\|U_0\big\|_{\wtp} \left\|w^1-w^2\right\|_{\cnt^\vrho_{r,p}(T)}, \quad ~~\text{\rm for \aall}\,\, \om \in \Om\,,
\end{align*}
}
\noindent where the last inequality follows from the $\cnt^\vrho$-norm of $w^1-w^2$ and the fact that 

$$\|w^1_s - w^2_s\|+\|v^1_s-v^2_s\|\le C s^{\vrho}\|w^1 - w^2\|_{\cnt^\vrho}.$$
This proves the first assertion \eqref{eq:Aalp_T12_U0_lip}. The second estimate follows by a similar argument. Indeed,
\begin{align*}
(A^1_{n, t})^\alpha \int_0^t\left[\Tst^{1,n} - \Tst^{2,n}\right] F(s) ds &= \int_0^t \int_s^t (A^1_{n, t})^\alpha \Ttaut^{1,n} A^1_{n,\tau}\left[(A^1_{n,\tau})^{-1}-(A^2_{n,\tau})^{-1}\right] A^2_n(\tau) \Tstau^{2,n} F(s) d\tau \: ds \\[0.2cm]
= & \int_0^t (A^1_{n, t})^\alpha \Ttaut^{1,n} A^1_{n,\tau} \left[(A^1_{n,\tau})^{-1} - (A^2_{n,\tau})^{-1} \right] A^2_{n,\tau} \int_0^\tau \Tstau^{2,n} F(s) \: ds \: d\tau .
\end{align*}
\noindent By Theorem 3.10~\cite{Yagi09}, since 
\begin{align*}
\left\|A^2_{n,\tau} \int_0^\tau \Tstau^{2,n} F(s) d s\right\|_{\V_p} \leq C \tau^{1-1}\|F\|_{\Y^{1,\vsig}((0, T];\V_p)}
\end{align*}
and, for $n \rightarrow \infty$, converges to a constant function $c(\tau)$ on $(0, T]$, then letting $n \rightarrow \infty$, we obtain that
\begin{align*}
(A^1_t)^\alpha \int_0^t\left[\Tst^1-\Tst^2\right] F(s) d s = \int_0^t (A^1_t)^\alpha \Tstau^{1} (A^1_{\tau})^{1-\vartheta} (A^1_{\tau})^{\vartheta}\left[(A^1_{\tau})^{-1} - (A^2_{\tau})^{-1} \right] c(\tau) d \tau .
\end{align*}
\noindent Hence,
{\allowdisplaybreaks
\begin{align*}
&\left\| (A^1_t)^\alpha \int_0^t\left[\Tst^1 - \Tst^2\right] F(s) d s\right\|_{\V_p} \leq \int_0^t \left\|(A^1_t)^\alpha \Tstau^{1} (A^1_{\tau})^{1-\vartheta} \right\|_{\L(\V_p)} \left\| (A^1_{\tau})^{\vartheta}\left[(A^1_{\tau})^{-1} - (A^2_{\tau})^{-1} \right] \right\|_{\L(\V_p)} \left\| c(\tau) \right\|_{\V_p} \, d\tau \\[0.2cm]
&\leq C \int_0^t(t-\tau)^{\vartheta-\alpha-1}  \Big( \big\|w^1_{\tau}-w^2_{\tau}\big\|_{\wrp} + \big\|v^1_{\tau}-v^2_{\tau}\big\|_{\wtp} \Big)  d\tau ~\|F\|_{\Y^{1,\vsig}((0, T] ; X)} \quad (\rm{c.f.} ~\eqref{eq:Atheta_T_Aphi_theta_bdd} ~~\text{\rm and } ~ Lemma~\ref{lem:Au_property}) \\[0.2cm]
& \leq C \left\|w^1 - w^2 \right\|_{\cnt^\vrho_{r,p}(T)} \|F\|_{\Y^{1,\vsig}} \int_0^t(t-\tau)^{\vartheta-\alpha-1} \tau^{\vrho} d \tau \quad \text{(by Lemma~\ref{lem:Lip_ineq_normal} and the $\cnt^\vrho$-norm of $w^1-w^2$)} \\[0.2cm]
& \leq C t^{\vrho+1+\vartheta-\alpha-1} \left\|w^1 - w^2 \right\|_{\cnt^\vrho_{r,p}(T)} \|F\|_{\Y^{1,\vsig}} \\[0.2cm]
& \leq C t^{1-\alpha+\vrho} \left\|w^1 - w^2 \right\|_{\cnt^\vrho_{r,p}(T)} \|F\|_{\Y^{1,\vsig}}.
\end{align*}}
where, like above, we have used the estimate 

$$\|w^1_\tau - w^2_\tau\|+\|v^1_\tau - v^2_\tau\|\le C \tau^{\vrho}\|w^1 - w^2\|_{\cnt^\vrho}.$$
This establishes the second inequality \eqref{eq:Aalp_T12_F_lip}. 
\noindent Finally, 
\begin{align*}
\begin{aligned}
\left\|(A^2_t)^\alpha \int_0^t \Tst^2\left[F^1-F^2\right] d s\right\|_{\V_p} & \leq \int_0^t \left\|(A^2_t)^\alpha \Tst^2\left[F^1-F^2\right]\right\|_{\V_p} d s \le \int_0^t \left\|(A^2_t)^\alpha \Tst^2\right\|_{\L(\V_p)} \left\|F^1-F^2\right\|_{\V_p} d s \\[0.2cm]
&\hspace*{-3cm} \leq \int_0^t (t-s)^{-\alpha} \Big( \big\|u^1_s-u^2_s\big\|_{\wrp} + \big\|v^1_s-v^2_s\big\|_{\wtp} \Big) d s, \quad (\rm{c.f.} ~  \eqref{eq:Atheta_T_Aphi_theta_bdd} ~~\text{\rm and } ~ Lemma~\ref{lem:Fuv_property}) \\[0.2cm]
&\hspace*{-3cm}\leq C \Big(\left\|u^1-u^2\right\|_{\mathcal{C}^{\vrho}\left([0, T] ; \wrp \right)}+\left\|w^1-w^2\right\|_{\mathcal{C}^{\vrho}\left([0, T] ; \wrp \right)}\Big) \int_0^t (t-s)^{-\alpha} s^{\vrho}\, d s \\[0.2cm]
&\hspace*{-3cm}\leq C t^{1-\alpha+\vrho} \Big(\left\|u^1-u^2\right\|_{\mathcal{C}^{\vrho}\left([0, T] ; \wrp \right)}+\left\|w^1-w^2\right\|_{\mathcal{C}^{\vrho}\left([0, T] ; \wrp \right)}\Big), \quad 0 \leq t \leq T,
\end{aligned}
\end{align*}
where the $v$-difference is converted to a $w$-difference by Lemma~\ref{lem:Lip_ineq_normal} (applied to the inputs $w^1,w^2$). 
This establishes the final inequality \eqref{eq:Aalp_T2_F12_lip}.
\end{proof}
\noindent We are now ready to prove that the mapping $\Phi$ introduced in Lemma~\ref{lem:sem_lin_aprx_soln} is closed and contractive on $\Zb(T)\subset\mathcal{Y}(T)$.
\begin{lem}[Closedness and contraction]\label{lem:contraction}
Let $T>0$.
\begin{enumerate}
\item[\textup{(i)}] { {The set $\Zb(T)$ is closed in $\mathcal{Y}(T)$, hence complete with respect to $\|\cdot\|_{\mathcal{Y}(T)}$.}}
\item[\textup{(ii)}] For sufficiently small $T>0$, the mapping $\Phi:\Zb(T) \to \Zb(T)$ defined in Lemma~\ref{lem:sem_lin_aprx_soln} is a contraction on $\Zb(T)$ with respect to $\|\cdot\|_{\mathcal{Y}(T)}$.
\end{enumerate}
\end{lem}

\begin{proof}
\noindent {{First, we show that $\Zb(T)$ is closed in $\mathcal{Y}(T)$. Let $(w^n)_{n\ge1} \subset \Zb(T)$ be a Cauchy sequence in $\mathcal{Y}(T)$ and let $w \in \mathcal{Y}(T)$ denote its limit. Since $\mathcal{Y}(T)$ is a Banach space, $w^n \to w$ in $L_p(\Om;\ztp(T))$ as well as in both components of $\X(T)$. It remains to verify the defining constraints of $\Zb(T)$.
For a.a.\ $\om \in \Om$, $\|w^n(\om)-w(\om)\|_{\zsp(T)} \to 0$, hence $\|w^n(\om,0)-w(\om,0)\|_{\wsp} \to 0$. Since $w^n(0)=u_0$ for every $n$, we obtain $w(0)=u_0$.
Since $\vsig<\vrho$, the map

$$
\mathcal{C}^{\vrho}([0,T];\wrp) \ni f \mapsto [f]_{\vsig} \in \R_+
$$
is continuous with respect to $\|\cdot\|_{\cnt^{\vrho}([0,T];\wrp)}$. Indeed, for $f,g \in \cnt^{\vrho}([0,T];\wrp)$ and $0 \le s < t \le T$,

$$
\frac{\|f(t)-f(s)\|_{\wrp}}{(t-s)^{\vsig}} - \frac{\|g(t)-g(s)\|_{\wrp}}{(t-s)^{\vsig}} \le \frac{\|f(t)-f(s)-g(t)+g(s)\|_{\wrp}}{(t-s)^{\vsig}} \le 2 (t-s)^{\vrho-\vsig} \|f-g\|_{\cnt^{\vrho}([0,T];\wrp)}.
$$
Taking the supremum over $0 \le s < t \le T$ yields $|[f]_{\vsig}-[g]_{\vsig}| \le 2 T^{\vrho-\vsig} \|f-g\|_{\cnt^{\vrho}([0,T];\wrp)}$. Because $w^n \to w$ in $L_p(\Om;\cnt^{\vrho}_{r,p}(T))$, there exists a subsequence, still denoted by $(w^n)$, such that $w^n(\om) \to w(\om)$ in $\cnt^{\vrho}([0,T];\wrp)$ for a.a.\ $\om \in \Om$. Therefore $[w^n(\om)]_{\vsig} \to [w(\om)]_{\vsig}$ for a.a.\ $\om \in \Om$, and Fatou's lemma gives

$$
\E\big( {[w]_{\vsig}}^p \big) \le \liminf{n\to\infty} \E\big( {[w^n]_{\vsig}}^p \big) \le 1.
$$
Thus, $w \in \Zb(T)$, and $\Zb(T)$ is closed in $\mathcal{Y}(T)$.}}

\smallskip\noindent Next, we show that the mapping $w\mapsto\Phi(w)=U(w)$ maps $\Zb(T)$ into itself. To this end, we note that
{\allowdisplaybreaks
\begin{align*}
U_t = & U_0+\left[e^{-t A_{w}(0)}-1\right] U_0+\left[\Tzt(w)-e^{-t A_w(0)}\right] A_w(0)^{-\beta} A_w(0)^\beta U_0 \\[0.1cm]
& +\int_0^t \Tst(w)[F(s)-F(t)] d s+\int_0^t\left[\Tst(w)-e^{-(t-s) A_w(s)}\right] d s F(t) \\[0.1cm]
& +\int_0^t\left[e^{-(t-s) A_w(s)}-e^{-(t-s) A_w(0)}\right] d s F(t)-\left[e^{-t A_{w}(0)}-1\right] A_{w}(0)^{-1} F(t).
\end{align*}}
Then, using the continuous embedding $\D(A^{\beta})\hookrightarrow\W^s_p$ (\cf ~ Remark 16.6 of~\cite{Yagi09}), we obtain
\begin{align*}
\left\|\left[e^{-t A_{w}(0)}-1\right] U_0\right\|_{\wsp} &\leq C\left\|\left[e^{-t A_{w}(0)}-1\right] A^{\beta}_{w}(0) U_0\right\|_{\V_p} \leq C\left\|\left[e^{-t A_{w}(0)}-1\right] A^{\beta-1}_{w}(0) \right\|_{\L(\V_p)}  \left\| A_{w}(0) U_0\right\|_{\V_p} \\
    &\leq C\left\|\left[e^{-t A_{w}(0)}-1\right]A^{\beta-1}_w\right\|_{\L(\V_p)} \const_{\Zb}  \quad ~~\text{\rm for \aall } \om \in \Om\,.
\end{align*}
The $C_0$-continuity property of the semigroup $e^{-tA_w(0)}$ gives 
\begin{align}\label{eq:C0_cont_eAw}
\left\|\left[e^{-t A_{w}(0)}-1\right] A^{\beta}_{w}(0) U_0 \right\|_{\V_p} &= \left\|\int_0^t A^{1}_w(0) e^{-s A_w} A^{\beta-1}_{w}(0) A_{w}(0) U_0 ds\right\|_{\V_p}\\[0.2cm]
&\le C \|A_w(0) U_0\|_{\V_p} \int_0^t \tau^{-\beta} d\tau =  \frac{C}{1-\beta} \|U_0\|_{\wtp} t^{1-\beta} \\[0.2cm]
\implies \left\|\left[e^{-t A_{w}(0)}-1\right] U_0\right\|_{\wsp} &\to 0 \quad ~~\text{\rm as} \quad t \to 0 \quad ~~\text{\rm almost surely}.
\end{align}
Utilizing \eqref{eq:Atheta_T_minus_E_Aphi} and the embedding $\D(A^{\beta})\hookrightarrow \wsp$ (\cf.~\eqref{eq:U_Soln_embds}) obtained in Theorem~\ref{thm:U_Soln}, we can apply ${\big\|f\big\|_{\wsp} \le C\big\|A^{\beta}_w f\big\|_{\V_p}}$ to verify that
{\allowdisplaybreaks
\begin{align*}
&\E \left[\left\|\left[\Tzt(w)-e^{-t A_w(0)}\right] A_w(0)^{-1} A^{}_w(0) U_0\right\|_{\wsp} \right]\\[0.2cm] &\leq C_2\left\|A^{\beta}_w(t)\left[\Tzt(w)-e^{-t A_w(0)}\right] A_w(0)^{-1}\right\|_{\L(\V_p)} \E \left[\left\|A_w(0) U_0\right\|_{\V_p}\right] \leq C t^{1+\vsig-\beta}  \E \left[\left\|A_w(0) U_0\right\|_{\V_p}\right]. 
\end{align*}}
By~\eqref{eq:AthetaT_bdd} and the embedding~\eqref{eq:U_Soln_embds}, we obtain
\begin{align*}
\E \left[ \left\|\int_0^t \Tst(w)[F(s)-F(t)] d s\right\|_{\wsp} \right] &\leq C \int_0^t(t-s)^{\vsig-\beta} s^{-\vsig} ds ~~ \E\left[\sup{0 \le s<t} \frac{s^{\vsig}\|F(t)-F(s)\|_{\V_p}}{(t-s)^\vsig}\right] \\[0.2cm]
&\leq C \: t^{1-\beta} \: \text{Beta}(1-\vsig,1+\vsig-\beta) \: \E \left[\|F\|_{\Y^{1,\vsig}}\right],
\end{align*}
where $\text{Beta}(x,y)$ is the beta function. Similarly, by~\eqref{eq:Atheta_T_minus_E_Aphi} and the embeddings~\eqref{eq:U_Soln_embds} we get
\begin{align*}
\E \left[\left\|\int_0^t\left[\Tst(w)-e^{-(t-s) A_w(s)}\right] ds F(t)\right\|_{\wsp}\right] \leq C \int_0^t (t-s)^{\vsig - \beta} ds \: \E\left[\big\|F\big\|_{\Y^{1,\vsig}}\right] &\leq C \: t^{1+\vsig-\beta} \: \E \left[ \big\|F\big\|_{\Y^{1,\vsig}} \right].
\end{align*}
By Lemma 3.2 in~\cite{Yagi09} and the embedding~\eqref{eq:U_Soln_embds}, 
\begin{align*}
\begin{aligned}
\E \left[\left\|\int_0^t\left[e^{-(t-s)A_w(s)} - e^{-(t-s)A_w(0)} \right] ds F(t) \right\|_{\wsp}\right]  &\leq C \int_0^t (t-s)^{-\beta} s^\vsig ~ ds ~ \E \left[\big\|F\big\|_{\Y^{1,\vsig}} \right] \\[0.2cm]
&\leq C t^{1+\vsig-\beta} \text{Beta}(1+\vsig,1-\beta) \E \left[\big\|F\big\|_{\Y^{1,\vsig}} \right]\,.
\end{aligned}
\end{align*}
Finally, by~\eqref{eq:bded_sob_embd}, \eqref{eq:dom_embds}, and \eqref{eq:weighHol_cont}, we can write
\begin{align*}
\E \left[\left\|\left[e^{-t A_{w}(0)}-1\right] A_{w}(0)^{-1} F(t)\right\|^p_{\wsp}\right] &\leq C_2\left\|\left[e^{-t A_{w}(0)}-1\right]A^{\beta-1}\right\|_{\L(\V_p)} \E \left[\|F\|^p_{\Y^{1,\vsig}}\right] \\[0.2cm]
\overset{\text{similar to } \eqref{eq:C0_cont_eAw}} {\implies}\bigg\|\left[e^{-t A_{w}(0)}-1\right] A_{w}(0)^{-1} F(t)\bigg\|_{\wsp} &\to 0 \quad ~~\text{\rm as} \quad t \to 0.
\end{align*}
 We have therefore verified that
\begin{align*}
\E \left[\big\|U(t)-U(0)\big\|_{\wsp}\right]  &\leq ~ \E \left[ \left\| h(t) \right\|\right] + C~t^{(1-\beta)} \left(\E \left[ \big\|F\big\|_{\Y^{1, \vsig}}\right] + \E \left[\big\|U_0\big\|_{\wsp}\right] \right), \quad 0<t \leq T,
\end{align*}
where $h(t):[0,T]\to \R_+$ is a continuous function such that $h(t) \to 0$ as $t \to 0$. Hence, choosing the ball radius $\const_{\Zb}$ large enough depending on $\|U_0\|_{\wtp}$ and then restricting to a sufficiently small time horizon $T>0$ (depending on $\const_{\Zb}$, $\|U_0\|$, and the a priori bound on $F$), we obtain
\begin{align}
\E \left[ \big\|U\big\|^p_{\zsp} \right] \leq \const^p_{\Zb}.
\end{align}
Furthermore, by Theorem~\ref{thm:U_Soln}, we have 
\begin{align}\label{eq:exp_uwtp_norm}
\E \left[\big\|U\big\|^p_{\wtp}\right] \le C ~ \E \left[\big\|A_w(t) U(t)\big\|^p_{\V_p}\right] \leq C\left(\E \left[\big\|F\big\|^p_{\Y^{1, \vsig}} \right] + \E \left[\big\|U_0\big\|^p_{\wtp}\right] \right) .
\end{align}
To verify the H\"older condition, from (4.14) in~\cite{Yagi09} we notice that
\begin{align*}
U(t)= & \Tst(w) \T_{0,s}(w) U_0+\Tst(w) \int_0^s \T_{\tau,s}(w) F(\tau) d \tau + \int_s^t \T_{\tau,t}(w) F(\tau) d \tau .
\end{align*}
Therefore, for $0 \leq s<t \leq T$ ,
\begin{align*}
U(t)-U(s)&= {\left[\Tst(w)-1\right] U(s)+\int_s^t \T_{\tau,t}(w) F(\tau) d \tau } \\[0.2cm]
&= \left\{\left[\Tst(w)-e^{-(t-s) A_w(s)}\right] \right. + \left. \left[e^{-(t-s) A_w(s)}-1\right]\right\} A_w(s)^{-1} A_w(s) U(s)+\int_s^t \T_{\tau,t}(w) F(\tau) d \tau.
\end{align*}

\noindent Then, by~\eqref{eq:Atheta_T_minus_E_Aphi} and~\eqref{eq:C0_cont_eAw}, we obtain
\begin{align*}
& \left \| A^{\alpha}_w(t)\left[\Tst(w)-e^{-(t-s) A_w(s)}\right] A_w(s)^{-1} \right\|_{\L(\V_p)} \leq C(t-s)^{1+\vsig-\alpha}, \\[0.1cm]
& \left\|A^{\alpha}_w(s)\left[e^{-(t-s) A_w(s)}-1\right] A_w(s)^{-1}\right\|_{\L(\V_p)} \leq C (t-s)^{1-\alpha},
\end{align*}
respectively, while by~\eqref{eq:AthetaT_bdd} we get,
\begin{align*}
\left\|\int_s^t A_w(t)^\alpha U_w(t, \tau) F(\tau) d \tau\right\|_{\V_p} &\leq C ~ \big\|F\big\|_{\Y^{1,\vsig}} \int_s^t(t-\tau)^{-\alpha} d \tau \quad \leq C ~ \big\|F\big\|_{\Y^{1,\vsig}} (t-s)^{1-\alpha}.
\end{align*}
In view of Theorem~\ref{thm:U_Soln} and~\eqref{eq:exp_uwtp_norm}, we therefore obtain that
\begin{align*}
\big\|U(t)-U(s)\big\|_{\wrp} &\leq C\left( \big\|F\big\|_{\Y^{1,\vsig}} + \big\|U_0\big\|_{\wsp} \right)(t-s)^{1-\alpha} \\[0.2cm]
\implies \frac{\big\|U(t)-U(s)\big\|^p_{\wrp}}{(t-s)^{\vsig}} &\leq C^p\left(\big\|F\big\|_{\Y^{1,\vsig}} + \big\|U_0\big\|_{\wsp} \right)(t-s)^{(1-\alpha - \vsig)},\quad \text{ for } t>s.
\end{align*}

\noindent Hence, for $\vsig < \vrho = 1-\alpha$,
\begin{align*}
\E \left[\left(\sup{0\le s < t \le T} \frac{\big\|U(t)-U(s)\big\|_{\wrp}}{|t-s|^{\vsig}} \right)^p\right] \leq C ~ \left(\E \left[\big\|F\big\|_{\Y^{1,\vsig}}\right] + \E \left[\big\|U_0\big\|_{\wsp}\right] \right) ~ T^{(\vrho-\vsig)}, \quad 0 \leq s \leq t \leq T.
\end{align*}
Thus, for $T>0$ sufficiently small we get 
\begin{align*}
    \E \left[\left(\sup{0\le s < t \le T} \frac{\big\|U(t)-U(s)\big\|_{\wrp}}{|t-s|^{\vsig}} \right)^p\right] \leq 1\,.
\end{align*}
Hence, $U \in \Zb(T) \cap \kploc$.\\

\noindent Next, we look at the contraction property. Specifically, let $w^1, w^2 \in \Zb$ and  
$$
U^1:= u(w^1),\,\,\,\, U^2 := u(w^2),\,\,\,\, F^1 := F[w^1] = F\big(U^1,v(w^1)\big),\,\,\,\, F^2 := F[w^2] = F\big(U^2,v(w^2)\big).
$$ 
Then,
\begin{align*}
    U^1_t - U^2_t = [\T^1_{0,t} - \T^2_{0,t}] u_0 + \int_0^t [\Tst^1 - \Tst^2]F_1(s) ds + \int_0^t \Tst^2[F^1(s) - F^2(s)] ds.
\end{align*}
We now estimate each of these terms. From Lemma~\ref{lem:Op_Sg_Lip_ineq} with $\vartheta=1$ and $\alpha = \beta$ (as per Theorem~\ref{thm:U_Soln}) we obtain that for $0 \leq t \leq T$
\begin{align*}
\left\|(A^1_t)^{\beta}\left[U^1_t - U^2_t\right]\right\|_{\V_p} & \leq C T^{1+\vrho-\beta}\Big(\big\|F\big\|_{\Y^{1,\vsig}}+\big\|U_0\big\|_{\wrp} + 1\Big)\Big(\big\|U^1-U^2\big\|_{\cnt^\vrho_{r,p}(T)}+\big\|w^1-w^2\big\|_{\cnt^\vrho_{r,p}(T)}\Big).
\end{align*}
Hence
\begin{align}\label{eq:U12_Bnorm}
\big\|U^1_t-U^2_t\big\|_{\zsp} \leq C\left(\big\|F\big\|_{\Y^{1,\vsig}} + \big\|U_0\big\|_{\wrp} + 1\right) T^{1+\vrho-\beta}\Big(\big\|U^1 - U^2\big\|_{\cnt^\vrho_{r,p}(T)}+\big\|w^1 - w^2\big\|_{\cnt^\vrho_{r,p}(T)}\Big) .
\end{align}
The $\cntrrp$ norm of $U_1-U_2$ is also estimated in a quite similar way by applying Lemma~\ref{lem:Op_Sg_Lip_ineq}. Indeed,
\begin{align*}
& \left\|(A^1_t)^{\alpha} \int_s^t\left[\Tst^1-\Tst^2\right] F(s) d s\right\|_{\V_p} \leq C (t-s)^{1+\vrho-\alpha}\big\|F\big\|_{\Y^{1,\vsig}}\big\|w^1 - w^2\big\|_{\cnt^\vrho_{r,p}(T)}, \quad 0 \leq s < t \leq T.
\end{align*}
Finally, for  $0 \leq s < t \leq T$,
\begin{align*}
\begin{aligned}
\left\|(A^2_t)^{\alpha} \int_s^t \Tst^2\left[F^1-F^2\right] d s\right\|_{\V_p} \le C (t-s)^{1+\vrho-\alpha} \Big(\left\|U^1-U^2\right\|_{\cnt^\vrho_{r,p}(T)}+\left\|w^1-w^2\right\|_{\cnt^\vrho_{r,p}(T)}\Big).
\end{aligned}
\end{align*}

\noindent Altogether, we get that
\begin{align*}
\begin{aligned}
\big\|(U^1_t - U^1_s) - (U^2_t-U^2_s)\big\|_{\wrp} \le & \,\, C (t-s)^{1+\vrho-\alpha}\big\|F\big\|_{\Y^{1,\vsig}}\big\|w^1 - w^2\big\|_{\cnt^\vrho_{r,p}(T)} \quad 0 \leq s < t \leq T \\[0.2cm]
&+ C (t-s)^{1+\vrho-\alpha} \Big(\big\|U^1-U^2\big\|_{\cnt^\vrho_{r,p}(T)}+\big\|w^1-w^2\big\|_{\cnt^\vrho_{r,p}(T)}\Big) .
\end{aligned}
\end{align*}

\noindent Thus, for $0<S \leq T$, we have that
\begin{align}
\label{eq:U12_Crhonorm}
\begin{aligned}
\quad \sup{0 \leq s < t \leq S} (t-s)^{-\vrho} \Big\|(U^1(t) - U^1(s)) - (U^2(t) - U^2(s)) \Big\|_{\wrp}  &\leq C T^{\vrho}\left(\big\|F\big\|_{\Y^{1,\vsig}} + \big\|U_0\big\|_{\wrp} + 1\right) \\
&\hspace*{-1cm} \times \Big(\big\|U^1 - U^2\big\|_{\cnt^\vrho_{r,p}(S)}+\big\|w^1 - w^2\big\|_{\cnt^\vrho_{r,p}(S)}\Big).
\end{aligned}
\end{align}
Therefore, from ~\eqref{eq:U12_Bnorm} and~\eqref{eq:U12_Crhonorm} it yields that
\begin{align*}
\left\|U^1-U^2\right\|_{\X{(T)}} \leq C T^{\vrho}\left(\big\|F\big\|_{\Y^{1,\vsig}}+\big\|U_0\big\|_{\wrp} + 1\right)\Big(\big\|U^1 - U^2\big\|_{\X(T)}+\big\|w^1 - w^2\big\|_{\X(T)}\Big), \quad w^1, w^2 \in \Zb(T).
\end{align*}
Absorbing the $\big\|U^1-U^2\big\|_{\X(T)}$ term on the left-hand side for $T>0$ sufficiently small, we obtain
\begin{align*}
\left\|U^1-U^2\right\|_{\X{(T)}} \leq C T^{\vrho}\left(\big\|F\big\|_{\Y^{1,\vsig}}+\big\|U_0\big\|_{\wrp} + 1\right)\big\|w^1 - w^2\big\|_{\X(T)}, \quad w^1, w^2 \in \Zb(T).
\end{align*}
By~\eqref{eq:U12_Bnorm} and Theorem~\ref{thm:U_Soln} (which places $U^i(\om)\in\bU(T)\subset\ztp$ for a.a.\ $\om$), the same estimates also control $\|U^1-U^2\|_{L_p(\Om;\ztp(T))}$, hence
\begin{align*}
\left\|U^1-U^2\right\|_{\mathcal{Y}(T)} \leq C T^{\vrho}\left(\big\|F\big\|_{\Y^{1,\vsig}}+\big\|U_0\big\|_{\wrp} + 1\right)\big\|w^1 - w^2\big\|_{\mathcal{Y}(T)}, \quad w^1, w^2 \in \Zb(T).
\end{align*}
Hence, $\Phi$ is a contraction mapping on $(\Zb(T),\|\cdot\|_{\mathcal{Y}(T)})$, provided that $T > 0$ is sufficiently small.
\end{proof}
We are now in a position to state and prove the main existence and uniqueness result for the macroscopic system~\eqref{eq:mult_glioma_sys}.
\begin{thm}[Existence and uniqueness of solution to the macroscopic system]
\label{thm:macro_existence}
Let Assumptions~\ref{asm:init_cond} and~\ref{asm:model_params} hold, let $p>n$, and let $r,s,s_1,s_2,\vartheta$ and the H\"older exponent $\vsig$ be as in Theorem~\ref{thm:U_Soln} (in particular~\eqref{eq:U_Soln_vsig_Yagi}). Then, there exists $T^* \in (0,T]$ such that the macroscopic system~\eqref{eq:mult_glioma_sys}, together with the exponential Ornstein-Uhlenbeck process~\eqref{eq:noise}, admits a unique pair $(u,v)$ satisfying

$$
u \in \X(T^*),\qquad v(\om,t,\cdot)=v_0(\cdot)\exp\!\left(-\int_0^t\big(\kappa_3\psi(u_s(\cdot))+R_Q(\xi_s)\big)\,ds\right) \quad \text{a.a.\ }\om\in\Om,
$$
{{Moreover, $(u,v)$ satisfies~\eqref{eq:pre_cancer}--\eqref{eq:pre_normal} in the sense of Theorem~\ref{thm:U_Soln}, pathwise for almost every $\om\in\Om$. Furthermore, $u(t,x)\ge 0$ for $\avery x\in\dom$ and $\aall t\in[0,T^*]$, for a.a.\ $\om\in\Om$.}}
\end{thm}
\begin{proof}
Let $\Phi$ be the mapping defined in Lemma~\ref{lem:sem_lin_aprx_soln}. By Lemma~\ref{lem:contraction}-(i), $(\Zb(T^*),\|\cdot\|_{\mathcal{Y}(T^*)})$ is a complete metric space. By Lemma~\ref{lem:contraction}-(ii), for sufficiently small $T^*>0$, $\Phi$ maps $\Zb(T^*)$ into itself and is a contraction with respect to $\|\cdot\|_{\mathcal{Y}(T^*)}$. The Banach fixed point theorem therefore yields a unique fixed point
$u\in\Zb(T^{\ast})\subset\mathcal{Y}(T^{\ast})$ such that $\Phi(u)=u$. The corresponding tissue density $v$ is then given by Lemma~\ref{lem:sol_normal}. Therefore, the pair $(u,v)$ solves the macroscopic system~\eqref{eq:mult_glioma_sys}. 

{Finally, nonnegativity of $u$ follows from the invariant set $K:=\{w\in\wsp: w\ge 0\ \text{a.e.\ in }\dom\}$. By Assumption~\ref{asm:init_cond}, $u_0\in K$, and $K$ is invariant for~\eqref{eq:ACP_Can}: at $u=0$ the reaction term in~\eqref{eq:pre_cancer} vanishes, while $g(0,v)=\kappa_2\zeta(0,v)=0$ (cf.~\eqref{eq:pre_cancer}) nullifies the taxis flux in~\eqref{eq:ACP_Ops}, and the no-flux boundary conditions in~\eqref{eq:pre_cancer} yield zero outward normal transport on $\{u=0\}$. Therefore $\Phi(w)\in K$ for every $w\in\Zb(T^{\ast})$. By the invariant set method~\cite{Yagi09}, the fixed point $u=\Phi(u)$ satisfies $u(t,x)\ge 0$ for $\avery x\in\dom$ and $\aall t\in[0,T^{\ast}]$, for a.a.\ $\om\in\Om$.}
\end{proof}

\section{Discussion}\label{sec:discussion}
In this work, we analyzed the random PDE--ODE model proposed in~\cite{conte2026randomness} for glioma dynamics, which couples a reaction-diffusion-taxis PDE for tumor cells with an ODE for normal tissue, both involving stochasticity in the depletion terms that characterize radiotherapeutic effects. The stochastic process therein is of OU type, thus satisfies a corresponding stochastic differential equation (SDE). Here, we only proved local existence of solutions; extensions to global existence results for this highly nonlinear problem turn out to be highly challenging. Approaches relying - as in the deterministic framework - on more general fixed point theorems would require compactness arguments which cannot simply be translated from the deterministic to the stochastic case. The authors of~\cite{Zhigun2016} considered a prototype degenerate PDE-SDE system and used simultaneous estimates for Malliavin and weak Sobolev derivatives within a compactness argument to prove global existence of weak solutions, but it is not clear how this method could work for the more complex system studied in this note. Directly considering SPDEs with (space-) time-dependent white noise would have been even more involved, both from analytical and numerical perspectives, the latter raising, a.o., the issue of positivity preserving.

The model analyzed here is one of the very few existing continuous frameworks describing tumor growth and invasion and accounting for stochasticity. In~\cite{hiremath2015stochastic}, a random PDE-ODE model for tumor and proton dynamics was considered and local existence of the solution was proved. In~\cite{hiremath2016stochastic}, a method similar to the one used here was applied to a random PDE-ODE system with cross diffusion, being able to achieve a global result, however under rather generous assumptions on the involved coefficient functions. In~\cite{hiremath2018coupled} a complex PDE-SDE system describing acid-mediated tumor invasion with pH-taxis was investigated, and global existence and uniqueness were proved - again under adequate assumptions. 

The present model is more complex, as it involves haptotaxis instead of chemotaxis and features more intricate motility terms originating in the model deduction described in~\cite{conte2026randomness}. The SDE model for glioma growth introduced in~\cite{kroos2019} also investigated various therapy combinations via TCP/UTCP effects, but the focus there was on the interplay between stem cells and their differentiated counterparts, without accounting for spatial spread. These spatial aspects are, however, crucial for glioma diagnostics and therapy planning. Therefore, our model offers a means to connect therapy assessment with tumor margin predictions, and its analysis is a first step towards better understanding the mathematical structure of such settings.

\subsection*{Acknowledgment}
The research of MC has been carried out under the auspices of the National Group of Mathematical-Physics (GNFM-INdAM). MC acknowledges support from the European Union - Next Generation EU, Mission 4, Component 1 (CUP: E13C24002380006) and the National Group of Mathematical Physics (GNFM-INdAM) (CUP E5324001950001).  MC has been also partially supported by the State Research Agency of the Spanish Ministry of Science and FEDER-EU, project PID2022-137228OB-I00 (MICIU/AEI /10.13039/501100011033) and by Modeling Nature Research Unit, Grant QUAL21-011 funded by Consejería de Universidad, Investigaci\'on e Innovaci\'on (Junta de Andalucía). SAH acknowledges support from the Department of Mathematics at RPTU Kaiserslautern. This work was conducted in association with the CZS project \emph{AI-CARE}, which focuses on the development of PDE-based AI models and scRNA-seq data integration for discovering novel targets for GBM treatment.

\subsection*{Competing interests}
The authors declare that the research was conducted in the absence of any commercial or financial relationships that could be construed as a potential conflict of interest.

\counterwithin{figure}{section}
\counterwithin{table}{section}
\counterwithin{equation}{section}
\renewcommand\thefigure{\thesection\arabic{figure}}
\renewcommand\thetable{\thesection\arabic{table}}
\renewcommand\theequation{\thesection.\arabic{equation}}

\appendix

\section{Appendix} \label{sec:appendix}
In this section, we provide the detailed derivations of some of the results used in the main text.

\subsection{Exponential OU Process}
In this section, we provide the closed-form solution to the SDE characterizing the Ornstein-Uhlenbeck process $\xi_t$ employed in our model. Specifically, consider the process $\xi_t$ defined by the following SDE: 
\begin{equation}\label{eq:Exp_OU_proc}
    \left\{\begin{aligned} 
    d \xi_t & = \left(\mu_R-\xi_t\right) d t+\sigma_R \xi_t d W_t\,, \\[0.2cm] 
    \xi(0) & =\xi_0 \text { (given), }
    \end{aligned} \quad \left(\,\mu_R, \sigma_R>0\right).\right. 
\end{equation}
The SDE in \eqref{eq:Exp_OU_proc} admits the following explicit solution:
\begin{align}\label{eq:Exp_OU_proc_solved}
    \xi_t &= e^{-t \sigma_{\xi} + {\sigma}_R W_t}\left(\xi_0+\mu_R \int_0^t e^{s \sigma_{\xi}} e^{-\sigma_R W_s} d s\right), \quad \sigma_{\xi} := \left(1+\frac{\sigma_R^2}{2}\right).
\end{align}
Based on the construction of the solution, \eqref{eq:Exp_OU_proc_solved} is the unique strong solution of \eqref{eq:Exp_OU_proc}, hence $(\xi_t)_{t\ge0}$ is $(\Fil_t)$-adapted. Moreover, if $\xi_0\ge 0$ and $\mu_R>0$, then $\xi_t\ge 0$ almost surely for all $t\ge 0$, since the prefactor $e^{-t\sigma_{\xi}+\sigma_R W_t}>0$ and the integral term in~\eqref{eq:Exp_OU_proc_solved} is nonnegative. These properties of $\xi_t$, together with the explicit representation \eqref{eq:Exp_OU_proc_solved}, will be used repeatedly in the existence and regularity arguments for the macroscopic system developed in the main text.

\subsection{Moment estimates}
Based on the explicit solution~\eqref{eq:Exp_OU_proc_solved}, we derive moment estimates for the process $\xi_t$. To this end, it is convenient to rewrite the solution in a form that is better suited to subsequent estimates. In particular, exploiting the independence of increments of the Wiener process $W_t$ and the fact that $W_t-W_s \sim \mathcal{N}(0,t-s)$, we can write
\begin{align}
    \xi_t &= e^{-t \sigma_{\xi} + {\sigma}_R W_t}\left(\xi_0+\mu_R \int_0^t e^{s \sigma_{\xi}} e^{-\sigma_R W_s} d s\right) \\[0.1cm]
    &= \xi_0 e^{-t \sigma_{\xi}} e^{{\sigma}_R W_t} + \mu_R \int_0^t e^{-(t-s) \sigma_{\xi}} e^{\sigma_R (W_t-W_s)} d s \\[0.1cm]
    &\sim \xi_0 e^{-t \sigma_{\xi}} e^{{\sigma}_R \sqrt{t} w} + \mu_R \int_0^t e^{-(t-s) \sigma_{\xi}} e^{\sigma_R \sqrt{(t-s)} w} d s, \quad w \sim \mathcal{N}(0,1).
\end{align}

\begin{lem}[]\label{lem:mgf}
Let $w \sim \mathcal{N}(0,\sigma^2)$; then its moment-generating function $M_w(t)$ is given by $e^{\tfrac{1}{2}\sigma^2 t^2}$. 
\end{lem}
\begin{proof}
This is the standard Gaussian moment-generating function (MGF) which can be obtained by completing the square in the Gaussian integral (see, e.g., \cite{DaPrato1992}).
\end{proof}

\noindent Using the moment-generating function derived above, we deduce the following moment estimates.
\begin{lem} Let $\xi = (\xi_t)_{t\ge0}$ be the process defined in \eqref{eq:Exp_OU_proc_solved}. The following relations for its moments hold:
\begin{align*}
    \E\left[(d\xi_t)^2\right] &= \sigma_R^2  \E\left[\xi_t^2\right] dt\,, \\[0.2cm]
    \E\left[\xi_t\right] &= \mu_R + e^{-t}(\xi_0 - \mu_R)\,, \\[0.2cm]
    \E\left[\xi^2_t\right] &= \xi_0^2 e^{-t(2 - \sigma^2_R)} + {2 \mu_R \xi_0  \over (1 - {{\sigma}^2_R }) } \left(e^{-t} - e^{-t(2 - {{\sigma}^2_R})}\right) \\[0cm]
     & \quad+ {\mu^2_R \over (1 - \sigma^2_R)} \left(1 - e^{-t} - e^{-t(1 - \sigma^2_R)} + e^{-t(2 - \sigma^2_R)} \right).
\end{align*}
\end{lem}
\begin{proof}
The proof uses It\^o calculus, independence of Brownian increments, and Lemma \ref{lem:mgf}. We report here only the nontrivial identities.

From equation $d\xi_t=(\mu_R-\xi_t)dt+\sigma_R\xi_t dW_t$ and the identity $(dW_t)^2=dt$, we get
\begin{align*}
    \E[(d\xi_t)^2]=\sigma_R^2\E[\xi_t^2]dt\,,
\end{align*}
while, taking expectation in \eqref{eq:Exp_OU_proc_solved} and using
the fact that $\E\left[e^{\sigma_R(W_t-W_s)}\right]=e^{\frac12\sigma_R^2(t-s)}$, yields
\begin{align*}
    \E[\xi_t]=\mu_R+e^{-t}(\xi_0-\mu_R).
\end{align*}
For the second moment, it is convenient to decompose $\xi_t$ as $\xi_t=P_t+Q_t$ where
\begin{align*}
    P_t:=\xi_0e^{-t\sigma_\xi}e^{\sigma_RW_t},
    \qquad
    Q_t:=\mu_R\int_0^t e^{-(t-s)\sigma_\xi}e^{\sigma_R(W_t-W_s)}ds.
\end{align*}
Then, after straightforward evaluation of one-dimensional exponential integrals, we can obtain
\begin{align*}
    \E\left[P_t^2\right]&=\xi_0^2e^{-t(2-\sigma_R^2)},\\[0.15cm]
    \E\left[Q_t^2\right]&=\frac{\mu_R^2}{1-\sigma_R^2}\Big(1-e^{-t}-e^{-t(1-\sigma_R^2)}+e^{-t(2-\sigma_R^2)}\Big),\\[0.15cm]
    \E\left[P_tQ_t\right]&=\frac{\mu_R\xi_0}{1-\sigma_R^2}
    \Big(e^{-t}-e^{-t(2-\sigma_R^2)}\Big)\,.
\end{align*}
Collecting all three contributions gives
\begin{align}
    \E\left[\xi^2_t\right] &=  \E\left[P^2_t\right] + 2\E\left[P_t Q_t\right] + \E[Q^2_t] \notag \\[0.15cm]
    &= \xi_0^2 e^{-t(2 - \sigma^2_R)} +  {2 \mu_R \xi_0  \over (1 - {{\sigma}^2_R }) } \left(e^{-t} - e^{-t(2 - {{\sigma}^2_R})}\right) +  {\mu^2_R \over (1 - \sigma^2_R)} \left(1 - e^{-t} - e^{-t(1 - \sigma^2_R)} + e^{-t(2 - \sigma^2_R)} \right).
\end{align}
\end{proof}

\noindent These results immediately lead to the following estimates.
\begin{cor}
 Let $\sigma^2_R < 1$. Then, the first and second moments of $(\xi_t)_{t\ge0}$ satisfy the following inequalites:
\begin{align}
    \label{eq:exp_ou_bd}
    \begin{aligned}
    \E\left[\xi_t\right] &\le \mu_R + |\xi_0| \\[0.15cm]
    \E\left[\xi^2_t\right] &\le  \xi_0^2 +  {2 \mu_R \xi_0  \over (1 - {{\sigma}^2_R}) }  +  {\mu^2_R \over (1 - \sigma^2_R)}.
    \end{aligned}
\end{align}
\end{cor}

\noindent Now we derive higher moment estimates for the process defined in~\eqref{eq:Exp_OU_proc_solved}.
\begin{lem}[]
\label{lem:n_mgf}
Let $n > 1$ and $\sigma_R^2 < \frac{1}{n-1}$. Then, the $n$-th moment of the process
\begin{align*}
   Q_t = \mu_R \int_0^t e^{-(t-s) \sigma_{\xi}} e^{\sigma_R (W_t-W_s)} \, d s
\end{align*}
satisfies
\begin{align}
    \E[Q_t^n] = \mu^n_R \, e^{-q_n t} \prod_{k=1}^n \frac{e^{p_k t} - 1}{p_k},
\end{align}
where $\sigma_{\xi} := 1 + \frac{\sigma^2_R}{2}$, $q_n := n\!\left(1 - (n-1)\frac{\sigma^2_R}{2}\right)$, and $p_k := 1-(k-1)\sigma^2_R > 0$, for $k = 1,\ldots,n$. In particular, $p_n = q_n - q_{n-1}$, $p_1 = 1$, and $\sigma_\xi - \frac{\sigma_R^2}{2} = 1$.
\end{lem}

\begin{proof}
We first expand $Q_t^n$ as an $n$-fold integral over the simplex 

$$S_n(t) := \{0 \le s_1 \le \cdots \le s_n \le t\}$$ to obtain
\begin{align*}
    Q_t^n = n!\,\mu^n_R \int\limits_{S_n(t)} \prod_{j=1}^n e^{-(t-s_j)\sigma_{\xi}} e^{\sigma_R (W_t-W_{s_j})}\, d\mathbf{s},
\end{align*}
where $d\mathbf{s} = ds_1\cdots ds_n$.

\noindent Then, setting $s_{n+1} := t$, since $t - s_j = \sum\limits_{k=j}^{n}(s_{k+1}-s_k)$, interchanging the order of summation gives
\begin{align*}
\sum_{j=1}^n (t-s_j) = \sum_{k=1}^n k\,(s_{k+1}-s_k), \qquad \sum_{j=1}^n (W_t-W_{s_j}) = \sum_{k=1}^n k\,(W_{s_{k+1}}-W_{s_k}).
\end{align*}
Consequently,
\begin{align*}
\prod_{j=1}^n e^{-(t-s_j)\sigma_\xi + \sigma_R(W_t-W_{s_j})} = \prod_{k=1}^n e^{-k\sigma_\xi(s_{k+1}-s_k)}\,e^{k\sigma_R(W_{s_{k+1}}-W_{s_k})}.
\end{align*}

\noindent Since $\{W_{s_{k+1}}-W_{s_k}\}_{k=1}^n$ are independent with $W_{s_{k+1}}-W_{s_k} \sim \mathcal{N}(0,s_{k+1}-s_k)$, applying Fubini's theorem and Lemma~\ref{lem:mgf} gives
\begin{align*}
\E\left[Q_t^n\right] = n!\,\mu^n_R \int\limits_{S_n(t)} \prod_{k=1}^n e^{-q_k(s_{k+1}-s_k)}\, d\mathbf{s},
\end{align*}
where $q_k = k\!\left(\sigma_\xi - k\tfrac{\sigma_R^2}{2}\right) = k\!\left(1-(k-1)\tfrac{\sigma_R^2}{2}\right)$. Evaluating the simplex integral iteratively from the outermost variable inward and using $p_k = q_k - q_{k-1}$ at each step yields
\begin{align}
    \E\left[Q_t^n\right] = \mu^n_R \, e^{-q_n t} \prod_{k=1}^n \frac{e^{p_k t} - 1}{p_k}.
\end{align}
As a consequence, the second moment takes the explicit form
\begin{align*}
    \E[Q_t^2] &= \frac{\mu^2_R}{1 - \sigma^2_R} \left(1 - e^{-t} - e^{-t(1 - \sigma^2_R)} + e^{-t(2 - \sigma^2_R)} \right),
\end{align*}
and, for $\sigma_R^2 < \frac{1}{3}$, the fourth moment satisfies 
\begin{align*}
\E\left[Q_t^4\right] &= \mu^4_R \, e^{-q_4 t} \prod_{k=1}^4 \frac{e^{p_k t}-1}{p_k} \le 4(2\mu_R)^4.
\end{align*}
\end{proof}
\noindent
We can now state the almost-sure boundedness and H\"older-regularity of the process $(\xi_t)_{t \ge 0}$.
\begin{lem}[]
\label{lem:xi_asure_bded}
Let $m\in2\N$, $m > 2$, and fix a finite horizon $T>0$. Then all moments of $(\xi_t)_{t\in[0,T]}$ up to order $m$ are finite (with constants allowed to depend on $T$, $m$, $\sigma_R$, $\mu_R$, and $\xi_0$). In particular, there exists a finite constant $\const_{\xi}=\const_{\xi}(T,m)$ such that

$$ \sup{t\in[0,T]}|\xi_t(\om)| \le \const_{\xi} \quad \text{\rm for a.a.\ }\om\in\Om. $$
Moreover, the sample paths satisfy
\begin{align*}
\E\big[|\xi_t - \xi_s|^m\big] \le C(T,m,\sigma_R,\mu_R,\xi_0)\,|t-s|^{m/2}, \qquad s,t \in [0,T],
\end{align*}
and, consequently, $\Om \ni \om \mapsto \xi_t(\om) \in \cnt^{\vsig}([0,T];\R)$ for every $\vsig < \tfrac{1}{2}$.
\end{lem}
\begin{proof}
Writing $\xi_t = P_t + Q_t$, with $P_t$ and $Q_t$ as above, finiteness of moments of every order on the compact interval $[0,T]$ follows from the explicit representation~\eqref{eq:Exp_OU_proc_solved}: the factors $e^{\sigma_R W_t}$ have Gaussian moments of all orders, and the time integrals remain bounded on $[0,T]$. For each fixed $t$, Markov's inequality gives $|\xi_t|<\infty$ a.s.. The uniform bound $\sup{t\in[0,T]}|\xi_t(\om)|<\infty$ for a.a.\ $\om\in\Om$ is then deduced from~\eqref{eq:exp_ou_proc_phold} as follows.

\noindent
Setting $P_{t,r} := e^{-(t-r)\sigma_\xi + \sigma_R(W_t-W_r)}$, the explicit representation \eqref{eq:Exp_OU_proc_solved} gives $\xi_t = \xi_0 P_{t,0} + \mu_R\int_0^t P_{t,r}\,dr$, and hence 

$$\xi_t - \xi_s = \xi_0(P_{t,0}-P_{s,0}) + \mu_R\int_s^t P_{t,r}\,dr + \mu_R\int_0^s(P_{t,r}-P_{s,r})\,dr\,.$$ 
For each term, the factorization
\begin{align*}
P_{t,r} - P_{s,r} = e^{-s\sigma_\xi}e^{\sigma_R(W_s-W_r)}\big(e^{-\sigma_\xi(t-s)+\sigma_R(W_t-W_s)}-1\big)
\end{align*}
together with the elementary bounds $|e^x - 1| \le |x|e^{|x|}$ for the exponential differences and $\E[|W_t-W_s|^m]\le \tilde{C} |t-s|^{m/2}$, yields
\begin{align*}
\E\big[|P_{t,r}-P_{s,r}|^m\big] \le C\,|t-s|^{m/2}\,e^{C_m T}\,e^{-m(t-r)\max\Big(\bigl(1-(m-1)\sigma_R^2/2\bigr),0\Big)},
\end{align*}
where the factor $e^{C_m T}$ accounts for the (possibly growing) exponential when $1-(m-1)\sigma_R^2/2\le 0$. Integrating in $r\in[0,s]\subset[0,T]$ and applying Lemma~\ref{lem:mgf} therefore gives, for a constant $C=C(T,m,\sigma_R,\mu_R,\xi_0)$,
\begin{align}
\label{eq:exp_ou_proc_phold}
\E\big[|\xi_t - \xi_s|^m\big] \le C\,|t-s|^{m/2}.
\end{align}
Since $m/2 > 1$, the Kolmogorov--Chentsov theorem (see \cite{PrevotRockner2007}) provides a modification of $\xi$ with locally H\"older continuous paths of exponent $\vsig < \tfrac{1}{2} - \tfrac{1}{m}$. Sending $m\to\infty$ through even integers (admissible on each fixed $[0,T]$) and applying a standard Garsia--Rodemich--Rumsey estimate {(see, e.g., \cite{DaPrato1992})} to \eqref{eq:exp_ou_proc_phold} yields 

$$\E[\|\xi\|^m_{\cnt^\vsig}] < \infty \text{ for every $\vsig < \tfrac{1}{2}$.}$$ 
In particular, $\sup{t\in[0,T]}|\xi_t(\om)|\le \const_{\xi}$ for a.a.\ $\om\in\Om$ and some finite $\const_{\xi}$ depending on $T$.
\end{proof}
\noindent
As a consequence of Lemma \ref{lem:xi_asure_bded} we have the following result.
\begin{cor}[]
\label{cor:lp_hold_cont_proc}
Fix $T<\infty$. For every $m \in [1,\infty)$ and every $\vsig < {1\over2}$, the process $(\xi_t)_{t\in[0,T]}$ belongs to $L_m\big(\Om; \cnt^{\vsig}([0,T];\R)\big)$. In particular, under Assumption~\ref{asm:model_params}.\ref{Assumpt-3.4.1.} ($\sigma_R<1$), this holds for the H\"older exponents $\vsig<\min\{1/p,1/2\}$ used in Section~\ref{sec:analysis}.
\end{cor}
\begin{proof}
    By Lemma~\ref{lem:xi_asure_bded} we have $\sup{t\in[0,T]}|\xi_t(\om)|\le \const_{\xi}$ for a.a.\ $\om\in\Om$, hence
    
    $$\sup{t\in[0,T]} |\xi_t| \le \const_{\xi}, ~~\text{\rm almost surely.}$$
    Letting $N_{m} \subset \mathcal{F}$  such that $\P(N_{m}) = 0$ and $\sup{t\in[0,T]} |\xi_t(\om)| > \const_{\xi}$, we have that 
    \begin{align*} 
    \E\left[ (\sup{t\in[0,T]} |\xi_t|)^m \right] &= \E\left[ (\sup{t\in[0,T]} |\xi_t|)^m \indc_{N_{m}^c}\right] + \E\left[ (\sup{t\in[0,T]} |\xi_t|)^m \indc_{N_{m}}\right] \\[0.2cm]
    &= \E\left[ \const_{\xi}^m \right] \le \const_{\xi}^m.
    \end{align*}
    On the other hand, from Lemma~\ref{lem:xi_asure_bded} we also obtain that, for every $\vsig<\tfrac12$ (choosing an even integer $m>2/(1-2\vsig)$),
    \begin{align*}
        \sup{0 \le s\le t \le T} {|\xi_t - \xi_s| \over |t - s|^{\vsig} } \le C(T)  ~~\text{\rm almost surely.}
    \end{align*}
    Thus, letting $N_{\vsig} \subset \mathcal{F}$  such that $\P(N_{\vsig}) = 0$ and 
    
    $$\sup{t\in[0,T]} |\xi_t(\om)| > \const_{\xi} ~~\text{\rm or } \sup{0 \le s\le t \le T} {|\xi_t - \xi_s| \over |t - s|^{\vsig} } > C(T),$$ we have that 
    \begin{align*} 
    \E\left[ \|\xi\|_{\cnt^{\vsig}}^m \right] &= \E\left[ \|\xi\|_{\cnt^{\vsig}}^m \indc_{N^c_{\vsig}}\right] + \E\left[ \|\xi\|_{\cnt^{\vsig}}^m \indc_{N_{\vsig}}\right] \\[0.2cm]
    &= \E\left[ (\const_{\xi} + C(T))^m \right] \le (\const_{\xi} + C(T))^m.
    \end{align*}
\end{proof}

\subsection{Properties of {$\Tst$}}\label{sec:SG_Prop}

We now state (without proof) the properties of the evolution operator $\Tst$. 

\begin{thm}\label{thm:non_auto_evol_op} Let $A(t)$ be a family of sectorial operators on $\V_p$ satisfying the structural properties \eqref{eq:Au_struc_property}. Then, there exists a
unique evolution operator $\Tst$ satisfying 
\begin{align*}
    \Tst = \T_{r,t} \T_{s,r}, \quad 0 \le s \le r \le t \le T,\quad 
    \T_{s,s} = I, \ 0 \le s \le T.
\end{align*}
Furthermore, $A(t)\Tst$ is an $\L(\V_p)$-valued continuous function for $0 \le s \le t \le T$ with 

$$ \|A(t) \Tst\|_{\L(\V_p)} \le C|t-s|^{-1}, \quad 0 \le s \le t \le T, $$
$\Tst$ is an $\L(\V_p)$-valued strongly continuous function  with

$$ \|\Tst\|_{\L(\V_p)} \le C,\quad  0 \le s \le t \le T$$ 
and, for $t > s$, $\Tst$ is differentiable in $t$ with the derivative 
\begin{align*}
    \frac{\partial}{\partial t} \Tst &= -A(t) \Tst.
\end{align*}
\end{thm}
\begin{proof}
    Refer to Theorem 3.8, page 146 of \cite{Yagi09}.
\end{proof}

\noindent Let $\Tst$ be the evolution operator obtained by Theorem \ref{thm:non_auto_evol_op} (under the structural assumptions \eqref{eq:Au_struc_property}). Then the following properties hold true (cf. Section 8.1 page 154 \cite{Yagi09}).


\begin{enumerate}
    \item For any $\alpha \in[0, \vsig+ \vartheta)$,  then
        \begin{align}
        \label{eq:AthetaT_bdd}
        \begin{aligned}
            \mathcal{R}\left(\Tst\right) & \subset D\left(A_t^\alpha\right), \quad 0 \leq s<t \leq T\,, \\[0.2cm]
            \left\|A_t^\alpha \Tst\right\| & \leq C_\alpha |t-s|^{-\alpha}, \quad 0 \leq s<t \leq T.
        \end{aligned}
        \end{align}

    \item For $\alpha \in[0, \vsig],$ $\Tst A^{\alpha}_s$ can be extended as a bounded operator on $X$, i.e.,
        \begin{align}
        \label{eq:TAtheta_bdd}
            \|\Tst A^{\alpha}_s\| & \le C_{\alpha} |t-s|^{-\alpha}, \quad 0 \leq s<t \leq T.
        \end{align} 
    \item For $\varphi \in[0,1],$ $\alpha \in[\varphi, \vsig+\vartheta]$,  then
    \begin{align}
        \label{eq:Atheta_T_Aphi_theta_bdd}
        \begin{aligned}
        \|A^{\alpha}_t \Tst A^{\varphi - \alpha}_s\| &\le C_{\alpha} |t-s|^{-\varphi}, \quad 0 \leq s<t \leq T \\[0.2cm]
        \overset{\varphi = \alpha = 1}{\implies} \|A_t \Tst\| &\le C_{\vartheta} |t-s|^{-1}, \quad 0 \leq s<t \leq T.
        \end{aligned}
    \end{align}
    
    \item For $\varphi \in [0,1]$, $\alpha \in [0, \vsig +\vartheta)$, and $k=0,1$, then
        \begin{align}
        \label{eq:Atheta_T_minus_E_Aphi}
        \begin{aligned}
        \left\|A_t^\alpha\left[\Tst-e^{-(t-s) A_t}\right] A_s^{-\varphi}\right\| &\leq C_\vartheta ~ (t-s)^{\varphi-\alpha+\vsig+\vartheta-1}, \quad 0<\vartheta \leq t \leq T\,, \\[0.2cm]
        \left\|A_t^k \Tst A_s^{-k}-e^{-(t-\beta) A_s}\right\| &\leq C ~ (t-s)^{\vsig+\vartheta-1}.
        \end{aligned}
        \end{align}
    \item For $\alpha \in (0,1)$ and $0 \le s \le T$,  then
    \begin{align}
        \label{eq:Atheta_T_A_minus_E}
        \left\|A_t^{\alpha} \Tst A^{-\alpha} - e^{-(t-s)A_s} \right\| \le \begin{cases}
                        C ~ \Gamma(1-\vartheta) (t-s)^{\vsig+\vartheta-1}, &   0 \leq s<t \leq T,\  ~\vartheta < 1 \\[0.2cm]
                        C ~ \ln((t-s)^{-1} +1)|t-s|^{\vsig}, & 0 \leq s<t \leq T,\  ~\vartheta = 1.\\
                    \end{cases}
    \end{align}
\end{enumerate}

\subsection{Function spaces}\label{def-theo}

\noindent
In the following (unless otherwise specified) $X$ will denote a Banach space.

\begin{define}[Bessel potential space, see \cite{ADAMS99}] \label{def:BesPotSpace} The space $\wsp(\R ^n)$ with $s\in \R$ and $1<p<\infty $ is the function space defined by 
\begin{align}
\wsp(\R^n) &:= \Big\{ f \in  \TD(\R^n) :  \IFT\big[(1 + |\xi|^2)^{\frac{s}{2}} \FT f\big] \in L_p(\R^n) \Big\} \\[0.2cm]
\|f\|_{\wsp(\R^n)} &:= \left\| \IFT \big[(1 + |\xi|^2)^{\frac{s}{2}}\FT f \big] \right\|_{L_p(\R^n)}, 
\end{align}
where $\TD(\R^n)$ denotes the dual of the Schwartz space of rapidly decreasing functions

$$\SW(\R^n):=\{u\in \cnt^{\infty }(\R^n):\forall \ \alpha ,\ \beta \in \N^n\ \sup{\xbf \in R^n}|\xbf^{\alpha }D^{\beta }u(\xbf )|<\infty \}$$ 
and 
$\FT $ is the usual Fourier transform operator.\\[0.5ex] 
Alternatively (see \cite{GRISV85}), it can also be defined as the space of all distributions $g$ on $\R^n$ such that 

$$G_s\ast g\in L_p(\R^n),$$ where $G_s := \IFT\left[(1 + |\xi|^2)^{\frac{s}{2}}\right]$ denotes the Bessel potential of order $s$ and $\ast $ the convolution product.\\
\end{define}
\begin{define}
{\allowdisplaybreaks Let $\dom$ be a bounded Lipschitz domain in $\R^n$. $u \in \wsp(\dom)$  is an equivalence class of functions $U \in \wsp(\R^n)$ such that $U|_{\dom} = u$.
That is, a function $u$ is said to be in $\wsp(\dom)$ if and only if there exists a function $U \in \wsp(\R^n)$ such that $U|_{\dom} = u$. \\
Moreover, $\wsp(\dom)$ 
endowed with the norm \[ \|u\|_{\wsp(\dom)} := \underset{U \in \wsp(\R^n)}{\inf{ U|_{\dom} = u,}} \|U\|_{\wsp(\R^n)} \] is a Banach space.}
\end{define}
\begin{rem} [\cite{GRISV85}]\label{rem:equiSSBPS} Let $W^k_p(\dom),\ k\in \N,\ p > 1$, denote the usual Sobolev space of order $k$ and regularity $p$. Then, we have the following relations between the Bessel potential spaces $\W^s_p$ and $W^s_p$.
\begin{enumerate}
\item For $s \in \R$, $\W^s_2(\dom) := \W^s(\dom) = W^s_2(\dom)$  in norm equivalence.
\item For $m \in \N$ and $1 < p < \infty$, $\W^m_p(\dom) = W^m_p(\dom)$ in norm equivalence.
\item $W^s_p(\R^n) \subseteq \wsp(\R^n), \quad 1 < p \le 2$.
\item $\wsp(\R^n) \subseteq W^s_p(\R^n), \quad 2 \le p < \infty$.\\[0.1cm]
This implies that: for $\dom \subset \R^n$ a bounded Lipschitz domain, $\wsp(\dom) \subset W^s_p(\dom)$, $2 < p < \infty$. Indeed, since $W^s_p(\overline{\dom})$ is defined 
in a way analogous to the definition of $\wsp(\dom)$ and due to the inclusion $\W^s(\R^n) \subset W^s_p(\R^n)$, then $U \in \wsp(\R^n)$ with $U|_{\dom} = u$ 
implies $U \in W^s_p(\R^n)$ with $U|_{\dom} = u$. Moreover, since 
\[ \|u\|_{\wsp(\dom)} 
:= \underset{U \in \wsp(\R^n)}{\inf{ U|_{\dom} = u,}} \|U\|_{\wsp(\R^n)} \ge \underset{U \in W^s_p(\R^n)}{\inf{ U|_{\dom} = u,}} \|U\|_{W^s_p(\R^n)} 
=: \|u\|_{W^s_p(\overline{\dom})},
\] 
we get $\|u\|_{\wsp(\dom)} \ge \|u\|_{W^s_p(\dom)}$. The claim follows from the fact that $W^s_p(\dom) = W^s_p(\overline{\dom})$ for a bounded Lipschitz domain $\dom$.
\item $\W^{s_1}_p(\R^n) \subset W^{s_2}_p(\R^n) \subset \W^{s_3}_p(\R^n) \subset W^{s_4}_p(\R^n)$, for $ s_1 > s_2 > s_3 > s_4$.\\[1ex]
As in item 4. this implies that for a bounded Lipschitz domain 
\begin{align}
\label{eq:bded_sob_embd}
\W^{s_1}_p(\dom) \subset W^{s_2}_p(\dom) \subset \W^{s_3}_p(\dom) \subset W^{s_4}_p(\dom), \quad s_1 > s_2 > s_3 > s_4 .    
\end{align}
\end{enumerate}
\end{rem}

\subsection{Domains of fractional powers of elliptic operators in $L_p$ spaces}

Let $\dom$ be a bounded domain in $\mathbb{R}^d$ with $\mathcal{C}^2$-boundary $\partial \dom$. Consider the bilinear form
\begin{align*}
a(u, v)=\sum_{i, j=1}^n \int_{\dom} a_{i j}(x) D_i u D_j \bar{v} d x+\int_{\dom} c(x) u \bar{v} d x, \quad u, v \in \W^1(\dom)
\end{align*}
on the space $\W^1(\dom)$. Let $a_{i j}(x),\, 1 \leq i, j \leq n$, be real-valued functions in $\dom$ satisfying the conditions
\begin{align*}
\begin{aligned}
& a_{i j} \in L_{\infty}(\dom), \quad 1 \leq i, j \leq n\,, \\
& \sum_{i, j=1}^n a_{i j}(x) \xi_i \xi_j \geq \lbd^{a} |\xi|^2, \quad \xi=\left(\xi_1, \ldots, \xi_n\right) \in \mathbb{R}^d, \text { a.e. } x \in \dom
\end{aligned}
\end{align*}
with some constant $\lbd^{a} > 0$, and $c(x)$ be a real-valued function in $\dom$ satisfying
\begin{align*}
c \in L_{\infty}(\dom) \text { and } c(x) \geq c_0>0 \text {, a.e. } x \in \dom .
\end{align*}

\noindent The operator $\mathcal{A}$ (resp. $A$ ) is regarded as a realization of the strongly elliptic operator 

$$ A := -\sum_{i, j=1}^n D_j\left[a_{i j}(x) D_i (\cdot)\right]+c(x)$$ generated by the bilinear form $a(\cdot,\cdot)$ in $\W^1(\dom)^*\left(\right.$resp. $L_2(\dom)$) under the Neumann-type boundary conditions
\begin{align*}
\frac{\partial u}{\partial {\nu}_A} \equiv \sum_{i, j=1}^n {\nu}_j(x) a_{i j}(x) D_i u=0 \quad \text { on } \partial \dom.
\end{align*}

\noindent We consider a realization $A_p$ of the strongly elliptic operator $A$ in the space $L_p(\dom),\  1<p<\infty$, under the Neumann-type boundary conditions (16.27) in \cite{Yagi09}. We want to characterize the domains of the fractional powers $A_p^\theta$ for $0 \leq \theta \leq 1$.

\noindent The domain of $A_p$ is given by
\begin{align*}
\mathcal{D}\left(A_p\right)=\left\{u \in \W_p^2(\dom) :\frac{\partial u}{\partial {\nu} _A}=0 \text { on } \partial \dom\right\}
\end{align*}
with the estimate
\begin{align*}
\tilde{\lbd}^a_p \|u\|_{\W_p^2} \leq\left\|A_p u\right\|_{L_p}, \quad u \in \mathcal{D}\left(A_p\right)
\end{align*}
$\tilde{\lbd}^a_p > 0$ being a positive constant. 
It also holds that
\begin{align}
\label{eq:dom_embds}
\begin{cases}\W_p^{2 \theta^{\prime \prime}}(\dom) \subset \mathcal{D}\left(A_p^\theta\right) \subset \W_p^{2 \theta^{\prime}}(\dom) & \text { if } 0 \leq \theta^{\prime}<\theta<\theta^{\prime \prime}<\frac{p+1}{2 p}, \\[0.3cm] 
\W_{p, N}^{2 \theta^{\prime \prime}}(\dom) \subset \mathcal{D}\left(A_p^\theta\right) \subset \W_{p, N}^{2 \theta^{\prime}}(\dom) & \text { if } \frac{p+1}{2 p}<\theta^{\prime}<\theta<\theta^{\prime \prime} \leq 1,\end{cases}
\end{align}
with continuous embeddings, where
\begin{align*}
    \W_{p, N}^{s} =\left\{u \in \W_p^s(\dom) :\frac{\partial u}{\partial {\nu}_A}=0 \text { on } \partial \dom\right\} \text{ for $s \in \left(1 + \frac{1}{p}, 2 \right]$ }.
\end{align*}

\begin{define}[Uniformly bounded function space, see Sec. 1.2.1 in \cite{Yagi09}]\label{def:bdBochSp} 
Let $[a, b]$ be a closed interval of $\R$ and $X$ a Banach space, 
then by $\B_{_X}:=\B([a, b];X)$ we denote the space of uniformly bounded functions on $[a, b]$ 
with values in $X$. 
\end{define}

\begin{define}[H\"older continuous function space, see Sec. 1.2.3 in \cite{Yagi09}] \label{def:timeHoldSp} 
Let $[a, b]$ be a closed interval of $\R$ and $X$ a Banach space, then for $0< \mu <1$,\, $\cnt^{\mu}_{\{a\}}([a, b];X)$ denotes the space of $X$-valued functions that are continuous on the $[a,b]$ and are $\mu$-H\"older continuous at least at $a$. It is endowed with the norm 
\[\|f\|_{\cnt^{\mu}_{\{a\}}} := \max_{t \in [a,b]} \|f\|_{_X} + \sup{t \in [a,b]} {\|f(t) - f(a)\|_{_X} \over (t - a)^{\mu}}. \]
Similarly, for $0< \mu <1$, we denote by $\cnt^{\mu}([a, b];X)$ the space of $X$-valued functions that are $\mu$-H\"older continuous on $[a,b]$. It is endowed with the norm
\[\|f\|_{\cnt^{\mu}} := \max_{t \in [a,b]} \|f\|_{_X} + \sup{s,t \in [a,b], s \ne t} {\|f(t) - f(s)\|_{_X} \over (t - s)^{\mu}}. \]
\end{define}

\begin{define}[Weighted H\"older continuous functions, see Sec. 1.2.4 in \cite{Yagi09}]\label{def:weighHol} The space $\G^{\eta,\rho}\big((a, b]; X\big)$ with 
$0 < \rho < \eta \le 1$ consists of $X-$valued functions on $(a, b]$ (resp. $[a, b]$) when $\eta < 1$ (resp. $\eta = 1$) with the following properties.\\[2ex]
1. When $\eta < 1$, the function $(t - a)^{1 - \eta} f(t)$ has a limit for $t \to a$.\\[1ex]
2. $f$ is a H\"older continuous function with the weight $(s - a)^{1 - \eta + \rho}$ and with exponent $\rho$, i.e.,
\[ \sup{a \le s < t \le b} \frac{ (s - a)^{1 - \eta + \rho} \|f(t) - f(s)\|_{X}}{|t - s|^{\rho}} 
= \sup{a < t \le b} \: \sup{a \le s < t} \frac{ (s - a)^{1 - \eta + \rho} \|f(t) - f(s)\|_{X}}{|t - s|^{\rho}} < \infty . \]
3. As $t \to a$, it holds that $w_f(t) \to 0$, where 
\begin{align} 
w_f(t) := \sup{a \le s < t}  \frac{ (s - a)^{1 - \eta + \rho} \|f(t) - f(s)\|_{X}}{|t - s|^{\rho}}\,. \label{eq:weighHol_cont}
\end{align}
4. The space $\G^{\eta, \rho}$ equipped with the norm \[ \|f\|_{\G^{\eta,\rho}} := \sup{a \le t \le b} (t - a)^{1-\eta}\|f(t)\|_X + \sup{a < t \le b} \: \sup{a \le s < t}  
\frac{ (s - a)^{1 - \eta + \rho} \|f(t) - f(s)\|_{X}}{|t - s|^{\rho}} \] is a Banach space.\\[1ex]
\end{define}

\begin{rem}$ $\\[0.1cm]
1. If $0 < \rho' < \rho < \eta \le 1$, then $\G^{\eta, \rho} \subset \G^{\eta,\rho'}$.\\[0.2cm]
2. If $\eta' < \eta$, then $\G^{\eta,\rho} \subset \G^{\eta',\rho}$, with $0< \rho < \eta' < \eta$.\\[0.2cm]
3. Let $g \in \cnt^{\rho}([0,b];X)$, with $g(0) = 0$, then for $0 < \rho < \eta < 1$ the function $f(t) := t^{\eta-1} g(t)$ belongs to $\G^{\eta, \rho}((0,b];X)$. 
\end{rem}

\begin{thm}[Theorem 3.9 and Theorem 3.10 from \cite{Yagi09}] \label{thm:ExiLinNonAuto}
Consider the Cauchy problem 
\begin{align}
\label{eq:LinNonAuto}
\left. 
\begin{array}{l}
{d \over dt} u + A(t) u = f(t) 	\quad ~~\text{\rm in } X, \quad t \in (0,T]\\[0.3cm]
u(0) = u_0\,.\\
\end{array} \right\}
\end{align}
Let $A(t)$ be a sectorial operator with a uniform spectral angle $\kappa < { \pi \over 2}$ and a uniform resolvent 
estimate \[ \|R_{\lambda}(A(t)\| \le {M \over |\lambda|}, \quad \kappa \notin \overline \Sigma_{\kappa},\ M \ge 1, t \in [0,T].\] Moreover, the 
domain of $A(t)$ may vary in time, but must satisfy the following conditions:
\begin{align}
\left.
\begin{array}{l}
 \D(A(t)) \subset \D(A^{\nu}(s)),\ s,t \in [0,T], \quad \nu \in (0,1] \\[0.3cm]
\|A^{\nu}(t) ( A^{-1}(t) - A^{-1}(s) ) \| \le C |t - s|^{\mu}, \quad \mu \in (0,1] \\[0.3cm]
 1 < \mu + \nu\,.\\[1ex]
\end{array} \right\}
\end{align}
Also, let $f \in \Fbs((0,T];X)$, with $\sigma < \min(\beta, \mu + \nu - 1)$.\\[0.5ex] 
Then, for every $u_0 \in \D(A(0)^{\beta})$ there exists a unique solution $u$ of \eqref{eq:LinNonAuto} given by
\[ u(t) = U(t,0) u_0 + \int_0^t U(t,s) f(s), \quad t \in [0,T], \]
where $U(t,s)$ denotes the evolution operator as per Theorem 3.8 in \cite{Yagi09}.
Also, $u$ satisfies the following estimates:
\[ u \in \cnt([0,T];X) \cap \cnt^1((0,T];X), \quad A(t)^{\beta} u \in \cnt([0,T];X),\quad  A(t) u \in \G^{\beta, \sigma}((0,T];X),\] with 
\[ \|u(t)\|_{\B_{_X}} + \|A(t)^{\beta} u\|_{\cnt_{_X}} \le \const ( \|A(0)^{\beta}u_0\|_X + \|f\|_{\Fbs}). \]
\[ \left\|{d \over dt} u(t) \right\|_{\Fbs} + \|A(t) u\|_{\Fbs} \le \const ( \|A(0)^{\beta} u_0\|_X + \|f\|_{\Fbs}), \]
where $\const $ denotes some positive constant. 
\end{thm}


\phantomsection
\printbibliography
\end{document}